\documentclass[11pt,a4paper]{article}
\usepackage[T1]{fontenc}
\usepackage[utf8]{inputenc}
\usepackage{lmodern,microtype}
\usepackage{amsmath,amssymb,amsthm,mathtools}
\usepackage[a4paper,margin=25mm]{geometry}
\usepackage{booktabs,array,graphicx,xcolor,enumitem}
\usepackage[colorlinks=true,linkcolor=blue!40!black,citecolor=blue!40!black,urlcolor=blue!40!black]{hyperref}
\ifdefined\pdfmapfile
\pdfmapfile{+lm.map}
\pdfmapfile{+cm.map}
\pdfmapfile{+cmextra.map}
\pdfmapfile{+symbols.map}
\pdfmapfile{+euler.map}
\fi
\hypersetup{pdftitle={Counting tree-child and reticulation visible networks},
 pdfsubject={Critical scaling, Airy constants and reticulation-visible defect asymptotics}}
\setlist[enumerate]{topsep=4pt,itemsep=2pt}
\setlist[itemize]{topsep=4pt,itemsep=2pt}
\numberwithin{equation}{section}
\newtheorem{theorem}{Theorem}[section]
\newtheorem{lemma}[theorem]{Lemma}
\newtheorem{proposition}[theorem]{Proposition}
\newtheorem{corollary}[theorem]{Corollary}
\theoremstyle{remark}

\newcommand{\Ai}{\operatorname{Ai}}

\newcommand{\RV}{\operatorname{RV}}

\title{Counting tree-child and reticulation visible networks}
\author{Miquel Pons Viver}
\date{24 November 2026}

\begin{document}
\maketitle

\begin{abstract}
For tree-child networks with $n$ leaves
and $k$ reticulations, we identify an entire critical scaling function
at $k=x\sqrt n$ and its first uniform correction. The function is a
Gaussian factor times the Laplace transform of Brownian excursion
area. An exact coefficient integral then yields the explicit amplitude
$\pi\Ai'(a_1)$ for the maximal-reticulation count, where $a_1$
is the largest negative Airy zero, together with three logarithmic
corrections. A summable defect expansion transfers these results to the
total number of tree-child networks.

For reticulation-visible networks, we construct an exact weighted
refinement over maximal tree-child networks. The total defect from
$3n-3$ reticulations splits into a skeleton defect, an arity excess,
and an ordinary-leaf count. The skeleton defect has a uniformly bounded
exponential moment, and the full partition function differs from its
maximal-skeleton sector by a bounded factor, smaller than $30$ at
activity one. A positive algebraic kernel and an explicitly solved
variational problem give the pressure with logarithmic error
$O(n^{2/3}\log^{2/3}n)$. This determines the total exponential growth,
the location of every mode, a full large-deviation principle, and
limiting structural densities. At activity $z=t/n$ with fixed $t>0$,
we also prove a $\operatorname{Poisson}(19t/4)$ limit and coefficient
asymptotics uniform for defects $q=o(n^{1/3})$.
\end{abstract}


\tableofcontents

\section{Introduction}
Counting phylogenetic networks involves both local restrictions on
reticulations and global constraints such as acyclicity. Even for
classes with a linear upper bound on the number of reticulations,
the asymptotic behaviour depends strongly on how that number grows
with the number of leaves. This paper studies two complementary
regimes and connects them through component decompositions.

For tree-child networks, an exact recurrence now follows from the
Pons--Batle enumerative identity, proved by Lin, Liu, Liu, Liu and Xin
using the enumeration of Chang, Fuchs, Liu, Wallner and Yu
\cite{PB21,Lin26,CFLLWY24}. The maximal case was already related to
words by Fuchs, Yu and Zhang, who established its Airy stretched
exponential up to a multiplicative constant \cite{FYZ21}. The leading
relation between the total count and the maximal count, and the
Poisson limit for the missing reticulations, were obtained in
\cite{CFLLWY24}. Recent work also establishes the optimal range
$k=o(\sqrt n)$ of the usual sparse-reticulation approximation
\cite{YZ26,FY26}. We use these results as enumerative starting points
and identify their analytic refinements.

For reticulation-visible networks, the component construction of
Chang and Fuchs reduces counting to local one-component weights on
tree-child component graphs \cite{CF24}. Its extremal case gives
$\RV_{n,3n-3}=T_{n,n-1}$. Away from that extremal case, however, a
component graph may have vertices of larger outdegree, and different
graphs support different numbers of expansions. An analysis of
skeletons alone therefore does not determine the network count. Our
refinement identity keeps the local multiplicities exactly and makes
it possible to control every skeleton sector under the weighted
counting measure.

\subsection{The principal results}
Write $T_{n,k}$ for the tree-child count and $M_n=T_{n,n-1}$.
Our critical-scale theorem identifies the function
\[
 F(x)=e^{-x^2/3}\mathbb E e^{-\sqrt2 x\mathcal A}.
\]
Here \(\mathbb E\) denotes expectation, and \(\mathcal A\) is the
area under a standard Brownian excursion, and proves
\[
 \frac{T_{n,k}}{\binom nk(2n+2k-3)!!}
 =F(x)-\frac{xF''(x)+x^2F'(x)}{2\sqrt n}+O_X(n^{-1}),
 \qquad x=\frac{k}{\sqrt n}\in[0,X].
\]
The exact coefficient solution and the Brownian-area identity are
proved in Section~\ref{sec:tc-critical}. This identifies the transition
function at the already known square-root threshold.

Let $a_1<0$ be the largest zero of $\Ai$, and set
\[
 D_n=(3n)^{-2/3}e^{a_1(3n)^{1/3}}(12/e^2)^n n^{2n},\qquad
 A=\frac{a_1^2}{6},\quad
 \beta=-a_1,\quad
 K=\pi\Ai'(a_1).
\]
The maximal and total counts satisfy
\begin{align*}
 M_n&=KD_n\exp\left\{A(3n)^{-1/3}+\beta(3n)^{-2/3}
                   +\frac{13}{6}(3n)^{-1}+O(n^{-4/3})\right\},\\
 T_n:=\sum_{k=0}^{n-1}T_{n,k}
 &=\sqrt e\,KD_n\exp\left\{A(3n)^{-1/3}
       +\frac{23\beta}{24}(3n)^{-2/3}
       +\frac{209}{96}(3n)^{-1}+O(n^{-4/3})\right\}.
\end{align*}
Sections~\ref{sec:tc-maximum} and \ref{tcm:tot-section-total} prove
these expansions. The Airy derivative appears because the nominal
first term of a cubic saddle vanishes at $a_1$.

For reticulation-visible networks put
\[
 V_{n,q}=\RV_{n,3n-3-q},\qquad
 Z_n(z)=M_n^{-1}\sum_{q=0}^{3n-3}V_{n,q}z^q.
\]
The exact refinement in Section~\ref{rvc:refinementsection} represents
this polynomial over maximal tree-child networks. If $Z_n^{(0)}$
denotes the sector with a maximal normalized component skeleton,
then Section~\ref{rvc:comparisonsection} proves
\[
 1\le Z_n(z)/Z_n^{(0)}(z)\le C(z)<\infty\quad(z>0),
 \qquad Z_n(1)<30 Z_n^{(0)}(1),
\]
uniformly in $n$. The skeleton defect has a uniform exponential
tail. Thus the other skeleton sectors affect the count by a bounded
factor while leaving the pressure unchanged.

The pressure is determined by the kernel
\[
 \mathcal C(x;z)=\frac{1-2zx}{(1-2zx)^{3/2}-x}.
\]
Let $\rho(z)$ be its positive pole, and let
$\xi(z)\mathcal C_x(\xi(z);z)=\mathcal C(\xi(z);z)$. Define
\[
 \Lambda(z)^{-1}=\frac12\left\{
 \int_{\xi(z)}^{\rho(z)}\frac{dx}{\mathcal C_x(x;z)}
 +\frac{\xi(z)}{\mathcal C_x(\xi(z);z)}\right\},
 \qquad \psi(z)=\log(\Lambda(z)/12).
\]
Sections~\ref{rva:laplace} and \ref{rva:pressuresection} prove
\[
 \log Z_n(z)=n\psi(z)
       +O_{\mathcal K}(n^{2/3}\log^{2/3}(n+2))
\]
uniformly for $z$ in positive compact sets. Consequently,
\[
 \log\RV_n=2n\log n+n\log(\Lambda(1)/e^2)
                 +O(n^{2/3}\log^{2/3}(n+2)),
\]
with $\Lambda(1)/e^2=13.9657067341\ldots$.
Every mode $q_n$ of $(V_{n,q})_q$ satisfies
$q_n/n\to\psi'(1)=1.1969529048\ldots$.
The same pressure gives a full large-deviation principle, an entropy
formula for counts in linear windows, and a leaf-marked extension
whose ordinary-leaf density is $0.2182255858\ldots$.

At a different scale, the $5/8$ law for leaves of maximal tree-child
networks in Section~\ref{sec:leaf-law} yields
\[
 Z_n(t/n)=e^{19t/4}+O_R(n^{-1/3})\quad(|t|\le R),\qquad
 V_{n,q}=M_n\frac{((19/4)n)^q}{q!}(1+o(1))
\]
uniformly for $q=o(n^{1/3})$. In particular, the leading coefficient
at defect two is $361/32$. The macroscopic defect distribution and
the microscopic Poisson limit arise from different activity scales.

\subsection{Conventions and scope}\label{sec:conventions}
All original networks in this paper are finite, simple directed
acyclic graphs. Their root has indegree zero and outdegree one;
leaves have indegree one and outdegree zero and are bijectively
labelled by $[n]=\{1,\ldots,n\}$. An internal tree vertex has degree
$(1,2)$, and a reticulation has degree $(2,1)$. There are no other
vertices. Networks are identified up to directed graph isomorphism
preserving the leaf labels. Their children are unordered.

A network is \emph{tree-child} if every non-leaf vertex has a child
which is not a reticulation. It is \emph{reticulation-visible} if,
for every reticulation, there is a leaf such that every path from
the root to that leaf passes through that reticulation. The sets
and counts are denoted by $\mathcal{TC}_{n,k}$, $T_{n,k}$ and
$\RV_{n,r}$ as above. Tree-child networks have at most $n-1$
reticulations. For example, after deleting incoming reticulation
arcs, every resulting tree component contains a distinct original
leaf, which gives this bound. The component argument in
Section~\ref{rvc:section} gives the reticulation-visible upper bound
$3n-3$ and its equality case.

The intermediate skeletons may have tree vertices of outdegree
larger than two. We call such a skeleton \emph{semibinary}: its
reticulations still have degree $(2,1)$, and each tree vertex has
indegree one and outdegree at least two. This change in degree
convention applies only to the explicitly defined intermediate
skeletons; all enumerated networks remain binary.

The proofs use the published counting recurrence and component
reconstruction, and the established encoding of maximal tree-child
networks. The new arguments are supplied in full. The logarithmic
RV estimate is not a multiplicative equivalent, and convergence of
the mode locations does not assert finite-$n$ unimodality. Numerical evaluations of the pressure are distinguished from the
rational interval certification of the bound $30$.

% BEGIN SECTION SOURCE: tc_critical.tex
\section{The critical scale for tree-child networks}\label{sec:tc-critical}

\subsection{An exact recurrence and its coefficient solution}
Let $T_{n,k}$ count the networks specified in Section~\ref{sec:conventions}.
The enumerative identity of Pons--Batle, established by
Lin, Liu, Liu, Liu and Xin together with the network enumeration of
Chang, Fuchs, Liu, Wallner and Yu, implies the recurrence
\begin{equation}\label{tc:countrec}
 T_{n,k}=(n-k+1)T_{n,k-1}
       +\frac{n(2n+k-3)}{n-k}T_{n-1,k},\qquad 1\le k<n.
\end{equation}
Here $T_{n,0}=(2n-3)!!$, $T_{n,n}=0$, and $T_{1,0}=1$.
For clarity, the exact identity used is
$T_{n,k}=n!c_{n-1,k}/(n-k)!$, where
$c_{a,b}=c_{a,b-1}+(2a+b-1)c_{a-1,b}$; see
\cite[Theorem~1.1 and Section~2.6]{Lin26} and
\cite[Theorem~3.4]{CFLLWY24}.
Only this counting identity is needed in this section.

Normalize the counts by
\begin{equation}\label{tc:normalization}
 R_{n,k}=\frac{T_{n,k}}{\binom nk(2n+2k-3)!!}.
\end{equation}
Then \eqref{tc:countrec} becomes
\begin{equation}\label{tc:harmonic}
 R_{n,k}=\frac{kR_{n,k-1}+(2n+k-3)R_{n-1,k}}{2(n+k)-3},
 \qquad R_{n,0}=1,\quad R_{n,n}=0.
\end{equation}
The following representation is a convenient form of the coefficient
formula in \cite{PB21}. We give a direct verification from the recurrence.

\begin{proposition}[Exact coefficient representation]\label{tc:exact}
Define $f_0=1$ and, for $i\ge1$, determine $f_i$ by
\begin{equation}\label{tc:gamma-rec}
 \sum_{j=0}^i\frac{\Gamma((3i+j-1)/2)}{j!}f_{i-j}=0.
\end{equation}
For $0\le j\le k$ put
\begin{equation}\label{tc:H}
 H_j(n,k)=(k)_{\underline j}
   \frac{\Gamma(n+k-(j+1)/2)}{\Gamma(n+k-1/2)},
 \qquad H_j(n,k)=0\quad(j>k).
\end{equation}
Then, throughout $0\le k\le n$,
\begin{equation}\label{tc:exact-sum}
 R_{n,k}=\sum_{j=0}^k f_jH_j(n,k).
\end{equation}
\end{proposition}
\begin{proof}
The coefficient of $f_i$ in \eqref{tc:gamma-rec} is
$\Gamma((3i-1)/2)>0$, so the coefficients are uniquely determined.
Each $H_j$ satisfies \eqref{tc:harmonic} in its interior. Indeed, with
$t=n+k-1/2$ and $a=j/2$, the two ratios to $H_j(n,k)$ are
\[
 \frac{H_j(n,k-1)}{H_j(n,k)}
  =\frac{k-j}{k}\frac{t-1}{t-a-1},\qquad
 \frac{H_j(n-1,k)}{H_j(n,k)}=\frac{t-1}{t-a-1}.
\]
The asserted relation follows immediately; cases where a falling
factorial vanishes follow by substitution. The sum has boundary value
$1$ at $k=0$. At $k=n$, multiply the sum by
$\Gamma(2n-1/2)/n!$ and put $j=n-\ell$; the result is precisely
\eqref{tc:gamma-rec} with $i=n$. The finite triangular recurrence and
these boundary values have a unique solution, proving the proposition.
\end{proof}

For comparison with the original notation of \cite{PB21},
\[
 f_{2r}=\frac{b_{2r}}{2^r(2r)!},\qquad
 f_{2r+1}=\frac{\sqrt\pi\,b_{2r+1}}{2^{r+1}(2r+1)!},
\]
where $b_0=1$ and
$b_i=-\sum_{j=1}^i\binom ij(3i-3+j)!!\,b_{i-j}/(3i-3)!!$.
Neither this notation nor a correspondence between individual networks
and those words is required for the proof above.

\subsection{The Brownian-area identity}
Let $e(t)$, $0\le t\le1$, be a standard Brownian excursion, with Brownian
variance parameter one, and write $\mathcal A=\int_0^1e(t)\,dt$.
Let $a_1>a_2>\cdots$ be the negative zeros of the Airy function $\Ai$.

\begin{theorem}[An entire Airy scaling function]\label{tc:airy}
The coefficients in \eqref{tc:gamma-rec} have the entire generating
function
\begin{equation}\label{tc:F}
 F(z)=\sum_{j\ge0}f_jz^j
     =e^{-z^2/3}\mathbb E e^{-\sqrt2 z\mathcal A}.
\end{equation}
For $|\arg z|<3\pi/4$, with the principal power,
\begin{equation}\label{tc:spectral}
 F(z)=2\sqrt\pi\,z e^{-z^2/3}
         \sum_{\nu\ge1}e^{a_\nu z^{2/3}}.
\end{equation}
In addition, for $j\ge1$,
\begin{equation}\label{tc:coefficient-bound}
 |f_j|\le 2\sqrt2\left(\frac{5e}{3j}\right)^{j/2}.
\end{equation}
\end{theorem}
\begin{proof}
We first check that the right-hand side of \eqref{tc:F} is entire.
The excursion has the law of the norm of a three-dimensional standard
Brownian bridge $B$; see \cite{Janson07}. Cauchy--Schwarz and Jensen give
\begin{equation}\label{tc:gaussian-area}
 \mathbb E e^{\mathcal A^2}
 \le\int_0^1\mathbb E e^{\|B(t)\|^2}\,dt\le2\sqrt2,
\end{equation}
because each coordinate has variance $t(1-t)\le1/4$.
Consequently $|F(z)|\le2\sqrt2 e^{5|z|^2/6}$.
Cauchy's coefficient inequality on the circle of radius $\sqrt{3j/5}$
proves \eqref{tc:coefficient-bound} for its Taylor coefficients.
We must identify those coefficients with \eqref{tc:gamma-rec}.

For real $z\ge0$, set
\[
 \mathcal J(z)=\int_0^\infty\frac{e^{-u^2}}{u^2}
       \bigl(e^{zu^4}F(zu^3)-1\bigr)\,du.
\]
We show that $\mathcal J$ is smooth from the right at zero and flat
there. For $x>0$, let
\[
 X_t=x+W_t-\sqrt2 t+\frac z{\sqrt2}t^2,
 \qquad q_z(x)=\mathbb P(X_t>0\text{ for every }t\ge0).
\]
The first-hitting density at zero is
\begin{equation}\label{tc:hitting-density}
 \frac{x}{\sqrt{2\pi}t^{3/2}}
 e^{-x^2/(2t)+\sqrt2x-t+zt^2-z^2t^3/3}
 \mathbb E\exp\left\{-\sqrt2z\int_0^tY_s^{x,t}\,ds\right\},
\end{equation}
where $Y^{x,t}$ is a three-dimensional Bessel bridge from $x$ to zero.
Here is the change-of-measure calculation. Conditional on the first
hitting time of zero being $t$, a Brownian path from $x$ is this Bessel
bridge. With drift $b(s)=-\sqrt2+\sqrt2zs$, integration by parts in the
Girsanov exponent gives
\[
 \int_0^t b(s)\,dY_s
       -\frac12\int_0^t b(s)^2\,ds
 =\sqrt2x-\sqrt2z\int_0^tY_s\,ds-t+zt^2-z^2t^3/3.
\]
Multiplication by the drift-free first-hitting density yields
\eqref{tc:hitting-density}. The change of measure is on each finite
time interval; the total hitting mass is then obtained by monotone
integration over $t$.

As $x\downarrow0$, the bridge area converges in distribution to
$t^{3/2}\mathcal A$. Subtract the hitting density for $z=0$ before
dividing by $x$ and passing to the limit. This subtraction is essential
at $t=0$. The bridge can be realized as the norm of a three-dimensional
Gaussian bridge with mean path $(1-s/t)(x,0,0)$, so its expected area
is at most $xt/2+C t^{3/2}$. Thus the difference, after division by $x$,
is dominated near zero by $C_z(1+t^{-1/2})$ for $0<x\le1$.
At infinity it is dominated by $C_z t^{-3/2}e^{-t/4}$, using
\[
 1-v+v^2/3\ge1/4\qquad(v\ge0).
\]
Since $q_0(x)=0$, dominated convergence and $t=u^2$ give
\begin{equation}\label{tc:flat-identity}
 \mathcal J(z)=-\sqrt{\pi/2}\lim_{x\downarrow0}\frac{q_z(x)}x.
\end{equation}
For $z>0$ the drift of $X$ is at most $-1/\sqrt2$ until
$T=(2z)^{-1}$. Pathwise comparison with this constant drift and its
first-hitting density give
\begin{equation}\label{tc:flat-bound}
 0\le-\mathcal J(z)
 \le\frac12\int_{1/(2z)}^\infty t^{-3/2}e^{-t/4}\,dt
 \le2(2z)^{3/2}e^{-1/(8z)}.
\end{equation}
To justify differentiation at zero, write the integrand after $t=u^2$
as the difference between $e^{-t}$ and
\[
 e^{-t+zt^2-z^2t^3/3}
       \mathbb E e^{-\sqrt2zt^{3/2}\mathcal A},
\]
times $t^{-3/2}/2$. Every positive-order derivative is bounded, uniformly
for $0\le z\le z_0$, by a polynomial in $t$ and excursion-area moments
times $e^{-t/4}$. At zero the derivative introduces at least a factor
$t^{3/2}$, so the resulting bound is also integrable there.
All moments exist by \eqref{tc:gaussian-area}. This proves smoothness
from the right. Equation~\eqref{tc:flat-bound} makes all the derivatives
zero. If $F(z)=\sum \widehat f_jz^j$, the coefficient of $z^i$ in
$\mathcal J$ is therefore
\[
 0=\frac12\sum_{j=0}^i
       \frac{\Gamma((3i+j-1)/2)}{j!}\widehat f_{i-j}.
\]
Uniqueness of \eqref{tc:gamma-rec} proves $\widehat f_j=f_j$.

Finally, the classical excursion-area transform
\cite[Equation~(80)]{Janson07}, with Laplace argument $\sqrt2z$,
gives \eqref{tc:spectral} for $z>0$. The series converges normally in
$|\arg z|<3\pi/4$ since $|a_\nu|$ grows as a positive constant times
$\nu^{2/3}$. Analytic continuation proves the stated complex identity.
\end{proof}

\subsection{A uniform transition and its first correction}
\begin{theorem}[Critical expansion]\label{tc:critical}
For every fixed $X<\infty$, uniformly over integers
$0\le k\le X\sqrt n$, and with $x=k/\sqrt n$,
\begin{equation}\label{tc:critical-expansion}
 R_{n,k}=F(x)-\frac{xF''(x)+x^2F'(x)}{2\sqrt n}+O_X(n^{-1}).
\end{equation}
The assertion also holds as a relative expansion on each such compact
interval because $F(x)>0$ for $x\ge0$.
\end{theorem}
\begin{proof}
For $j\ge2$ and $j\le k$, the gamma distribution and Jensen's inequality
give $\Gamma(t+j/2)/\Gamma(t)\ge t^{j/2}$, where
$t=n+k-(j+1)/2\ge n$. Thus $0\le H_j(n,k)\le x^j$.
For $j=1$, the elementary half-step gamma bound gives
$H_1(n,k)\le Cx$ with an absolute $C$.
Together with \eqref{tc:coefficient-bound}, these bounds show that
the part $j>C_0\log n$ of \eqref{tc:exact-sum}, or of any fixed
derivative series of $F$, is smaller than every prescribed power
of $n^{-1}$. Any fixed $C_0>0$ suffices for all sufficiently large $n$.

The case $j=0$ is exact, since $H_0=1$. For
$1\le j\le C_0\log n$, expand the falling factorial as a polynomial
in $x$:
\[
 n^{-j/2}(k)_{\underline j}
 =x^j-\frac{j(j-1)}{2\sqrt n}x^{j-1}
   +O\left(\frac{j^4(1+X)^j}{n}\right).
\]
The term with factor $j(j-1)$ is interpreted as zero when $j=1$.
This is uniform even when $k<j$: both the original product and its
polynomial expression then vanish. To bound the remainder, the
elementary symmetric polynomial of degree $\ell$ in
$0,1,\ldots,j-1$ is at most $(j(j-1)/2)^\ell/\ell!$.
Summing terms of degree at least two proves the bound, since
$j^2/\sqrt n=o(1)$.
The integral of the digamma function, or Stirling's formula with its
bounded remainder, gives in the same range
\[
 n^{j/2}\frac{\Gamma(n+k-(j+1)/2)}{\Gamma(n+k-1/2)}
 =1-\frac{jk}{2n}+O_X\left(\frac{(1+j)^2}{n}\right).
\]
Multiplying yields
\begin{equation}\label{tc:H-expansion}
 H_j(n,k)=x^j-
 \frac{j(j-1)x^{j-1}+jx^{j+1}}{2\sqrt n}
 +O_X\left(\frac{(1+j)^4(1+X)^j}{n}\right).
\end{equation}
For $j>k$ within this logarithmic range the same argument uses the zero
falling factorial, so the formula remains valid with our convention.
The displayed remainder is summable against $|f_j|$ by
\eqref{tc:coefficient-bound}. Substitution into
\eqref{tc:exact-sum} proves \eqref{tc:critical-expansion}.
\end{proof}

The first coefficients are $f_0=1$, $f_1=-\sqrt\pi/2$, and $f_2=1/12$.
In particular, $R_{n,k}\to1$ when $k=o(\sqrt n)$, while for
$k/\sqrt n\to x>0$ it tends to $F(x)<1$. The optimal subcritical
range itself is already available in the recent enumeration literature
\cite{YZ26,FY26}. Theorem~\ref{tc:critical} describes its boundary and
the first uniform correction. It is a scaling expansion, not a finite
radius of convergence: $F$ is entire. The spectral formula also gives
\[
 F(x)=2\sqrt\pi\,x e^{-x^2/3+a_1x^{2/3}}
 \left(1+O(e^{-(a_1-a_2)x^{2/3}})\right)\qquad(x\to+\infty).
\]
This last relation concerns the scaling function; it does not extend
the uniformity range in Theorem~\ref{tc:critical} to growing $X$.

% END SECTION SOURCE: tc_critical.tex

\begin{figure}[tbp]
 \centering
 \includegraphics[width=\textwidth]{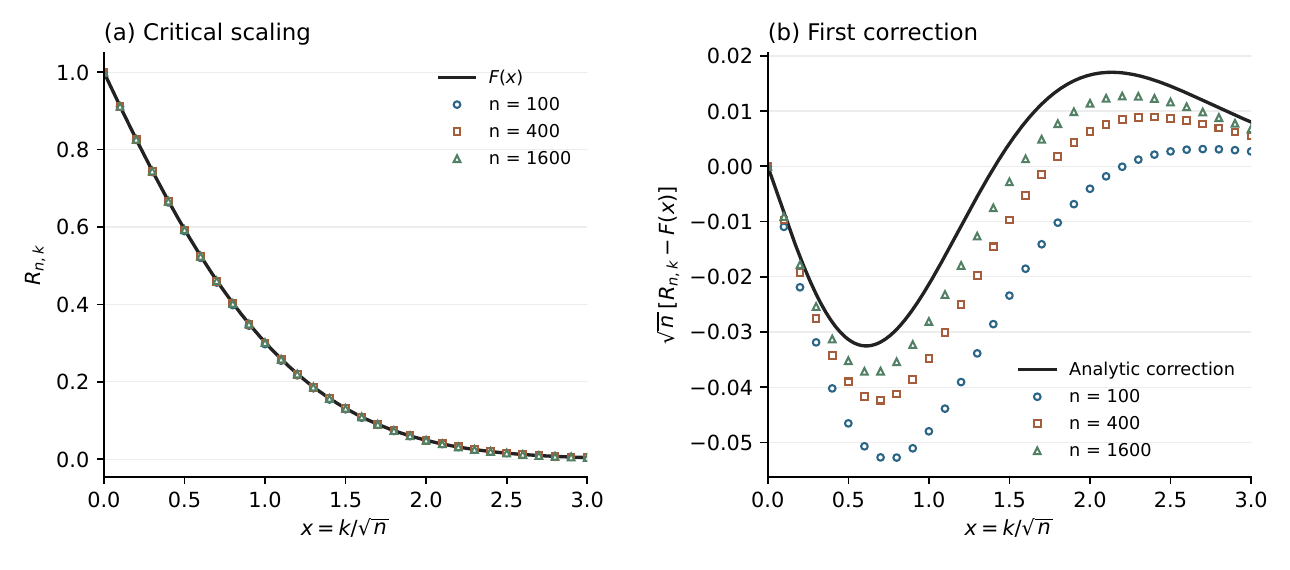}
 \caption{The critical tree-child regime. Left: the entire scaling
 function $F$ and normalized exact counts from the recurrence.
 Right: $\sqrt n(R_{n,k}-F(k/\sqrt n))$ approaches the analytic
 correction $-(xF''(x)+x^2F'(x))/2$. The points are deterministic
 recurrence evaluations; the curves involve no fitted parameters.}
 \label{fig:critical}
\end{figure}
% BEGIN SECTION SOURCE: tc_maximum_total.tex
% Article section. Requires amsmath, amssymb, amsthm, booktabs, hyperref.
% External references: Theorem tc:airy and formula tc:exact-sum.
% Bibliography keys: FYZ21, Janson07, PB21, Lin26, CFLLWY24, DLMF.
\providecommand{\Ai}{\operatorname{Ai}}

\section{Maximal reticulation: an explicit Airy amplitude}\label{sec:tc-maximum}\label{tcm:max-section-maximum}
Recall that $M_n=T_{n,n-1}$, and let $a=a_1$ be the largest negative zero of $\Ai$. Set
\[
D_n=(3n)^{-2/3}e^{a(3n)^{1/3}}\left(\frac{12}{e^2}\right)^n n^{2n}.
\]
Fuchs--Yu--Zhang proved $M_n=\Theta(D_n)$ \cite{FYZ21}. We refine their estimate, using the exact coefficient formula \eqref{tc:exact-sum} and the identity proved in Theorem~\ref{tc:airy}:
\begin{equation}\label{tcm:max-FAiry}
F(z)=\sum_{j\ge0}f_jz^j=e^{-z^2/3}\mathbb E e^{-\sqrt2 z\mathcal A},
\end{equation}
where $\mathcal A$ is standard Brownian excursion area. Throughout this section $F(z)=\sum_{j\ge0}f_jz^j$. The classical spectral formula for the excursion transform \cite[Eq.~(80)]{Janson07} gives
\begin{equation}\label{tcm:max-spectral}
F(z)=2\sqrt\pi\,z\,e^{-z^2/3}\sum_{\nu\ge1}e^{a_\nu z^{2/3}},
\qquad |\arg z|<3\pi/4,
\end{equation}
with the principal power. Normal convergence proves the complex extension from positive real $z$.

\begin{theorem}\label{tcm:maximum}
The maximal-reticulation count has the positive amplitude
\begin{equation}\label{tcm:max-K}
\boxed{K=\pi\Ai'(a_1)>0.}
\end{equation}
More precisely,
\begin{align}\label{tcm:max-main}
\log\frac{M_n}{D_n}={}&\log\left(\pi\Ai'(a_1)\right)
+\frac{a_1^2}{6}(3n)^{-1/3}
-a_1(3n)^{-2/3}
+\frac{13}{6}1(3n)^{-1}+O(n^{-4/3}).
\end{align}

\end{theorem}

\subsection{An exact coefficient integral}
Put $m=n-1$ and $c_m=M_{m+1}/(m+1)!$. The exact normalized coefficient formula is
\begin{equation}\label{tcm:max-finite}
c_m=\frac{2^{2m}}{\sqrt\pi}\sum_{j=0}^m
\frac{f_j}{(m-j)!}\Gamma\left(2m+\frac12-\frac j2\right).
\end{equation}
This follows by setting $n=m+1$ and $k=m$ in \eqref{tc:exact-sum}, and separating the parity-dependent double factorials. The identity $M_n=n!c_{n-1}$, equivalently $M_n=n!c_{n-1,n-1}$ in the word notation of the preceding section, was already proved in \cite[Proposition~2]{FYZ21}. Define
\[
Q_m(z)=\int_0^\infty u^{3m}e^{-u^2+zu}\,\mathrm d u.
\]
The gamma integral and Cauchy's formula give the exact identity
\begin{equation}\label{tcm:max-exact-int}
\boxed{c_m=\frac{2^{2m+1}}{\sqrt\pi}\frac1{2\pi i}
\oint \frac{F(z)Q_m(z)}{z^{m+1}}\,\mathrm d z.}
\end{equation}
Indeed, $[z^m]F(z)e^{zu}=\sum_{j=0}^mf_ju^{m-j}/(m-j)!$ and
$\int_0^\infty u^{4m-j}e^{-u^2}\,\mathrm d u=\frac12\Gamma(2m+\frac12-\frac j2)$.
The interchange is justified on every fixed circle by the integrable bound $u^{3m}e^{-u^2+r u}$.

\subsection{Global localization of the coefficient integral}
Choose the circle $|z|=\sqrt{3m}$ and write $z=\sqrt m\,\zeta$.

\subsubsection{The left semicircle is exponentially small}
The excursion is the norm of a three-dimensional Brownian bridge \cite{Janson07}. Cauchy--Schwarz and Jensen imply
\[
\mathbb E e^{\mathcal A^2}\le\int_0^1\mathbb E e^{\|B(t)\|^2}\,\mathrm d t\le2\sqrt2,
\]
because each bridge coordinate has variance $t(1-t)\le1/4$. Therefore, for $\Re z\le0$,
\[
|F(z)|\le2\sqrt2\exp\left(-\frac{\Re z^2}{3}+\frac{(\Re z)^2}{2}\right)
\le2\sqrt2 e^{|z|^2/3}.
\]
Also $|Q_m(z)|\le Q_m(0)=\frac12\Gamma((3m+1)/2)$. Stirling's formula shows that the contribution of this semicircle to \eqref{tcm:max-exact-int} is at most
\[
C\left(\frac{12}{e}\right)^m m^m e^{-\delta m},
\qquad \delta=\frac32\log2-\frac12>0.
\]
It is exponentially smaller than the final answer, including its stretched exponential.

\subsubsection{The inner Gaussian saddle}
For $\Re\zeta\ge0$, $|\zeta|=\sqrt3$, define
\begin{align*}
U(\zeta)&=\frac{\zeta+\sqrt{\zeta^2+24}}4,\\
q(\zeta)&=-U(\zeta)^2+\zeta U(\zeta)+3\log U(\zeta),\\
g(\zeta)&=\left(2+\frac3{U(\zeta)^2}\right)^{-1/2}.
\end{align*}
All branches are positive at $\zeta=\sqrt3$. Rotating the $u$-ray to $\arg U(\zeta)$ is valid: $|\arg U|<\pi/4$ uniformly. Along this ray the real phase is strictly concave, with its unique maximum at $|U|$. The Gaussian saddle expansion is uniform on the semicircle and in a fixed complex neighbourhood of $\sqrt3$:
\begin{equation}\label{tcm:max-Q-exp}
Q_m(\sqrt m\zeta)=m^{3m/2}\sqrt{2\pi}\,g(\zeta)e^{m q(\zeta)}
\left(1+\frac{b(\zeta)}m+O(m^{-2})\right),
\end{equation}
where, with $\kappa_0=2+3/U^2$,
\[
b(\zeta)=-\frac9{4U^4\kappa_0^2}+\frac{15}{2U^6\kappa_0^3}.
\]
To check the coefficient, write
$H(v;\zeta)=-v^2+\zeta v+3\log v$ for the inner phase.
At $v=U$ its derivatives are
$H''=-\kappa_0$, $H'''=6/U^3$, and $H''''=-18/U^4$.
The two Gaussian contributions are
$H''''/(8\kappa_0^2)$ and $5(H''')^2/(24\kappa_0^3)$.
After rotation, the real second derivative in the radial variable $r$ is
$-2\cos(2\arg U)-3/r^2<0$, uniformly bounded away from zero near
the saddle. The tails are uniformly separated from its maximum.
Taylor's formula on a fixed saddle neighbourhood, followed by Gaussian
integration, therefore gives the stated uniform $O(m^{-2})$ remainder.
The same argument applies in a fixed complex neighbourhood of $\sqrt3$;
there the coefficients and the remainder in \eqref{tcm:max-Q-exp} are
analytic in $\zeta$.

\subsubsection{A unique cubic maximum on the right semicircle}
Set
\begin{equation}\label{tcm:max-phase}
\phi(\zeta)=q(\zeta)-\zeta^2/3-\log\zeta.
\end{equation}
At $\zeta_0=\sqrt3$,
\[
\phi(\zeta_0)=\log3-1,\quad\phi'(\zeta_0)=\phi''(\zeta_0)=0,
\quad\phi'''(\zeta_0)=-\frac{16\sqrt3}{81}.
\]
There is also a global check. Write $\zeta=\sqrt3e^{i\theta}$, $U=Re^{i\eta}$, $t=R^2$. From $\zeta=2U-3/U$ one obtains
\begin{equation}\label{tcm:max-monotonephase}
\frac{\,\mathrm d}{\,\mathrm d\theta}\Re\phi(\sqrt3e^{i\theta})
=\frac23\left(t-\frac9t\right)\sin(2\eta)<0,
\qquad 0<\theta\le\pi/2.
\end{equation}
Here $3/2\le t<3$ and $0<\eta<\pi/4$. Indeed, $\Re\zeta=(2R-3/R)\cos\eta\ge0$ gives $t\ge3/2$, while $|\zeta|^2=4t+9/t-12\cos(2\eta)=3$ and $\cos(2\eta)\le1$ give $t\le3$, with equality only at $\theta=0$. The same identity gives $\cos(2\eta)\ge3/4$.
For the derivative in \eqref{tcm:max-monotonephase}, stationarity gives
$q'(\zeta)=U$, and hence
$\phi'(\zeta)=U-2\zeta/3-1/\zeta$. Substituting
$\zeta U=2U^2-3$ and $\zeta^2=4U^2-12+9/U^2$ into
$\Re(i\zeta\phi'(\zeta))$ gives the displayed expression.
By conjugation the other half is symmetric. Near zero,
\[
\Re\phi(\sqrt3e^{i\theta})=\phi(\sqrt3)-\frac{16}{81}\theta^4+O(\theta^6).
\]
Thus every fixed arc separated from $\theta=0$ is exponentially negligible.

On the right semicircle \eqref{tcm:max-spectral} is uniformly dominated by its first zero, since $\Re z^{2/3}\ge |z|^{2/3}/2$. The error from the other zeros is exponentially small in $m^{1/3}$. This remains true after integration: the inequalities
\[
\Re\phi(\sqrt3e^{i\theta})\le\phi(\sqrt3)-c\theta^4,
\qquad
\Re(a z^{2/3})\le a(3m)^{1/3}+C m^{1/3}\theta^2
\]
show, by completing a square in $\theta^2$, that the modulus integral loses only a polynomial factor relative to the principal exponential. That factor cannot offset the spectral gap $a_1-a_2>0$.

\subsection{The local Airy integral and its first nonzero term}
Substituting \eqref{tcm:max-spectral} and \eqref{tcm:max-Q-exp} into the localized integral gives, up to exponentially small errors,
\begin{equation}\label{tcm:max-reduced}
c_m=4^{m+1}\sqrt{2\pi}\,m^{m+1/2}\frac1{2\pi i}
\int g(\zeta)e^{m\phi(\zeta)+a m^{1/3}\zeta^{2/3}}
\left(1+\frac{b(\zeta)}m+O(m^{-2})\right)\,\mathrm d\zeta.
\end{equation}
Put
\[
\epsilon=m^{-1/3},\qquad \kappa=\frac{3^{7/6}}2,
\qquad \zeta=\sqrt3-\kappa\epsilon t,
\qquad \rho=\frac{3^{2/3}}2.
\]
The exponent becomes
\[
m\phi(\zeta)+a m^{1/3}\zeta^{2/3}
=m(\log3-1)+a(3m)^{1/3}+\frac{t^3}3-at+O(\epsilon).
\]
The original local arc is deformable to the two cubic descent rays.
After orientation reversal it induces the usual Airy contour from
$e^{-i\pi/3}\infty$ to $e^{i\pi/3}\infty$.
The minus sign in $\,\mathrm d\zeta=-\kappa\epsilon\,\mathrm d t$ is thereby
removed. The next lemma justifies this deformation and the expansion,
including the error after the leading coefficient vanishes.

\begin{lemma}[Termwise local expansion]\label{tcm:local-expansion}
After factoring out the exponential in \eqref{tcm:max-reduced}, its local integrand can be expanded through $\epsilon^4$ and integrated term by term, with integrated remainder $O(\epsilon^5)$. In particular, after the leading vanishing term is removed, the relative remainder is $O(\epsilon^4)$.
\end{lemma}
\begin{proof}
Choose a fixed small disc about $\zeta_0=\sqrt3$.
There is an analytic coordinate $\chi$ with
$\phi(\zeta)-\phi(\zeta_0)=\chi(\zeta)^3/3$ and
$\chi'(\zeta_0)=-1/\kappa$.
For the two endpoints of the original local arc,
$\Re\chi^3<0$ by \eqref{tcm:max-monotonephase}; moreover, their
arguments lie in the decay sectors bordering the positive and negative
imaginary axes. In the $\chi$-plane they can be joined to rays of
arguments $\pm\pi/3$ inside those sectors. The connecting pieces stay
away from zero, so their real phase is smaller than
$\phi(\zeta_0)$ by a fixed positive amount.
Pulling the connections back and, if necessary, making a further local
deformation lets us take the central portions to be straight rays in
the linear coordinate $\zeta=\zeta_0-\kappa\epsilon t$.
All functions are analytic in this disc, so no singularity is crossed.

On these rays the cubic term is $-|t|^3/3$.
The remaining phase terms have modulus at most
$C\epsilon|t|^4$, while the nonconstant part of the stretched
exponential is at most $C|t|$ whenever $|\epsilon t|$ stays in
the fixed disc. Shrinking the disc gives the uniform majorant
$\exp(-c|t|^3+C|t|)$. The higher analytic coefficients are polynomials
in $t$ multiplied by the corresponding powers of $\epsilon$.

For example, one may first restrict to $|t|\le m^{1/100}$, where ordinary Taylor expansion with remainder gives an integrable polynomial majorant times $e^{-c|t|^3+C|t|}$. The rest of each ray is smaller than every power of $m^{-1}$ by the same majorant. The Gaussian expansion \eqref{tcm:max-Q-exp} is analytic and uniform on the disc. These bounds justify every finite termwise expansion and give the stated $O(\epsilon^5)$ integrated remainder. They also justify replacing the finite descent rays by infinite Airy rays.
\end{proof}

The order-$\epsilon$ correction to the local amplitude is particularly revealing:
\begin{equation}\label{tcm:max-P1}
P_1(t)=3^{2/3}\left(-\frac{t}{18}-\frac{a t^2}{12}+\frac{5t^4}{36}\right).
\end{equation}
Using the Airy contour identity \cite[Eq.~9.5.4]{DLMF},
\[
\frac1{2\pi i}\int_{\mathcal C}t^r e^{t^3/3-at}\,\mathrm d t=(-1)^r\Ai^{(r)}(a),
\]
the zeroth-order integral is $\Ai(a)=0$. Furthermore
$\Ai''(a)=0$ and $\Ai^{(4)}(a)=2\Ai'(a)$, so \eqref{tcm:max-P1} integrates to
\begin{equation}\label{tcm:max-firstnonzero}
\frac1{2\pi i}\int_{\mathcal C}P_1(t)e^{t^3/3-at}\,\mathrm d t
=3^{-1/3}\Ai'(a).
\end{equation}
This is the origin of the multiplicative constant, including its normalization; it is not fitted to counting data.

\subsection{Evaluation of the first corrections}
To make the calculation reproducible, use $w=1-\zeta/\sqrt3=\rho\epsilon t$ and
\[
v(w)=\frac{1-w+\sqrt{(1-w)^2+8}}4,\qquad U=\sqrt3v.
\]
Direct expansion yields
\begin{align*}
\phi(\zeta)-\phi(\sqrt3)
&=\frac8{27}w^3+\frac{20}{81}w^4+\frac{146}{729}w^5
+\frac{3283}{19683}w^6+\frac{2812}{19683}w^7+O(w^8),\\
\frac{g(\zeta)}{g(\sqrt3)}
&=1-\frac w9-\frac{2w^2}{81}+\frac{7w^4}{4374}+O(w^5),\\
b(\zeta)&=-\frac{17}{972}-\frac{8w}{729}+O(w^2).
\end{align*}
Also $a m^{1/3}\zeta^{2/3}=a3^{1/3}\epsilon^{-1}(1-w)^{2/3}$. Multiplying these finite expansions, the coefficients after integration against $e^{t^3/3-at}/(2\pi i)$ are as follows. Each table entry is divided by $\Ai'(a)$.
\begin{center}
\begin{tabular}{cc}
\toprule
Order & Integrated coefficient divided by $\Ai'(a)$\\
\midrule
$\epsilon^0$ & $0$\\[2pt]
$\epsilon^1$ & $3^{-1/3}$\\[2pt]
$\epsilon^2$ & $3^{1/3}a^2/18$\\[2pt]
$\epsilon^3$ & $a^4/216$\\[2pt]
$\epsilon^4$ & $3^{2/3}(a^6-108)/11664$\\
\bottomrule
\end{tabular}
\end{center}
Here is a finite recurrence that verifies every entry of the table.
Let $p_j$ be the coefficient of $w^j$ in the displayed phase expansion
and put, for $1\le j\le4$,
\[
 E_j(t)=p_{j+3}\rho^{j+3}t^{j+3}
 +a3^{1/3}\binom{2/3}{j+1}(-\rho)^{j+1}t^{j+1}.
\]
The coefficients of the exponential correction are determined by
\begin{equation}\label{tcm:max-exp-recurrence}
 e_0(t)=1,\qquad
 e_j(t)=\frac1j\sum_{r=1}^j rE_r(t)e_{j-r}(t).
\end{equation}
Multiply $\sum_{j=0}^4e_j(t)\epsilon^j$ by the displayed expansion
of $g/g(\sqrt3)$ and by $1+\epsilon^3b(\zeta)$, retaining degrees
at most four; call the resulting polynomials $P_j(t)$.
To integrate them, define polynomial pairs by
\[
 A_0(x)=1,\quad B_0(x)=0,\qquad
 A_{r+1}=A_r'+xB_r,\quad B_{r+1}=A_r+B_r'.
\]
Then $\Ai^{(r)}=A_r\Ai+B_r\Ai'$ follows by induction from
$\Ai''=x\Ai$, and at the zero $a$ one obtains
\begin{equation}\label{tcm:max-integral-recurrence}
 \frac1{2\pi i\Ai'(a)}\int_{\mathcal C}
       P_j(t)e^{t^3/3-at}\,\mathrm d t
 =\sum_{r\ge0}[t^r]P_j(t)(-1)^rB_r(a).
\end{equation}
Equations \eqref{tcm:max-exp-recurrence} and
\eqref{tcm:max-integral-recurrence}, using the rational phase coefficients
above, give the table by finite polynomial arithmetic.

Write $A=a^2/6$. After factoring out $3^{-1/3}\Ai'(a)\epsilon$, the table becomes
\[
1+\frac{A}{3^{1/3}}\epsilon+\frac{A^2}{2\,3^{2/3}}\epsilon^2+
\left(\frac{A^3}{18}-\frac1{36}\right)\epsilon^3+O(\epsilon^4).
\]
Therefore \eqref{tcm:max-reduced} gives
\begin{equation}\label{tcm:max-am}
\boxed{\begin{aligned}
c_m={}&2\,3^{1/3}\sqrt{2\pi}\,\Ai'(a)
\left(\frac{12}{e}\right)^m m^{m-1/6}e^{a(3m)^{1/3}}\\
&\times\exp\left(A(3m)^{-1/3}-\frac1{36m}+O(m^{-4/3})\right).
\end{aligned}}
\end{equation}
The prefactor follows from
$4\kappa\sqrt{2\pi/3}\,3^{-1/3}=2\,3^{1/3}\sqrt{2\pi}$.

Finally $M_n=n!c_{n-1}$. Stirling's formula and $m=n-1$ give
\[
K=3^{2/3}\frac{2\,3^{1/3}\sqrt{2\pi}\,\Ai'(a)\sqrt{2\pi}}{12}
=\pi\Ai'(a).
\]
The shift of the stretched exponential gives $-a(3n)^{-2/3}=\beta(3n)^{-2/3}$. At order $n^{-1}$ the three contributions are
\[
\frac1{12}+\frac23-\frac1{36}=\frac{13}{18}.
\]
This proves \eqref{tcm:max-main}. Notice that the $m^{-2/3}$ term in the logarithm of \eqref{tcm:max-am} is zero; the nonzero $n^{-2/3}$ term comes entirely from the index shift in the stretched exponential.

\section{Summation over the reticulation defect}\label{tcm:tot-section-total}
Let $T_n=\sum_{k=0}^{n-1}T_{n,k}$, including the trees at $k=0$.
Retain $D_n$, $K$, $A=a_1^2/6$, and $\beta=-a_1>0$ from the preceding section.
We use the maximum-count expansion proved in Theorem~\ref{tcm:maximum}:
\begin{equation}\label{tcm:tot-input}
M_n=KD_n\exp\left(A(3n)^{-1/3}+\beta(3n)^{-2/3}+\frac{13}{6}(3n)^{-1}+O(n^{-4/3})\right).
\end{equation}
The exact word/network identity and the recurrence below are the Pons--Batle identity, initially conjectured in \cite{PB21} and now proved; see \cite{Lin26}. The leading relation $T_n\sim\sqrt e M_n$ and the Poisson$(1/2)$ limit of the reticulation defect were proved in \cite[Theorem 1.9 and Corollary 1.11]{CFLLWY24}. Our calculation refines the relation between the total and the maximum.

\begin{theorem}\label{tcm:total}
As $n\to\infty$,
\begin{equation}\label{tcm:tot-ratio}
\boxed{\frac{T_n}{M_n}=\sqrt e\,\exp\left(-\frac{\beta}{24}(3n)^{-2/3}+\frac1{288n}+O(n^{-4/3})\right).}
\end{equation}
In particular,
\begin{equation}\label{tcm:tot-total}
\boxed{T_n=K_{\rm tot}D_n\exp\left(A(3n)^{-1/3}+\frac{23\beta}{24}(3n)^{-2/3}+\frac{209}{96}(3n)^{-1}+O(n^{-4/3})\right),}
\end{equation}
where
\[
K_{\rm tot}=\sqrt e\,K=\sqrt e\,\pi\Ai'(a_1).
\]
No additional amplitude constant is introduced by summing over reticulations.
\end{theorem}
The three terms here are the three corrections of orders $n^{-1/3},n^{-2/3},n^{-1}$, in the same logarithmic convention as \eqref{tcm:tot-input}. In particular, multiplying \emph{every term} of the maximum-count approximation by $\sqrt e$ would not give the correct second and third corrections.

\subsection{An exact recurrence at fixed defect}
Let $c_{m,k}$ be the number of Pons--Batle words. Then
\begin{equation}\label{tcm:tot-words}
c_{m,k}=c_{m,k-1}+(2m+k-1)c_{m-1,k},\qquad
T_{n,k}=\frac{n!}{(n-k)!}c_{n-1,k},
\end{equation}
with $c_{m,0}=(2m-1)!!$ and $c_{m-1,m}=0$. Thus $c_{m,m}=c_{m,m-1}$. Put
\[
c_m=c_{m,m},\qquad b_{m,d}=\frac{c_{m,m-d}}{c_m}\quad(0\le d\le m).
\]
Here $b_{m,d}$ denotes the normalized boundary row. With $m=n-1$,
\begin{equation}\label{tcm:tot-defectratio}
\frac{T_{n,n-1-d}}{M_n}=\frac{b_{m,d}}{(d+1)!}.
\end{equation}
Positivity of \eqref{tcm:tot-words} gives the elementary bound
\begin{equation}\label{tcm:tot-tailbound}
0\le b_{m,d}\le1,\qquad
0\le\frac{T_{n,n-1-d}}{M_n}\le\frac1{(d+1)!}.
\end{equation}
This weaker-than-optimal bound is already sufficient for all tails below.

For $1\le d\le m-1$, equation \eqref{tcm:tot-words} gives
\begin{equation}\label{tcm:tot-boundaryrec}
\boxed{b_{m,d+1}=b_{m,d}-\gamma_{m,d}b_{m-1,d-1},\qquad
\gamma_{m,d}=(3m-d-1)\frac{c_{m-1}}{c_m}.}
\end{equation}
The initial values in the defect are exactly $b_{m,0}=b_{m,1}=1$. In particular, $T_{n,n-2}=M_n/2$ for $n\ge2$.

\begin{lemma}\label{tcm:tot-adjacent-lemma}
The input \eqref{tcm:tot-input} implies
\begin{equation}\label{tcm:tot-adjacent}
\frac{c_{m-1}}{c_m}=\frac1{12m}
\left(1+\beta(3m)^{-2/3}+\frac2{3m}+O(m^{-4/3})\right).
\end{equation}
Consequently, uniformly for $1\le d\le m-1$,
\begin{equation}\label{tcm:tot-cm}
\gamma_{m,d}=\frac14+\frac\beta4(3m)^{-2/3}+\frac{1-d}{12m}
+O\big((d+1)m^{-4/3}\big).
\end{equation}
For all sufficiently large $m$, also $0<\gamma_{m,d}\le1$ throughout this range.
\end{lemma}
\begin{proof}
Since $c_m=M_{m+1}/(m+1)!$, Stirling's formula and \eqref{tcm:tot-input} give, for a constant $C_*$,
\begin{equation}\label{tcm:tot-logam}
\log c_m=C_*+m\log(12/e)+(m-1/6)\log m
+a3^{1/3}m^{1/3}+A(3m)^{-1/3}-\frac1{36m}+O(m^{-4/3}).
\end{equation}
Taking the backward difference yields
\[
\log\frac{c_m}{c_{m-1}}
=\log(12m)-\beta(3m)^{-2/3}-2(3m)^{-1}+O(m^{-4/3}).
\]
Here the two remainder terms are bounded separately: no differentiation or additional smoothness of an $O$-term is assumed. Exponentiating proves \eqref{tcm:tot-adjacent}. Multiplication by $3m-d-1$ gives \eqref{tcm:tot-cm}. The assertion $\gamma_{m,d}\le1$ follows from \eqref{tcm:tot-adjacent} and $3m-d-1<3m$.
\end{proof}

\subsection{The expansion and its uniform remainder}
Define, for $d\ge0$,
\begin{align}\label{tcm:tot-uvw}
u_d&=\frac{d+1}{2^d},&
v_d&=-\frac{(d+1)d(d-1)}{6\,2^d},&
w_d&=\frac{(d+1)d(d-1)(d-2)}{36\,2^d},
\end{align}
and put
\begin{equation}\label{tcm:tot-P}
P_{m,d}=u_d+\beta v_d(3m)^{-2/3}+w_d m^{-1}.
\end{equation}
The polynomial factors in \eqref{tcm:tot-uvw} vanish when their orders exceed $d$; in particular $P_{m,0}=P_{m,1}=1$.

\begin{lemma}[Uniform, summable error]\label{tcm:tot-uniform}
There are constants $C>0$ and an integer $m_0\ge2$ such that
\begin{equation}\label{tcm:tot-error}
|b_{m,d}-P_{m,d}|\le\frac{C\,2^d}{(m-d)^{4/3}}
\qquad(m-d\ge m_0).
\end{equation}
In particular, for $0\le d\le m/2$ and $m\ge2m_0$,
\begin{equation}\label{tcm:tot-weightederror}
\left|\frac{T_{m+1,m-d}}{M_{m+1}}-
\frac{1}{2^d d!}\left(1-\frac{\beta d(d-1)}6(3m)^{-2/3}
+\frac{d(d-1)(d-2)}{36m}\right)\right|
\le C'\frac{2^d}{(d+1)!}m^{-4/3}.
\end{equation}
\end{lemma}
\begin{proof}
Substitution into \eqref{tcm:tot-boundaryrec} gives the coefficient equations
\begin{align*}
u_{d+1}&=u_d-\tfrac14u_{d-1},\\
v_{d+1}&=v_d-\tfrac14v_{d-1}-\tfrac14u_{d-1},\\
w_{d+1}&=w_d-\tfrac14w_{d-1}+\tfrac{d-1}{12}u_{d-1}.
\end{align*}
They hold identically for \eqref{tcm:tot-uvw}; $u_0=u_1=1$ and $v_0=v_1=w_0=w_1=0$ fix the solutions. These identities can also be read from
\[
\sum_{d\ge0}u_dz^d=(1-z/2)^{-2},\quad
\sum_{d\ge0}v_dz^d=-\tfrac14z^2(1-z/2)^{-4},\quad
\sum_{d\ge0}w_dz^d=\tfrac1{12}z^3(1-z/2)^{-5}.
\]
Consider the residual
\[
\Delta_{m,d}=P_{m,d+1}-P_{m,d}+\gamma_{m,d}P_{m-1,d-1}.
\]
Use \eqref{tcm:tot-cm}, the three coefficient identities, and
$(m-1)^{-2/3}=m^{-2/3}+O(m^{-5/3})$ and
$(m-1)^{-1}=m^{-1}+O(m^{-2})$. All surviving products are bounded by
\[
|\Delta_{m,d}|\le C_1(d+1)^5 2^{-d}m^{-4/3}
\le C_0m^{-4/3},
\qquad 1\le d\le m-1.
\]
The constants do not depend on $d$: polynomial multiples of $2^{-d}$ are bounded. The factor $(d+1)$ in the remainder of \eqref{tcm:tot-cm} is absorbed by the displayed fifth-degree polynomial.

Set $e_{m,d}=b_{m,d}-P_{m,d}$. Then $e_{m,0}=e_{m,1}=0$ and
\[
e_{m,d+1}=e_{m,d}-\gamma_{m,d}e_{m-1,d-1}-\Delta_{m,d}.
\]
Choose $m_0$ so that $\gamma_{m,d}\le1$ for $m\ge m_0$, and choose $C\ge C_0$. Induct on $d$, simultaneously for every admissible $m$. If $m-d-1\ge m_0$, the inductive bounds give
\begin{align*}
|e_{m,d+1}|
&\le\frac{C(2^d+2^{d-1})}{(m-d)^{4/3}}+C_0m^{-4/3}\\
&\le\frac{C(2^d+2^{d-1})+C_0}{(m-d-1)^{4/3}}
\le\frac{C2^{d+1}}{(m-d-1)^{4/3}},
\end{align*}
because $d\ge1$. This proves \eqref{tcm:tot-error}. Divide by $(d+1)!$ and use $m-d\ge m/2$ to obtain \eqref{tcm:tot-weightederror}.
\end{proof}

For each fixed $d$, the lemma gives the more familiar relative form
\begin{equation}\label{tcm:tot-fixed}
\frac{T_{n,n-1-d}}{M_n}=
\frac1{2^d d!}\left(1-\frac{\beta d(d-1)}6(3n)^{-2/3}
+\frac{d(d-1)(d-2)}{36n}+O_d(n^{-4/3})\right).
\end{equation}
However, it is \eqref{tcm:tot-weightederror}, not the fixed-$d$ notation $O_d$, that justifies the next step.

\subsection{Summation over all defects}
It is useful to retain a marker for the defect. Define
\[
\mathcal G_n(z)=\sum_{d=0}^{n-1}\frac{T_{n,n-1-d}}{M_n}z^d.
\]
\begin{proposition}\label{tcm:defect-generating-function}
For every fixed $R>0$, uniformly on $|z|\le R$,
\begin{equation}\label{tcm:tot-marker}
\boxed{\mathcal G_n(z)=e^{z/2}
\left(1-\frac{\beta z^2}{24}(3n)^{-2/3}+\frac{z^3}{96}(3n)^{-1}+O_R(n^{-4/3})\right).}
\end{equation}
\end{proposition}
\begin{proof}
First sum over $d\le m/2$, with $m=n-1$. Lemma \ref{tcm:tot-uniform} bounds the sum of absolute errors by
\[
C' m^{-4/3}\sum_{d\ge0}\frac{(2R)^d}{(d+1)!}=O_R(m^{-4/3}).
\]
For $d>m/2$, \eqref{tcm:tot-tailbound} gives a factorial tail, uniformly on the same disk. The tails of the approximating sequences $P_{m,d}/(d+1)!$ are also smaller than every inverse power of $m$, since they are bounded by polynomial multiples of $(R/2)^d/(d+1)!$. We may therefore extend all approximating sums to infinity with an error $O_R(m^{-4/3})$.

For $r\ge0$, the elementary identity
\[
\sum_{d\ge0}\frac{(d)_{\underline r}}{d!}\left(\frac z2\right)^d
=\left(\frac z2\right)^r e^{z/2}
\]
gives
\begin{align*}
\sum_{d\ge0}\frac{u_dz^d}{(d+1)!}&=e^{z/2},\\
\sum_{d\ge0}\frac{v_dz^d}{(d+1)!}&=-\frac{z^2}{24}e^{z/2},\\
\sum_{d\ge0}\frac{w_dz^d}{(d+1)!}&=\frac{z^3}{288}e^{z/2}.
\end{align*}
Finally, replacing $m$ by $n-1$ changes $m^{-2/3}$ by $O(n^{-5/3})$ and $m^{-1}$ by $O(n^{-2})$, both absorbed by the claimed remainder. The factor $e^{z/2}$ is nonzero and bounded away from zero on the fixed disk.
\end{proof}

\begin{proof}[Proof of Theorem \ref{tcm:total}]
Set $z=1$ in \eqref{tcm:tot-marker}. Since the square of its first correction is $O(n^{-4/3})$, taking a logarithm proves \eqref{tcm:tot-ratio}. Combining with \eqref{tcm:tot-input} yields
\[
\beta-\frac\beta{24}=\frac{23\beta}{24},\qquad
\frac{13}{18}+\frac1{288}=\frac{209}{288}.
\]
The leading constant is multiplied by $\sqrt e$, which proves \eqref{tcm:tot-total}.
\end{proof}

% END SECTION SOURCE: tc_maximum_total.tex

% BEGIN SECTION SOURCE: rv_combinatorics.tex
% Self-contained combinatorial sections for the joint article.
% Published inputs: CF24 (component representation and local recurrence),
% FYZ21 (maximal tree-child coding; a direct proof is also supplied below).
% All labels in this file have prefix rvc:.

\section{Reticulation-visible networks and component skeletons}
\label{rvc:section}

Throughout this part, $n\ge2$. Networks are rooted, simple directed
acyclic graphs with a planted root of outdegree one and leaves labelled
by $[n]$. In a binary network, a tree vertex has indegree one and
outdegree two, and a reticulation has indegree two and outdegree one.
A vertex is \emph{visible} if some leaf is reached from the root only
through that vertex. A network is reticulation-visible if every
reticulation is visible. A network is tree-child if every non-leaf
vertex has a child which is not a reticulation.

Write $\operatorname{RV}_{n,r}$ for the number of binary
reticulation-visible networks with $n$ leaves and $r$ reticulations,
and set
\begin{equation}
 m=n-1,\qquad M_n=T_{n,n-1},\qquad
 V_{n,q}=\operatorname{RV}_{n,3n-3-q},\qquad
 Z_n(z)=\frac1{M_n}\sum_{q=0}^{3n-3}V_{n,q}z^q.
 \label{rvc:normalization}
\end{equation}
Thus $q$ measures the distance from the largest possible reticulation
number. The parameter $z$ is an activity for this defect.

\subsection{The published reduction and its normalization}

We use the component construction of Chang and Fuchs~\cite{CF24}.
Remove the incoming arcs of every reticulation and compress each
remaining tree component to a vertex, retaining its labelled leaves.
Two arcs from one component to another are represented by one marked
arc. Their component theorem characterizes the resulting unmarked
graphs of reticulation-visible networks: they are tree-child graphs
of indegree at most two, with no reticulation having a single child
which is itself a tree vertex. Their reconstruction replaces each
internal component vertex by a one-component network; marked outgoing
arcs specify its reticulated terminal ports. This characterization
and reconstruction are the published inputs used below.

For our purposes it is convenient to normalize the component graph.
Split every vertex of indegree two and outdegree at least two into a
reticulation of outdegree one followed by a tree vertex. The new arc
is called \emph{artificial}. Restore the planted root, which is not
included among the tree vertices. The resulting skeleton $H$ is a
semibinary tree-child network: all reticulations have degree $(2,1)$,
whereas a tree vertex $v$ can have any outdegree $c(v)\ge2$.
Write $b(v)\ge1$ for its number of nonreticulation children, and
$\ell(v)$ for the number of these children which are leaves.

This normalization is reversible: contract every
reticulation--tree arc. Such arcs are exactly the artificial arcs,
because the unnormalized component graph excludes reticulations
whose sole child is a tree vertex. A reticulation followed directly
by a leaf is left unchanged. In the normalized skeleton the marking
rules are therefore as follows. Arcs from a tree vertex to a tree
vertex are marked; arcs from a tree vertex to a leaf are optionally
marked. Arcs entering reticulations, artificial arcs, and the planted
arc are unmarked. In particular, artificial arcs do not contribute
an additional reticulation in the reconstruction.

We record a useful elementary observation.

\begin{lemma}\label{rvc:rigidity}
A tree-child network with fixed leaf labels has no nontrivial
automorphism preserving those labels.
\end{lemma}

\begin{proof}
From any vertex $v$, repeatedly choose a nonreticulation child until
a leaf is reached. All vertices after $v$ on this path have indegree
one. An automorphism fixes the labelled leaf and hence, by following
unique parents backwards, fixes $v$. This applies to every vertex.
\end{proof}

Consequently the children at each skeleton vertex are distinguishable.
Choosing once and for all an ordering of these children identifies
them with labelled ports. The number of local replacements depends
only on the number of ports and on which of them are reticulated,
not on this auxiliary ordering.

\subsection{One-component weights and the defect identity}

A one-component network is a binary network in which every
reticulation is followed directly by a leaf. Let $w_{c,j}$ count such
networks on $c$ labelled ports, with a prescribed set of $j$ ports
below reticulations. This is the quantity denoted by $M_{c,j}$ in
\cite{CF24}; it is different from our $M_n$.
The published recurrence and initial values are
\begin{align}
w_{c,0}&=(2c-3)!!,\nonumber\\
w_{c,1}&=(c-1)(2c-3)!!,\nonumber\\
w_{c,j}&=(c+j-2)w_{c,j-1}+(j-1)w_{c,j-2}\nonumber\\
&\quad+\frac12\sum_{d=1}^{j-1}\binom{j-1}{d}(2d-1)!!
\bigl(w_{c-d,j-1-d}-w_{c+1-d,j-1-d}\bigr).
\label{rvc:localrecurrence}
\end{align}
The recurrence applies for $2\le j\le c$, and we use $(-1)!!=1$.
In particular,
$w_{2,0}=w_{2,1}=1$ and $w_{2,2}=3$.

At a normalized tree vertex $v$, omit the optional mark on exactly
$s_v$ of its $\ell(v)$ leaf arcs. The local replacement then has
weight $w_{c(v),b(v)-s_v}$. There are
$\binom{\ell(v)}{s_v}$ choices for the omitted marks. All other
nonreticulation outgoing arcs are marked. The component
reconstruction of~\cite{CF24} makes these local choices independent.

Let $\kappa$ and $t$ be the numbers of reticulations and tree vertices
of $H$, and define
\begin{equation}
 \delta(H)=n-1-\kappa,\qquad
 e(H)=\sum_{v\text{ tree}}(c(v)-2),\qquad
 s=\sum_v s_v.
 \label{rvc:threeparameters}
\end{equation}
The parameter $s$ has an intrinsic interpretation in the reconstructed
RV network: it is the number of original leaves with a tree parent.
Indeed, a leaf port with an omitted mark has a tree parent in its
one-component replacement, while a marked leaf port has a
reticulation parent. A leaf already below a skeleton reticulation
continues to have a reticulation parent.

\begin{lemma}\label{rvc:budget}
For every normalized skeleton and every reconstruction,
\begin{equation}
 t=n+\kappa-1-e,\qquad
 \sum_{v\text{ tree}}(b(v)-1)=\delta,\qquad
 q=\delta+e+s.
 \label{rvc:budgeteq}
\end{equation}
In particular $\delta,e,s\ge0$.
\end{lemma}

\begin{proof}
Counting arcs by tails and by heads gives
$1+2t+e+\kappa=t+n+2\kappa$, proving the first identity.
There are $n+t$ nonreticulation vertices other than the planted
root. One is reached from the planted root, and $\kappa$ are
children of reticulations. Thus
$\sum_v b(v)=n+t-\kappa-1$, which proves the second identity.
There are $\sum_v b(v)-s$ marked arcs. The reconstructed network
therefore has
\[
 r=\kappa+\sum_v b(v)-s=n+t-1-s
\]
reticulations. Subtracting this from $3n-3$ proves the last identity.
Nonnegativity follows from $c(v)\ge2$, $b(v)\ge1$, and the definition
of $s$.
\end{proof}

We can keep track of these statistics simultaneously. Put
\[
 Z_n(z,y;v)=\frac1{M_n}\sum_{N\text{ binary RV on }[n]}
                  z^{q(N)}y^{s(N)}v^{\delta(N)},
 \qquad Z_n(z,y)=Z_n(z,y;1).
\]
Here $\delta(N)$ refers to the canonical normalized component
skeleton of $N$. In particular $Z_n(z)=Z_n(z,1;1)$.
The component formula and Lemma~\ref{rvc:budget} give the exact identity
\begin{equation}
 Z_n(z,y;v)=\frac1{M_n}\sum_H
       v^{\delta(H)}z^{\delta(H)+e(H)}
       \prod_{a\text{ tree in }H}
       \left(\sum_{j=0}^{\ell(a)}\binom{\ell(a)}j
                   w_{c(a),b(a)-j}(zy)^j\right).
 \label{rvc:skeletonsum}
\end{equation}
The sum is over normalized skeletons with leaf set $[n]$.
The case $q=0$ requires $\delta=e=s=0$. Its skeleton is a binary
tree-child network with $n-1$ reticulations; every local weight is
$w_{2,1}=1$. We recover the maximal-reticulation identity
of~\cite{CF24},
\begin{equation}
 \operatorname{RV}_{n,3n-3}=M_n,
 \label{rvc:maximalcount}
\end{equation}
as well as the upper bound $3n-3$ itself.

\subsection{Uniform bounds for one component}

The following bound will control all local replacements, including
vertices of arbitrarily large outdegree.

\begin{lemma}\label{rvc:localbounds}
For $c\ge2$ and $0\le j\le c$, let
$U_{c,j}=(2c+2j-3)!!/2^j$. Then
\begin{equation}
 \frac35U_{c,j}\le w_{c,j}\le U_{c,j}.
 \label{rvc:localboundseq}
\end{equation}
\end{lemma}

\begin{proof}
For the lower bound, unfold each terminal reticulation into two
distinctly labelled copies of its leaf. This gives a binary tree
on $c+j$ leaves. Conversely, merging the prescribed pairs gives a
valid one-component network whenever neither copy pair forms a
cherry. A prescribed pair is a cherry in a fraction
$1/[2(c+j)-3]$ of the trees. The union bound shows that at least
\[
 1-\frac{j}{2(c+j)-3}\ge\frac35
\]
of the trees are admissible. The orbits under the $j$ independent
copy interchanges are the one-component networks. Each orbit has
size at most $2^j$, proving the lower bound. No freeness of this
action is assumed.

For the upper bound, adjoining a new ordinary leaf at the root is
an injection proving $w_{a+1,k}\ge w_{a,k}$. Hence the difference
sum in~\eqref{rvc:localrecurrence} is nonpositive, and
\[
 w_{c,k}\le(c+k-2)w_{c,k-1}+(k-1)w_{c,k-2}.
\]
Fix $c$ and write $R_k=w_{c,k}/U_{c,k}$. Then
\[
 R_k\le a_kR_{k-1}+b_kR_{k-2},\qquad
 a_k=\frac{c+k-2}{c+k-3/2},\qquad
 b_k=\frac{k-1}{(c+k-3/2)(c+k-5/2)}.
\]
For $k\le c-1$, $a_k+b_k\le1$. Also,
$R_0=1$ and $R_1=1-1/(2c-1)$. For $c\ge3$, the $k=2$ bound gives
$R_2\le1-\varepsilon_c$, where
$\varepsilon_c=2c/(4c^2-1)$; the same bound holds for $R_1$.
Induction yields $R_1,\ldots,R_{c-1}\le1-\varepsilon_c$.
At the last step,
\[
 a_c+b_c=1+\frac1{(4c-3)(4c-5)}.
\]
Since $\varepsilon_c\ge1/[(4c-3)(4c-5)]$, this still gives
$R_c\le1$. For $c=2$ the assertion follows from the three initial
values displayed after~\eqref{rvc:localrecurrence}.
\end{proof}

\section{An exact expansion over maximal tree-child networks}
\label{rvc:refinementsection}

Let $\mathcal M_n$ be the set of binary tree-child networks on $[n]$
with $m=n-1$ reticulations. Each of their tree vertices has exactly
one reticulation child and one nonreticulation child. Consequently
their tree components are $n$ paths, each ending at a leaf; a path
can have length zero if a reticulation is followed directly by a
leaf. There are $2m$ tree vertices in total. We call the identity
of a tree vertex's reticulation child its \emph{port}. Each
reticulation identity occurs as a port exactly twice.

\subsection{Refinement, inversion, and multiplicity}

\begin{proposition}\label{rvc:refinement}
For $B\in\mathcal M_n$, let $\Pi(B)$ be the set of partitions of
the sequences of tree vertices on its tree paths into consecutive
nonempty blocks. A path with no tree vertices has a unique empty
partition. Refinements are counted as pairs $(B,\pi)$, up to
isomorphisms preserving the leaf labels and the block partition;
the blocks carry no additional labels.

A normalized skeleton $H$ has exactly
\begin{equation}
 \prod_{a\text{ tree in }H}\mathcal N(c(a),b(a)),
 \qquad
 \mathcal N(c,b)=b\frac{(c+b-2)!}{2^{b-1}},
 \label{rvc:fibresize}
\end{equation}
maximal refinements $(B,\pi)$ whose blocks invert to $H$. Conversely,
every pair $B\in\mathcal M_n$, $\pi\in\Pi(B)$ determines a normalized
skeleton, denoted by $Q(B,\pi)$.
For a block with $t$ ports containing $r$ repeated port pairs,
the corresponding skeleton vertex has
\begin{equation}
 b=r+1,\qquad c=t-r+1.
 \label{rvc:blockparameters}
\end{equation}
\end{proposition}

\begin{proof}
At a skeleton vertex with $c$ children, choose one of its $b$
nonreticulation children as the continuation. Insert a new
reticulation above each of the other $b-1$ nonreticulation children.
Replace the vertex by a path of $t=c+b-2$ binary tree vertices.
Its ports consist of each old reticulation child once and each new
reticulation twice. Arrange this multiset in any order, and join
the last tree vertex to the chosen continuation. Each order is
valid: consecutive copies form a triangle with distinct parents,
not parallel arcs. The construction preserves acyclicity and the
tree-child condition. Its number of choices is exactly
$b\,t!/2^{b-1}$.

The operations at different vertices are independent. They add
$\sum_a(b(a)-1)=\delta(H)$ reticulations, giving a network in
$\mathcal M_n$. The paths replacing the original vertices supply
the distinguished blocks.

For the inverse, fix $(B,\pi)$. For a tree vertex $x$, let $\beta(x)$
be its block, and put $\beta(x)=x$ for a leaf. For each reticulation
$h$, write $p_1(h),p_2(h)$ for its parents and $d(h)$ for its child.
Before making any contraction, define
\[
 R_{\mathrm{int}}
 =\{h:\beta(p_1(h))=\beta(p_2(h))\}.
\]
The vertices of $Q(B,\pi)$ are the planted root $\rho$, the original
labelled leaves, the blocks, and the reticulations outside
$R_{\mathrm{int}}$. Let $d(\rho)$ be the child of the root, and let
$\sigma(I)$ be the nonreticulation child of the last tree vertex
of $I$. Its arcs are specified simultaneously by
\begin{align*}
 \rho&\longrightarrow\beta(d(\rho)),\\
 I&\longrightarrow\beta(\sigma(I))
       &&(I\in\pi),\\
 \beta(p_i(h))&\longrightarrow h
       &&(h\notin R_{\mathrm{int}},\ i=1,2),\\
 h&\longrightarrow\beta(d(h))
       &&(h\notin R_{\mathrm{int}}),\\
 I&\longrightarrow\beta(d(h))
       &&(h\in R_{\mathrm{int}},\ I=\beta(p_1(h))).
\end{align*}
All vertices and incidences on the right are read in the original
network $B$. Equivalently, contract the original tree arcs internal
to the blocks and suppress precisely the reticulations in
$R_{\mathrm{int}}$, retaining and reconnecting their children.
The newly created reconnection arcs are retained.

To check that these arcs are distinct, remove the incoming
reticulation arcs in $B$. Each reticulation starts its own path
component, containing its child. Neither of its parents lies in
that component, since this would create a directed cycle.
Thus the children reconnected to a block lie in original components
distinct from one another and from the component of the block.
Contractions identify only tree vertices in the same original
component, so these children remain distinct. The continuation
$\beta(\sigma(I))$ lies in the component of $I$, and is therefore
distinct from all its reconnected children. The surviving
reticulation children are distinct from these nonreticulation
children, and their two parents are different blocks by the
definition of $R_{\mathrm{int}}$. Hence no parallel arcs are created.

Contracting an arc whose lower endpoint has indegree one cannot
create a directed cycle. Removing a reticulation with a single
distinct parent cannot create one either. Every block keeps its
nonreticulation continuation, so the inverse graph is a normalized
tree-child skeleton. A block with $t$ ports and $r$ internal pairs
has $t-2r$ remaining reticulation children, its continuation, and
the $r$ reconnected children. This proves~\eqref{rvc:blockparameters}.
Since $t\ge1$ and $r\le\lfloor t/2\rfloor$, its outdegree
$c=t-r+1$ is at least two.

The two constructions are inverse, including their local choices.
For the vertex obtained from a block $I$, the continuation is
$\beta(\sigma(I))$; its other nonreticulation children are exactly
the children reconnected from the reticulations internal to $I$.
These children identify the new reticulations in the refinement.
Reading the reticulation children along the directed path of $I$
recovers the multiset order of its ports. Thus both the continuation
and the port order are recovered before counting isomorphism classes.
By Lemma~\ref{rvc:rigidity}, children of a fixed labelled skeleton
cannot be permuted by an automorphism. Hence there is no additional
symmetry divisor in~\eqref{rvc:fibresize}; its factor $2^{b-1}$
only identifies the two copies of each repeated port. The product
therefore counts exactly the fibre of $Q$ over the isomorphism
class of $H$.
\end{proof}

If the blocks are $I$, formula~\eqref{rvc:blockparameters} gives
\begin{equation}
 \delta=\sum_I r(I),\qquad
 e=\sum_I\bigl(t(I)-r(I)-1\bigr),\qquad
 \delta+e=\sum_I(t(I)-1).
 \label{rvc:blockbudget}
\end{equation}
For a block $I$, let $\ell(I)$ count the ordinary children which are
leaves after inversion. These are its continuation, if the block
ends its path, and the leaf children of its internally repeated
reticulations. In particular $0\le\ell(I)\le r(I)+1$. This statistic
depends only on $B$ and the block, not on the partition of any
other path: a reconnected child which is a tree vertex belongs to
a nonempty block and remains a tree vertex after inversion.
Consequently a nonfinal block can have $\ell(I)>0$, and a final
block can have $\ell(I)>1$. If an empty path is reconnected after
suppressing its initial reticulation, its leaf is included in
$\ell(I)$ for the block containing both parents.

\begin{theorem}[Exact maximal-refinement identity]
\label{rvc:exactidentity}
Define
\begin{equation}
 F_{t,r,\ell}(z,y)=
 \frac{2^r z^{t-1}}{(r+1)t!}
       \sum_{j=0}^{\ell}\binom\ell j
                 w_{t-r+1,r+1-j}(zy)^j.
 \label{rvc:fullblock}
\end{equation}
With $\Pi(B)$ as in Proposition~\ref{rvc:refinement}, and with an
empty path contributing the empty product $1$, we have
\begin{equation}
 Z_n(z,y;v)=\frac1{M_n}\sum_{B\in\mathcal M_n}
       \sum_{\pi\in\Pi(B)}
       \prod_{I\in\pi}v^{r(I)}F_{t(I),r(I),\ell(I)}(z,y).
 \label{rvc:exactpartition}
\end{equation}
\end{theorem}

\begin{proof}
For a normalized skeleton $H$, put
\[
 D(H)=\prod_{a\text{ tree in }H}\mathcal N(c(a),b(a))
\]
and let $W(H)$ be its unnormalized weight in
\eqref{rvc:skeletonsum}, namely
\[
 W(H)=v^{\delta(H)}z^{\delta(H)+e(H)}
       \prod_{a\text{ tree in }H}
       \left(\sum_{j=0}^{\ell(a)}\binom{\ell(a)}j
                   w_{c(a),b(a)-j}(zy)^j\right).
\]
Proposition~\ref{rvc:refinement} gives
$|Q^{-1}([H])|=D(H)$, where the fibre consists of pairs $(B,\pi)$
with their block partitions retained.

Fix one such pair. Its blocks correspond to the tree vertices
of $H$, with
\[
 t(I)=c(a)+b(a)-2,\qquad r(I)=b(a)-1,\qquad \ell(I)=\ell(a).
\]
Consequently,
\[
 \frac{2^{r(I)}}{(r(I)+1)t(I)!}
 =\frac1{\mathcal N(c(a),b(a))}.
\]
Together with~\eqref{rvc:blockbudget}, this gives the pointwise
identity
\[
 \prod_{I\in\pi}v^{r(I)}F_{t(I),r(I),\ell(I)}(z,y)
 =\frac{W(H)}{D(H)}.
\]
In particular, the total weight of a pair is constant on its fibre.
Grouping the pairs by their inverse skeleton now yields
\begin{align*}
 &\sum_{B\in\mathcal M_n}\sum_{\pi\in\Pi(B)}
       \prod_{I\in\pi}v^{r(I)}F_{t(I),r(I),\ell(I)}(z,y)\\
 &\qquad=\sum_H\sum_{(B,\pi)\in Q^{-1}([H])}\frac{W(H)}{D(H)}\\
 &\qquad=\sum_H\frac{|Q^{-1}([H])|}{D(H)}W(H)
 =\sum_H W(H)
 =M_nZ_n(z,y;v).
\end{align*}
Dividing by $M_n$ proves~\eqref{rvc:exactpartition}.
\end{proof}

Let $Z_n^{(0)}(z,y)$ be the part of this sum with $\delta=0$.
Equivalently, no block contains a repeated port. The superscript
does not require the inverse skeleton to be binary: contractions
and omitted leaf marks remain allowed. Put
\begin{equation}
 a_t=\frac{(2t-1)!!}{t!}=2^{-t}\binom{2t}{t},\qquad
 \omega_t(z)=t a_tz^{t-1},\qquad
 \tau_t(z,y)=(t+zy)a_tz^{t-1}.
 \label{rvc:referenceweights}
\end{equation}
When $r=0$, formula~\eqref{rvc:fullblock} equals $\omega_t$ for a
nonfinal block and $\tau_t$ for a final block. Define the relaxed
model $\widetilde Z_n^{(0)}(z,y)$ by allowing repeated ports while
retaining these reference weights. Write $G_0(B;z,y)$,
$\widetilde G(B;z,y)$, and $G(B;z,y)$ for the corresponding
partition sums at a fixed $B$. Then their averages over
$\mathcal M_n$ are $Z_n^{(0)}$, $\widetilde Z_n^{(0)}$, and $Z_n$.
When $y=1$ we omit it from the notation.

\subsection{The path kernel and its positive representation}

In the relaxed model a path contributes a weight depending only
on its length. Let $C_j(z,y)$ be this weight for a path of length
$j$, with $C_0(z,y)=1$. The series of nonfinal block weights is
$\Omega(x;z)=x(1-2zx)^{-3/2}$. The final-block series is
$\Omega(x;z)+y((1-2zx)^{-1/2}-1)$. Summing sequences of nonfinal
blocks followed by a final block gives
\begin{equation}
 C(x;z,y):=\sum_{j\ge0}C_j(z,y)x^j
 =\frac{1+y((1-2zx)^{-1/2}-1)}{1-x(1-2zx)^{-3/2}}.
 \label{rvc:bivariatekernel}
\end{equation}
For $y=1$ this simplifies to
\begin{equation}
 C(x;z)=\frac{1-2zx}{(1-2zx)^{3/2}-x},\qquad
 C_1(z)=1+z,\quad C_2(z)=1+4z+\tfrac32z^2.
 \label{rvc:kernel}
\end{equation}
If the tree-path lengths of $B$ are $\ell_0,\ldots,\ell_m$,
\begin{equation}
 \widetilde G(B;z,y)=\prod_{i=0}^m C_{\ell_i}(z,y).
 \label{rvc:pathproduct}
\end{equation}
No independence of these lengths is asserted or needed.

\begin{lemma}\label{rvc:stieltjes}
For $z>0$, let $u\in(0,1)$ solve $u^2+2zu^3=1$, and put
$\rho_z=u^3$, $\lambda_z=\rho_z^{-1}$ and $p_z=2u/(3-u^2)$.
Then
\begin{equation}
 C_j(z)=p_z\lambda_z^j+
 \int_0^{2z}t^j
 \frac{(2z-t)^{5/2}}
 {\pi\sqrt t\,[t+(2z-t)^3]}\,dt.
 \label{rvc:stieltjesformula}
\end{equation}
The measure on the right is positive and has total mass one.
Consequently,
\begin{equation}
 p_z\lambda_z^j\le C_j(z)\le\lambda_z^j,
 \qquad
 \frac{C_2(z)}{C_1(z)}\le\frac{C_{j+1}(z)}{C_j(z)}
       \le\lambda_z\quad(j\ge1).
 \label{rvc:coefficientbounds}
\end{equation}
Moreover $\rho_z<1/(2z)$ is the unique singularity of smallest
modulus of $C(x;z)$, and is a simple pole.
\end{lemma}

\begin{proof}
Consider $H(s)=C(1/s;z)/s$ on the principal sheet.
Its possible poles correspond to roots of $2zu^3+u^2-1=0$
with nonnegative real part. There is exactly one positive root.
The other roots, if real, are negative; if conjugate, their real
part is negative because their sum is $-1/(2z)-u$. Thus outside
the cut $[0,2z]$ the only pole is $\lambda_z$.
Its residue is $p_z$, and on the upper bank of the cut,
\[
 -\pi^{-1}\operatorname{Im}H(t+i0)
 =\frac{(2z-t)^{5/2}}
        {\pi\sqrt t\,[t+(2z-t)^3]}.
\]
This density is integrable at both endpoints. Cauchy's formula,
together with $H(s)\sim1/s$ at infinity, gives
\eqref{rvc:stieltjesformula} and total mass one. The support and
its atom imply the first inequalities in
\eqref{rvc:coefficientbounds}. Cauchy--Schwarz gives strict
log-convexity $C_j^2<C_{j-1}C_{j+1}$, so the coefficient ratios
increase to $\lambda_z$, proving the remaining inequalities.
Finally $\Omega(\rho_z;z)=1$ with positive derivative, whereas
$\rho_z<1/(2z)$. If another point of modulus $\rho_z$ also
satisfied $\Omega(x;z)=1$, equality in the triangle inequality
would force every positive-degree monomial of $\Omega$ to have
the same phase. Since all these degrees occur, that point must
be $\rho_z$ itself.
\end{proof}

\subsection{The exact maximal word average}
\label{rvc:wordencoding}

We give the particular word encoding needed here. It concerns
only maximal binary tree-child networks, as in~\cite{FYZ21}.
Let $\mathcal W_m$ consist of words on $[m]$ in which every letter
occurs three times, the second occurrences appear in the order
$1,\ldots,m$, and the third occurrences appear in the same order.

To encode $B\in\mathcal M_n$, process its root tree path first.
Then process tree paths below reticulations whose two parents have
already been visited, in the order in which their second parents
were encountered. Label the reticulations $1,\ldots,m$ in that
order. Write the port labels along each path, inserting the third
occurrence of $i$ immediately before processing the path below
reticulation $i$. Acyclicity guarantees an available unprocessed
component whenever one remains. One obtains a word in
$\mathcal W_m$, together with the ordered list of its $n$ terminal
leaf labels.

Conversely, cut a word immediately after each third occurrence
and delete these delimiter symbols. The first and second
occurrences remaining in its $m+1$ intervals form the port paths.
The path after the third occurrence of $i$
is placed below reticulation $i$, whose parents are the two port
positions carrying $i$. Both parents occur before that delimiter,
so this reconstruction is acyclic. Each tree vertex has its path
continuation as a nonreticulation child. The two constructions are
inverse, with the leaf labels assigned in arbitrary order to the
path ends. In particular,
\begin{equation}
 M_n=n!\,|\mathcal W_m|.
 \label{rvc:wordcount}
\end{equation}

For later use, erase the letter labels and denote a first, second,
or third occurrence by $F,S,T$, respectively. Every prefix satisfies
$\#F\ge\#S\ge\#T$, and each of the three total counts is $m$.
At a second occurrence, if $i$ first and $j$ second occurrences
have already appeared, there are $i-j$ choices for its unmatched
first occurrence. The next label is fixed by the order of the
second occurrences. The third occurrences require no further
choices. Thus a profile has weight equal to the product of $i-j$
over all its $S$ steps.

For a series $C(x)=\sum_{j\ge0}c_jx^j$ with $c_0=1$, let $A_m(C)$
be the sum of these profile weights multiplied by $c_j$ for each
interval of $j$ consecutive $F,S$ steps before the next $T$.
There are $m$ such intervals. The final symbol of every word is
necessarily a third occurrence, so the last, $(m+1)$st tree
component is empty and contributes $c_0=1$.
Write $C_0(x)=(1-x)^{-1}$ for the zero-activity kernel with all
coefficients equal to one; this subscript indicates the activity,
whereas the index in $C_j(z,y)$ indicates a path length. We have proved
\begin{equation}
 \boxed{\quad
 M_n=n!A_m(C_0),\qquad
 \widetilde Z_n^{(0)}(z,y)
 =\frac{A_m(C(\cdot;z,y))}{A_m(C_0)}.
 \quad}
 \label{rvc:wordaverage}
\end{equation}
This identity identifies the combinatorial input to the analytic
study of weighted words.

\section{A bounded comparison of the full and maximal sectors}
\label{rvc:comparisonsection}

We now show that allowing nonmaximal component skeletons changes
the partition function by only a bounded factor. The proof works
under the weighted law itself and does not assume concentration
properties of unweighted networks.

For $z>0$ define
\begin{align}
 \beta_z&=\lambda_z\frac{C_1(z)}{C_2(z)},&
 \theta_z&=\frac{3z}{1+3z},&
 \mathcal L_z&=\sum_{d\ge1}-\log(1-\theta_z^d),
 \label{rvc:comparisonconstants}\\
 \chi_z&=\max\{1,\tfrac58 e^{4z/3}\}.&&&
 \label{rvc:blockcomparisonconstant}
\end{align}
The series defining $\mathcal L_z$ converges. We will also prove
that $\chi_z$ may be replaced by one when $0<z\le33/32$.

\begin{theorem}\label{rvc:boundedratios}
For every $n\ge2$ and $z>0$,
\begin{align}
 1\le\frac{\widetilde Z_n^{(0)}(z)}{Z_n^{(0)}(z)}
       &\le e^{\beta_z\mathcal L_z},
       \label{rvc:relaxedratio}\\
 \frac{Z_n(z)}{\widetilde Z_n^{(0)}(z)}
       &\le e^{3z\beta_z(\chi_z-1)},
       \label{rvc:fullrelaxedratio}\\
 1\le\frac{Z_n(z)}{Z_n^{(0)}(z)}
       &\le e^{\beta_z\mathcal L_z+3z\beta_z(\chi_z-1)}.
       \label{rvc:fullratio}
\end{align}
These bounds are uniform on compact subsets of $z>0$. In particular,
the logarithms of the three partition functions differ by $O_z(1)$.
\end{theorem}

The proof is developed in the next three subsections. An improved
explicit constant at $z=1$ is given in
Corollary~\ref{rvc:thirty}.

\subsection{Factorial moments under the tilted law}

Write $W_z(B)=\widetilde G(B;z)=\prod_i C_{\ell_i}(z)$ and let
$\widetilde\mu_{n,z}$ be the probability measure on $\mathcal M_n$
proportional to $W_z(B)$. Denote by $\mathcal P(B)$ the
reticulations whose two parents lie in the same tree path. For
$a\in\mathcal P(B)$, let $d_a\ge1$ be the number of tree arcs
between its upper and lower parents.

\begin{lemma}[Weighted factorial moments]
\label{rvc:weightedmoments}
If $f:\mathbb N_{\ge1}\to[0,\infty)$ is nonincreasing and
$S_f=\sum_{d\ge1}f(d)<\infty$, then for every integer $r\ge0$,
\begin{equation}
 \mathbb E_{\widetilde\mu_{n,z}}
   \sum_{\substack{A\subset\mathcal P(B)\\|A|=r}}
       \prod_{a\in A}f(d_a)
 \le\frac{(\beta_zS_f)^r}{r!}.
 \label{rvc:factorialmoment}
\end{equation}
Consequently,
\begin{equation}
 \mathbb E_{\widetilde\mu_{n,z}}\sum_a f(d_a)\le\beta_zS_f,
 \qquad
 \mathbb E_{\widetilde\mu_{n,z}}\prod_a(1+f(d_a))
       \le e^{\beta_zS_f}.
 \label{rvc:momentconsequences}
\end{equation}
\end{lemma}

\begin{proof}
Mark a set $A$ of $r$ reticulations in $\mathcal P(B)$, and write
$u_a,v_a$ for the upper and lower parents of $a$ in the original
network $B$. Their pairs of parents are disjoint, since each tree
vertex has exactly one reticulation child; their intervening
intervals need not be disjoint. Fix an order $a_1,\ldots,a_r$ of $A$.

Move the upper parent of a marked pair to the beginning of the
root path: suppress it on its old tree continuation, retain the
reticulation and its lower parent, and insert a new tree vertex
below the planted root. Its children are this reticulation and
the child of the root \emph{after the suppression}. This also
covers the case in which the upper parent was already the first
root-path vertex: the underlying network may then be unchanged.
The operation preserves the maximal tree-child class and cannot
create a cycle. The path of origin never disappears because the
lower parent survives.

Let $P_0,\ldots,P_m$ be the original sequences of tree vertices
on the paths, with $P_0$ the root path, and put
$U_k=\{u_{a_1},\ldots,u_{a_k}\}$. Here $P_i\setminus U_k$ means
deletion from the sequence, preserving the order of the remaining
vertices. If $x_{a_j}$ is the new root vertex inserted at step $j$,
then after $k$ moves the path sequences are exactly
\begin{align*}
 P_0^{(k)}&=(x_{a_k},x_{a_{k-1}},\ldots,x_{a_1})
                  \mathbin{\Vert}(P_0\setminus U_k),\\
 P_i^{(k)}&=P_i\setminus U_k\qquad(i\ne0),
\end{align*}
where $\Vert$ denotes concatenation. Indeed, the next move removes
an original upper parent still present in its path and adds one
new vertex to the root prefix. No lower parent is moved, and the
two parents of each unprocessed pair remain in their original
path and relative order. This proves the invariant by induction
and shows that every one of the $r!$ orders is valid.

The new vertices form the first $r$ vertices of the output root
path, and their distinct reticulation children record the reverse
order $a_r,\ldots,a_1$. These identities are unambiguous: a
reticulation is identified by the labelled leaf at the end of its
child path, which is unchanged. Different orders of the same
marked set therefore remain distinguishable.

At each move from a nonroot path of length $\ell\ge2$, that length
decreases by one and the root length $\ell_0\ge1$ increases by one.
Thus, by~\eqref{rvc:coefficientbounds},
\[
 \frac{W_z(B_{\mathrm{old}})}{W_z(B_{\mathrm{new}})}
 =\frac{C_\ell/C_{\ell-1}}{C_{\ell_0+1}/C_{\ell_0}}
 \le\beta_z.
\]
The bounds apply with indices $\ell-1\ge1$ and $\ell_0\ge1$.
For a move within the root path all lengths are unchanged and
the ratio is one. The same coefficient bounds give
$C_2/C_1\le\lambda_z$, hence $\beta_z\ge1$, so these moves are
covered as well. Consequently the total distortion of weight
after all $r$ moves, with output denoted by $B'$, satisfies
$W_z(B)\le\beta_z^rW_z(B')$.

Now reverse the moves, in the order $r,r-1,\ldots,1$. At step $k$,
remove the first root-path vertex $x_{a_k}$. Its reticulation
child identifies $a_k$, and the other parent is the surviving
$v_{a_k}$. The unique upward tree path from this vertex determines
the path of origin. A reinsertion position above it is specified
by the resulting distance $d'_k\ge1$ between the two parents;
there is at most one position for each positive integer distance.

Write $(u_a,v_a)_B$ for the tree vertices strictly between the
two parents in their original path. For a true preimage, the
path invariant gives the exact identity
\[
 d'_k=d_{a_k}
       -\#\{i<k:u_{a_i}\in(u_{a_k},v_{a_k})_B\}.
\]
The moves with indices greater than $k$ have already been undone,
and precisely the earlier upper parents are still removed.
Each such parent inside the original interval shortens the
distance by one. If the path of origin is the root path, its
remaining new prefix lies above the true insertion position
and contributes no intervening arc. Thus
$1\le d'_k\le d_{a_k}$, including for crossing or nested intervals,
and monotonicity gives $f(d_{a_k})\le f(d'_k)$.

The output $B'$ together with $(d'_r,\ldots,d'_1)$ recovers the
original network, marked set, and order by these successive
reinsertions. Hence the map
\[
 (B,A,\text{order})\longmapsto(B',d'_r,\ldots,d'_1)
\]
is injective. For a fixed output, its admissible inverse choices
form a tree of depth $r$. At every node, distinct insertion
positions have distinct distances, so
\[
 \sum_{\text{admissible positions}}f(d')\le S_f.
\]
Induction on the depth bounds the total weight of the admissible
inverse lists by $S_f^r$, even though the available positions
may depend on all previous choices.

Summing over all $r!$ orders of each marked set and grouping
by the output network therefore gives
\[
 r!\sum_B W_z(B)
       \sum_{\substack{A\subset\mathcal P(B)\\|A|=r}}
                    \prod_{a\in A}f(d_a)
 \le\beta_z^rS_f^r\sum_{B'}W_z(B').
\]
Division proves~\eqref{rvc:factorialmoment}. Taking $r=1$ and
summing over $r$, respectively, gives
\eqref{rvc:momentconsequences}.
\end{proof}

For example, the number $N_J$ of same-path pairs at distance at
most $J$ satisfies
$\mathbb E_{\widetilde\mu}\binom{N_J}{r}\le(\beta_zJ)^r/r!$.
The summable version of the lemma, rather than an unweighted
estimate for $N_J$, is what yields a bounded comparison.

\subsection{Cuts and a direct repair of the relaxation}

For a fixed $B$, view the relaxed partitions as a probability
measure with weights $\omega_t$ on nonfinal blocks and
$\tau_t(z)=\tau_t(z,1)$ on final blocks. For $a,b\ge1$,
\begin{equation}
 \frac{\omega_{a+b}}{\omega_a\omega_b}\le3z,
 \qquad
 \frac{\tau_{a+b}}{\omega_a\tau_b}\le3z.
 \label{rvc:mergeodds}
\end{equation}
Indeed, after dividing the first ratio by $z$, increasing $a$
multiplies it by
$a(2a+2b+1)/[(a+b)(2a+1)]\le1$; the case $a=b=1$ gives three.
The second ratio is the first multiplied by
$(1+z/(a+b))/(1+z/b)\le1$.

Encode a partition by its cuts on the internal tree arcs.
Conditionally on all other cuts, the odds of not cutting a
specified arc versus cutting it are at most $3z$. This remains
true after requiring any collection of pairs to be separated:
adding a cut preserves all such requirements. It also remains
true if other arcs have previously been fixed to be uncut.
Successive conditioning therefore gives
\begin{equation}
 \mathbb P_B(\text{every arc of }I\text{ is uncut}
       \mid\text{previous pair separations})
 \le\theta_z^{|I|}
 \label{rvc:nocut}
\end{equation}
for any set $I$ of internal tree arcs. No independence assertion
is involved.

Impose the separations of the pairs in $\mathcal P(B)$ one at a
time. The parents of pair $a$ lie in one block precisely when
its $d_a$ intervening arcs are all uncut. Consequently
\begin{equation}
 \frac{G_0(B;z)}{W_z(B)}
 \ge\prod_{a\in\mathcal P(B)}(1-\theta_z^{d_a}).
 \label{rvc:softrepair}
\end{equation}
Average under $\widetilde\mu_{n,z}$ and use Jensen's inequality
with the nonincreasing summable function
$f(d)=-\log(1-\theta_z^d)$. Lemma~\ref{rvc:weightedmoments} gives
\[
 \frac{Z_n^{(0)}(z)}{\widetilde Z_n^{(0)}(z)}
 \ge\mathbb E_{\widetilde\mu}
       \exp\!\left\{-\sum_a f(d_a)\right\}
 \ge e^{-\beta_z\mathcal L_z},
\]
proving~\eqref{rvc:relaxedratio}.

\subsection{The complete model and a pointwise comparison}

For $r\ge1$, use Lemma~\ref{rvc:localbounds} in
\eqref{rvc:fullblock}. With $b=r+1$, $c=t-r+1$,
\[
 \frac{U_{c,b-j}}{U_{c,b}}
 =\prod_{i=0}^{j-1}\frac1{t+1/2-i}.
\]
Since $\ell\le\lfloor t/2\rfloor+1$, the sum of these ratios with
their binomial factors is at most
\[
 \left(1+\frac z{t+3/2-\ell}\right)^\ell\le e^{4z/3}.
\]
Let $F^{\mathrm{rel}}$ be the block's reference weight, either
$\omega_t$ or $\tau_t$. Then
\begin{equation}
 \frac{F_{t,r,\ell}(z,1)}{F^{\mathrm{rel}}}
 \le\frac{2t+1}{2(r+1)t}e^{4z/3}
 \le\frac58e^{4z/3}\le\chi_z,
 \label{rvc:fullblockbound}
\end{equation}
where $t\ge2r\ge2$. Blocks with $r=0$ have exactly their reference
weight.

For a relaxed partition $\pi$, let $B_{\mathrm{rep}}(\pi)$ be
the number of blocks containing an internal pair. Thus
\[
 \frac{G(B;z)}{W_z(B)}
 \le\mathbb E_{B,\mathrm{part}}\chi_z^{B_{\mathrm{rep}}(\pi)}.
\]
Expand this power by selecting repeated blocks, with a factor
$\chi_z-1$ for each selected block. Choosing in addition an
internal pair in every selected block only increases the sum.
Selected pairs occupy different blocks, so their intervals of
intervening arcs are disjoint. For a fixed selected set $A$, the
probability that these pairs lie inside distinct blocks is at
most $\theta_z^{\sum_{a\in A}d_a}$ by~\eqref{rvc:nocut}.
After this bound, dropping the disjointness requirement only
enlarges the sum. Hence
\begin{equation}
 \frac{G(B;z)}{W_z(B)}
 \le\prod_{a\in\mathcal P(B)}
                   (1+(\chi_z-1)\theta_z^{d_a}).
 \label{rvc:fullproduct}
\end{equation}
Apply~\eqref{rvc:momentconsequences} with
$f(d)=(\chi_z-1)\theta_z^d$, noting that
$\sum_{d\ge1}\theta_z^d=3z$. This proves
\eqref{rvc:fullrelaxedratio}. The inclusion
$Z_n^{(0)}\le Z_n$ completes the proof of
Theorem~\ref{rvc:boundedratios}.

We next justify the sharper comparison needed at activity one.

\begin{lemma}\label{rvc:pointwisesqueeze}
For $0<z\le33/32$, every complete block weight is at most its
relaxed reference weight. In particular,
\begin{equation}
 Z_n^{(0)}(z)\le Z_n(z)\le\widetilde Z_n^{(0)}(z).
 \label{rvc:squeeze}
\end{equation}
\end{lemma}

\begin{proof}
For $r=0$ the weights agree. If $r\ge1$,
\begin{equation}
 \frac{F_{t,r,\ell}(z,1)}{\omega_t}
 \le\frac{2t+1}{2(r+1)t}
   \sum_{j=0}^{\ell}\binom\ell j
             \frac{z^j}{(t+1/2)^{\underline j}}.
 \label{rvc:preciseblockratio}
\end{equation}
Here $\ell\le r$ in a nonfinal block, and $\ell\le r+1$ in a
final block; $a^{\underline j}$ denotes a falling factorial.
For $r=1,t=2$, the nonfinal ratio is at most $5/8+z/4\le1$.
The final ratio with $\omega_t$ is at most
$5/8+z/2+z^2/6\le1+z/2=\tau_t/\omega_t$; both comparisons
hold for $z\le3/2$.

For $r=1,t\ge3$, allowing $\ell=2$ gives
\[
 \frac{F}{\omega_t}
 \le\frac12+\frac1{4t}+\frac zt+\frac{z^2}{t(2t-1)}
 \le\frac7{12}+\frac z3+\frac{z^2}{15}.
\]
At $z=33/32$ this is $15329/15360<1$.
For $r=2,t=4$, the bound is
$3/8+z/4+z^2/14+z^3/105<1$; even at $z=21/20$ it equals
$29091/40000$. For $r=2,t\ge5$, bounding the binomial sum by
an exponential yields
\[
 \frac{F}{\omega_t}\le\frac{11}{30}e^{6z/7}
 \le\frac{11}{30}e^{99/112}<\frac{616}{625}<1.
\]
The last inequality follows from $e<3$ and
$e^{-13/112}\le(1+13/112)^{-1}$.
Finally, for $r\ge3$ and $t\ge2r$,
\[
 \frac{F}{\omega_t}\le\frac{13}{48}e^{8z/7}
 \le\frac{13}{48}e^{33/28}<\frac{91}{92}<1,
\]
using $e^{33/28}<3/(1-5/28)=84/23$.
These estimates apply to final blocks as well since
$\tau_t\ge\omega_t$. They exhaust the possibilities.
\end{proof}

\subsection{A certified constant at activity one}

The cut estimate can be sharpened by considering all subdivisions
of an interval simultaneously. Define
\begin{equation}
 U_0(z)=1,\qquad U_t(z)=\sum_{j=1}^t\omega_j(z)U_{t-j}(z),
 \qquad p_d(z)=\frac{\omega_{d+1}(z)}{U_{d+1}(z)}.
 \label{rvc:renewalcut}
\end{equation}
Thus $U_t$ sums compositions of $t$ with all block weights equal
to $\omega$, including the last block.

\begin{lemma}\label{rvc:sharpcut}
The sequence $p_d(z)$ is nonincreasing, satisfies
$p_d(z)\le\theta_z^d$, and may replace $\theta_z^d$ in the
conditional probability that a pair at distance $d$ lies in
one block after any previous pair separations. Consequently,
with
\begin{equation}
 \mathcal L_z^*=\sum_{d\ge1}-\log(1-p_d(z)),
 \label{rvc:sharpL}
\end{equation}
we have
\begin{equation}
 \frac{\widetilde Z_n^{(0)}(z)}{Z_n^{(0)}(z)}
 \le e^{\beta_z\mathcal L_z^*}.
 \label{rvc:sharpratio}
\end{equation}
\end{lemma}

\begin{proof}
Fix every cut outside the $d$-arc interval between a pair. If
leaving the interval uncut violates a previous separation, its
conditional probability is zero. Otherwise every subdivision
of the interval is permitted, because adding cuts preserves
separations. The block obtained without internal cuts may extend
$h$ ports to the left and $k$ to the right of the $d+1$ central
ports. An internal composition $(a_1,\ldots,a_s)$ produces block
lengths $h+a_1,a_2,\ldots,a_s+k$.

With all block weights $\omega$, the ratio of the uncut weight
to the weight of this subdivision cannot increase with $h$ or
$k$, since
$\omega_{t+1}/\omega_t=z(2+1/t)$ decreases with $t$.
If the last block is final, the additional factor in this ratio
is
\[
 \frac{1+z/(d+1+h+k)}{1+z/(a_s+k)}\le1.
\]
Thus the relative weight of each subdivision is at least its
value for the unextended interval with all weights $\omega$.
Summing these subdivisions bounds the conditional uncut
probability by $\omega_{d+1}/U_{d+1}=p_d$.

Increasing the last part injects compositions of $t$ into those
of $t+1$, and multiplies each weight by at least
$\omega_{t+1}/\omega_t$. Hence $U_{t+1}/\omega_{t+1}\ge
U_t/\omega_t$, proving monotonicity of $p_d$. Applying the
single-arc odds bound to the pure $\omega$ composition model
also gives $p_d\le\theta_z^d$.
The series in~\eqref{rvc:sharpL} therefore converges. Repeat
\eqref{rvc:softrepair} and Jensen's inequality with
$f(d)=-\log(1-p_d)$, using Lemma~\ref{rvc:weightedmoments}.
\end{proof}

\begin{corollary}\label{rvc:thirty}
For every $n\ge2$,
\begin{equation}
 \boxed{\quad
 Z_n^{(0)}(1)\le Z_n(1)\le\widetilde Z_n^{(0)}(1)
                 <30Z_n^{(0)}(1).
 \quad}
 \label{rvc:thirtyeq}
\end{equation}
In particular, a uniformly chosen binary RV network on $[n]$
has a maximal normalized skeleton with probability greater
than $1/30$.
\end{corollary}

\begin{proof}
Lemma~\ref{rvc:pointwisesqueeze} and~\eqref{rvc:sharpratio}
reduce the assertion to $e^{\beta_1\mathcal L_1^*}<30$.
Here is a finite rational certificate. The number $\lambda_1$
is the unique root larger than three of
$(\lambda-2)^3-\lambda=0$. This polynomial is increasing there
and has value $461081/125000000>0$ at $1761/500$; hence
\[
 \beta_1=\frac{4\lambda_1}{13}
       <\frac{1761}{1625}<\frac{271}{250}.
\]
All $p_d=p_d(1)$ in~\eqref{rvc:renewalcut} are rational. The
positive logarithm series and $p_d\le(3/4)^d$ give
\begin{align*}
 \mathcal L_1^*&\le
 \sum_{d=1}^{30}\left(
       \sum_{j=1}^{100}\frac{p_d^j}{j}
       +\frac{p_d^{101}}{101(1-p_d)}\right)
       +\frac{(3/4)^{31}}{(1/4)(1-(3/4)^{31})}
       <\frac{627}{200}.
\end{align*}
The last comparison is a finite rational inequality, with all
its entries specified by~\eqref{rvc:renewalcut}. To certify the
exponential as well, put $a=(271/250)(627/200)$. Then
\[
 e^a\le\sum_{j=0}^{20}\frac{a^j}{j!}
       +\frac{a^{21}}{21!}\frac1{1-a/22}<30.
\]
The geometric tail is valid because every subsequent term ratio
is at most $a/22<1$. These two finite inequalities prove the
claimed constant using rational arithmetic alone.
\end{proof}

\subsection{The skeleton defect has an exponential tail}

\begin{theorem}\label{rvc:tightdelta}
Fix $z>0$ and $1<v<\theta_z^{-2}$, and set
$a=\theta_z\sqrt v<1$. Under the RV distribution proportional
to $z^q$,
\begin{equation}
 \mathbb E_zv^\delta=\frac{Z_n(z,1;v)}{Z_n(z)}
 \le D_{z,v}:=
 \exp\left\{\beta_z\mathcal L_z^*
       +\frac{\beta_z\chi_z\sqrt v\,a}{(1-a)^3}\right\}.
 \label{rvc:deltamgf}
\end{equation}
In particular,
\begin{equation}
 \mathbb P_z(\delta\ge h)\le D_{z,v}v^{-h},\qquad
 \mathbb E_z\delta\le\frac{\log D_{z,v}}{\log v}=O_z(1),
 \label{rvc:deltatail}
\end{equation}
and all moments of $\delta$ are bounded uniformly in $n$.
The conclusions are uniform on compact subsets of $z>0$.
\end{theorem}

\begin{proof}
The marked identity~\eqref{rvc:exactpartition} multiplies a block
with $r$ internal pairs by $v^r$. Since $r\le t/2$,
\[
 \frac{v^rF_{t,r,\ell}}{F^{\mathrm{rel}}}
 \le1+\chi_zv^{t/2}\mathbf1_{\{r\ge1\}}.
\]
Expand by selected repeated blocks and select an internal pair
in each one, as in the proof of~\eqref{rvc:fullproduct}.
For a chosen pair at distance $d$, a containing block has length
$t=d+1+h+k$, where $h,k\ge0$ are its extensions. Distinct
selected blocks have disjoint sets of internal arcs. The
probability that specified such blocks appear is at most
$\theta_z^{\sum(t-1)}$, because all their internal arcs must
be uncut; the endpoint-cut requirements only decrease the
probability.

After this estimate, discard all restrictions of disjointness
and of fitting within the paths. The sum over possible blocks
containing one chosen pair is at most
\[
 f_v(d)=\sum_{h,k\ge0}\chi_zv^{(d+1+h+k)/2}\theta_z^{d+h+k}
       =\frac{\chi_z\sqrt v\,a^d}{(1-a)^2}.
\]
This is a nonincreasing summable function, with sum
$\chi_z\sqrt v\,a/(1-a)^3$. The weighted factorial-moment
lemma now yields
\[
 \frac{Z_n(z,1;v)}{\widetilde Z_n^{(0)}(z)}
 \le\exp\!\left\{\frac{\beta_z\chi_z\sqrt v\,a}{(1-a)^3}\right\}.
\]
Divide by $Z_n(z)\ge Z_n^{(0)}(z)\ge
e^{-\beta_z\mathcal L_z^*}\widetilde Z_n^{(0)}(z)$.
Markov's inequality proves the tail bound, and Jensen's inequality
applied to $v^\delta$ proves the mean bound. A single $v>1$ can
be chosen for all $z$ in any fixed positive compact set.
\end{proof}

The theorem concerns the missing reticulations of the normalized
skeleton. It does not imply that $\delta=0$ with probability
tending to one, nor does it bound the total defect $q$ or the
arity excess $e$ by a constant. It also makes no assertion that
the law of $\delta$, or the bounded ratio $Z_n/Z_n^{(0)}$, has
a limit.

\subsection{Uniformity for the leaf-marked kernel}

\begin{proposition}\label{rvc:bivariatecomparison}
The bounded comparisons and exponential-moment bounds above
remain valid for the leaf-marked functions $Z_n(z,y)$ and their
maximal and relaxed sectors, uniformly on compact subsets of
$(z,y)\in(0,\infty)^2$. One may replace $\beta_z$ and $\chi_z$
by
\begin{equation}
 \beta_{z,y}=\max(y,1/y)^2\beta_z,\qquad
 \chi_{z,y}=\max\{1,\tfrac58e^{4zy/3}\},
 \label{rvc:bivariateconstants}
\end{equation}
while keeping $\theta_z$, $p_d(z)$ and $\mathcal L_z^*$ unchanged.
If $zy\le33/32$, one may take $\chi_{z,y}=1$.
\end{proposition}

\begin{proof}
Equation~\eqref{rvc:bivariatekernel} shows that
$C_j(z,y)=A_j(z)+yB_j(z)$ with $A_j,B_j\ge0$ and
$A_j+B_j=C_j(z)$. Therefore
\[
 \min(1,y)C_j(z)\le C_j(z,y)\le\max(1,y)C_j(z).
\]
A nonroot relocation changes four such factors in a ratio,
so its distortion is at most $\max(y,1/y)^2\beta_z$.
A relocation within the root path still leaves the weight
unchanged. This proves the weighted factorial-moment lemma
with $\beta_{z,y}$, without assuming a Stieltjes representation
for the bivariate coefficients.

The terminal merge ratio is multiplied by
$(1+zy/(a+b))/(1+zy/b)\le1$. Thus the cut and sharp-renewal
arguments do not change. In the complete block formula the
optional-mark activity is $zy$, so
\eqref{rvc:fullblockbound} holds with $\chi_{z,y}$.
All comparisons and the proof of the exponential moment now
apply verbatim with these constants. Finally, the case analysis
in Lemma~\ref{rvc:pointwisesqueeze} depends only on the activity
of optional marks after $z^{t-1}$ cancels. It therefore applies
whenever $zy\le33/32$.
\end{proof}

Combining~\eqref{rvc:wordaverage} and
Theorem~\ref{rvc:boundedratios}, or its bivariate version, reduces
the logarithmic asymptotics of the full RV partition function to
the weighted-word problem with an additive error $O(1)$.
In particular, any uniform quantitative estimate for the
weighted-word average transfers to the complete RV model with
that same remainder, up to an additional bounded term.

% END SECTION SOURCE: rv_combinatorics.tex

% BEGIN SECTION SOURCE: tc_leaf_law.tex
\section{Leaves in a maximal tree-child skeleton}\label{sec:leaf-law}

Let $L=L(B)$ be the number of leaves whose parent is a tree vertex in a
uniform network $B\in\mathcal{TC}_{n,n-1}$. This statistic is different
from the corresponding statistic in an expanded reticulation-visible
network. The following estimate supplies the small-defect asymptotics.

\begin{theorem}[The $5/8$ law]\label{tc:leaflaw}
Put $m=n-1$ and let $a_m=M_{m+1}/(m+1)!$ be the number of maximal words
in Section~\ref{rvc:wordencoding}. Define
\begin{equation}\label{tcl:normalizer}
 \mathfrak z_m=\frac{a_m}{12^{m-1}m!}.
\end{equation}
For $m\ge1$ and $t\ge0$,
\begin{equation}\label{tcl:tail}
 \mathbb P\left(\left|L-\frac58m\right|>4+5t\right)
 \le\min\left\{1,\frac5{\mathfrak z_m}e^{-t^2/(4m)}\right\}.
\end{equation}
Consequently
\begin{equation}\label{tcl:mean}
 \mathbb E\left|\frac Ln-\frac58\right|=O(n^{-1/3}),
\end{equation}
and $L/n\to5/8$ in probability and in every fixed moment. For each
$\varepsilon>0$, the probability of a deviation $|L/n-5/8|>\varepsilon$
is at most $\exp\{-c_\varepsilon n+O(n^{1/3})\}$.
\end{theorem}

\begin{proof}
Under the maximal encoding in Section~\ref{rvc:wordencoding},
$L$ is the number of nonempty chains, equivalently the number of runs
of third occurrences. Each word has the same number $(m+1)!$ of leaf
label assignments, so the uniform network measure gives the uniform
word measure for $L$.

\paragraph{Insertion of the largest letter.}
Successively remove the largest letters of a word. If $d_s+1$ is the
length of its final run of third occurrences at size $s$, the resulting
history satisfies
\[
 d_1=0,\qquad 0\le d_s\le d_{s-1}+1.
\]
Conversely, to insert letter $s$ into a word of size $s-1$, place its
third occurrence at the end. Its second occurrence leaves exactly
$d_s$ old third occurrences to its right; this position is fixed once
$d_s$ is specified. The first occurrence has
\[
 N_s=3s-d_s-2
\]
possible positions before the second. The number is independent of
earlier insertion choices. Write $I_s=\mathbf1_{\{d_s\le d_{s-1}\}}$.
Insertion of the second occurrence creates $I_s$ new runs. Of the
$N_s$ first-occurrence positions, precisely
$s-d_s-I_s-L_{s-1}$ split an old run of third occurrences. Therefore,
conditional on the entire history and all previous insertions,
\begin{equation}\label{tcl:increment}
 L_s=L_{s-1}+I_s+B_s,\qquad
 \mathbb P(B_s=1\mid\text{history, previous insertions})
       =\frac{s-d_s-I_s-L_{s-1}}{3s-d_s-2}.
\end{equation}
The initial value is $L_1=1$. Uniform choices among these positions,
conditional on the history, produce the uniform word because their
number at each stage depends only on that history.

\paragraph{A reference walk.}
Under a reference probability $\mathbb P_0$, let
$d_s-d_{s-1}=X_s$ be independent with
\[
 \mathbb P_0(X_s=1)=\frac12,\qquad
 \mathbb P_0(X_s=-j)=2^{-j-2}\quad(j\ge0).
\]
This walk has mean zero and is allowed to become negative. A fixed
nonnegative history has probability $2^{d_m-2(m-1)}$ and represents
$\prod_{s=2}^m(3s-d_s-2)$ words. Hence the history law $\mathbb Q_m$
of a uniform word satisfies the exact identity
\begin{equation}\label{tcl:change-measure}
 \frac{d\mathbb Q_m}{d\mathbb P_0}
 =\frac{W_m}{\mathfrak z_m},\qquad
 W_m=\mathbf1_{\{d_s\ge0\ \forall s\}}2^{-d_m}
       \prod_{s=2}^m\left(1-\frac{d_s+2}{3s}\right).
\end{equation}
In particular $0\le W_m\le1$, $\mathfrak z_m=\mathbb E_0W_m$, and
$\mathbb Q_m(E)\le\mathfrak z_m^{-1}\mathbb P_0(E)$ for every event
depending only on the history. The maximal-count estimate of
\cite[Equation~(18)]{FYZ21}, or the sharper theorem in
Section~\ref{sec:tc-maximum}, gives
\begin{equation}\label{tcl:cost}
 \mathfrak z_m=\Theta\bigl(m^{-2/3}e^{a_1(3m)^{1/3}}\bigr),
 \qquad \log\mathfrak z_m^{-1}=O(m^{1/3}).
\end{equation}

\paragraph{A deterministic centering identity.}
Let $Y_s=L_s-5s/8$ and
$\xi_s=B_s-\mathbb E(B_s\mid\text{history, previous insertions})$.
Equation~\eqref{tcl:increment} gives exactly
\begin{equation}\label{tcl:centered}
 Y_s=\left(1-\frac1{3s}\right)Y_{s-1}
            +(I_s-1/2)+\xi_s+e_s,
 \qquad |e_s|\le\frac{d_s+29/24}{s}.
\end{equation}
Indeed, subtract $(s-L_{s-1})/(3s)$ from the conditional success
probability in \eqref{tcl:increment}. Its numerator, over the common
denominator $3s(3s-d_s-2)$, is
\[
 -d_s(2s+L_{s-1})+2(s-L_{s-1})-3sI_s.
\]
Its absolute value is at most $(d_s+1)/s$ after division; centering
adds $5/(24s)$. Define deterministic weights
\[
 \omega_{s,m}=\prod_{j=s+1}^m\left(1-\frac1{3j}\right).
\]
They belong to $(0,1]$ and satisfy
$\sum_{s=2}^m\omega_{s,m}/s=3(1-\omega_{1,m})\le3$.
With $H_m=\max(0,d_1,\ldots,d_m)$ and $Y_1=3/8$, iteration yields
\begin{equation}\label{tcl:deterministic}
 |Y_m|\le4+3H_m+
 \left|\sum_{s=2}^m\omega_{s,m}(I_s-1/2)\right|
 +\left|\sum_{s=2}^m\omega_{s,m}\xi_s\right|.
\end{equation}

\paragraph{Concentration.}
The reference increment has moment generating function
$\mathbb E_0e^{\theta X_s}=e^{2\theta}/(2e^\theta-1)$.
For $\theta\ge0$ its logarithm has second derivative at most $2$,
and its value and first derivative vanish at zero. The exponential
maximal inequality therefore gives
$\mathbb P_0(H_m>t)\le e^{-t^2/(4m)}$.
Under $\mathbb P_0$, the $I_s$ are independent Bernoulli variables
with parameter $1/2$. The exponential bound for a centered variable
whose range has length one, applied to the deterministic weights,
gives
\[
 \mathbb P_0\left(\left|\sum\omega_{s,m}(I_s-1/2)\right|>t\right)
 \le2e^{-2t^2/m}.
\]
Conditional on any admissible history, the $\xi_s$ are martingale
differences with range of length one, so the same bound holds for
$\sum\omega_{s,m}\xi_s$ under $\mathbb Q_m$, without a change-of-measure
factor. Apply \eqref{tcl:change-measure} to the first two events and
take their union in \eqref{tcl:deterministic}. This proves
\eqref{tcl:tail}.

For an explicit first-moment bound, put
$t_m=\sqrt{4m\log(5/\mathfrak z_m)}$. Integrating the tail gives
\[
 \mathbb E|L-5m/8|
 \le4+5\left(t_m+\frac{2m}{t_m}\right)=O(m^{2/3}),
\]
where the last assertion uses \eqref{tcl:cost}.
Division by $n=m+1$ proves \eqref{tcl:mean}. Boundedness of $L/n$
gives all fixed moments, and taking $t$ proportional to $m$ proves
the last assertion of the theorem.
\end{proof}

The proof also records an exact recurrence that can be used without
network generation. If $w_{s,d,\ell}$ counts maximal words with final
third-occurrence run $d+1$ and $\ell$ such runs, set
$w_{1,0,1}=1$ and all impossible indices to zero. Then
\begin{align}\label{tcl:exact-distribution}
 w_{s,d,\ell}={}&(2s+\ell-2)
 \left(w_{s-1,d-1,\ell}+\sum_{h\ge d}w_{s-1,h,\ell-1}\right)\notag\\
 &+(s-d-\ell+1)
 \left(w_{s-1,d-1,\ell-1}+\sum_{h\ge d}w_{s-1,h,\ell-2}\right).
\end{align}
The two coefficients are the nonsplitting and splitting position
counts in \eqref{tcl:increment}. Thus \eqref{tcl:exact-distribution}
computes the complete distribution, and not only its limiting mean.

% END SECTION SOURCE: tc_leaf_law.tex

% BEGIN SECTION SOURCE: rv_microscopic.tex
% Microscopic RV defect asymptotics.
% Place after rv_combinatorics.tex and the maximal-TC leaf law tc:leaflaw.
% No analytic pressure result is used here.

\section{Microscopic defects and the second defect layer}
\label{rvm:section}

We retain $m=n-1$, $M_n=T_{n,n-1}$,
$V_{n,q}=\operatorname{RV}_{n,3n-3-q}$, and
$Z_n(z)=M_n^{-1}\sum_qV_{n,q}z^q$.
In this section, an expectation without a subscript is uniform
over $\mathcal M_n$, the maximal binary tree-child networks.
For $B\in\mathcal M_n$, let $L=L(B)$ count leaves with a tree
parent, and let $\Delta=\Delta(B)$ count reticulations whose
parents are adjacent on one tree path. The maximal-tree-child
leaf law, Theorem~\ref{tc:leaflaw}, supplies
\begin{equation}
 \mathbb E\left|\frac Ln-\frac58\right|=O(n^{-1/3}),\qquad
 \mathbb P\{|L-5n/8|>x\}
 \le\min\{1,\exp(Cn^{1/3}-cx^2/n)\}
 \label{rvm:leafinput}
\end{equation}
for suitable positive constants $C,c$. Enlarging $C$ and reducing
$c$ absorbs an additive constant in the deviation threshold.
The leaf variable $L$ belongs to the maximal TC reference network;
it is distinct from the RV statistic $s$ in
Lemma~\ref{rvc:budget}.

\subsection{Binary resolutions and comparison of skeleton sectors}

For a normalized vertex with $d$ children, of which $b\ge1$ are
nonreticulation children, the number of binary tree-child
resolutions is
\begin{equation}
 B_{d,b}=\frac{(d+b-2)!}{2^{b-1}(b-1)!}.
 \label{rvm:binaryresolutions}
\end{equation}
To see this, first choose a binary tree on the $b$ ordinary
ports. Insert the other $d-b$ prescribed ports in a fixed order.
There are successively $2b-1,2b,\ldots,d+b-2$ admissible arcs,
including the planted arc of the local resolution and excluding
arcs ending at an already inserted reticulation port. Insertion
is reversed by deleting the last specified port. Multiplying
by $(2b-3)!!$ gives~\eqref{rvm:binaryresolutions}; the convention
for $b=1$ is $(-1)!!=1$.

We need a local bound adapted to these resolutions.

\begin{lemma}\label{rvm:threebound}
For $d\ge2$, $1\le b\le d$ and $0\le j\le b$,
\begin{equation}
 \frac{w_{d,j}}{B_{d,b}}
       \le 3^{d+b-3}(b-1)!.
 \label{rvm:threeboundeq}
\end{equation}
\end{lemma}

\begin{proof}
As in Lemma~\ref{rvc:localbounds}, the nonpositive correction
in~\eqref{rvc:localrecurrence} gives
$w_{d,k}\le(d+k-2)w_{d,k-1}+(k-1)w_{d,k-2}$.
For $k\ge1$ put
\[
 H_{d,k}=3^{d+k-3}\frac{(d+k-2)!}{2^{k-1}}.
\]
The product inequality
$(2d-3)!!\le3^{d-2}(d-2)!$ proves $w_{d,1}\le H_{d,1}$.
The $k=2$ recurrence and $d(d-1)+1\le3d(d-1)/2$
prove $w_{d,2}\le H_{d,2}$. For $k\ge3$, the sum of the
normalized induction factors is
\[
 \frac23+\frac{4(k-1)}{9(d+k-2)(d+k-3)}\le1.
\]
This proves $w_{d,k}\le H_{d,k}$ throughout $1\le k\le d$.
The majorant $H_{d,k}$ increases with $k$, and
$w_{d,0}\le w_{d,1}$. Hence $w_{d,j}\le H_{d,b}$ for $j\le b$.
Dividing by~\eqref{rvm:binaryresolutions} proves the claim.
\end{proof}

Let $\mathcal T_{n,k}$ denote the binary TC networks with $k$
reticulations. A vertex with two nonreticulation children will
be called free. A network in $\mathcal T_{n,m-\delta}$ has exactly
$\delta$ free vertices, by the binary degree balance.

\begin{lemma}\label{rvm:freeedges}
For $0\le\delta\le m$,
\begin{equation}
 T_{n,m-\delta}\le\frac{M_n}{2^\delta\delta!},
 \qquad M_n=2T_{n,m-1}.
 \label{rvm:freebound}
\end{equation}
Moreover,
\begin{equation}
 \mathbb E\Delta=1,\qquad \mathbb E\Delta^2\le2.
 \label{rvm:trianglemoments}
\end{equation}
\end{lemma}

\begin{proof}
Choose one of the $2\delta$ outgoing arcs of a free vertex.
Subdivide it by a new reticulation, and give that reticulation
a second parent by inserting a new tree vertex below the planted
root. Its other child is the previous root child. This is a
tree-child insertion and is injective: the new top vertex and
its reticulation child identify the inverse operation. Thus
$T_{n,k+1}\ge2(m-k)T_{n,k}$, proving the first bound by iteration.
For a maximal output, the top tree vertex necessarily has one
reticulation child, and deleting this insertion always produces
a network with one free vertex. The last step is therefore
surjective as well, proving $M_n=2T_{n,m-1}$.

Delete the arc from the upper parent to the reticulation of a
marked adjacent parent pair in a maximal network. Suppress the
upper parent and the reticulation, both of which now have degree
$(1,1)$. The lower parent becomes the unique free vertex
of a network in $\mathcal T_{n,m-1}$. Conversely, at this free
vertex choose one of its two child arcs, insert a reticulation
there, and insert its other parent immediately above the free
vertex. These two choices reconstruct precisely the marked
adjacent pair. Hence $\sum_{B\in\mathcal M_n}\Delta(B)
=2T_{n,m-1}=M_n$.

For the second moment, let $z\downarrow0$ in
Lemma~\ref{rvc:weightedmoments} with $r=2$ and
$f(d)=\mathbf1_{\{d=1\}}$. Its tilted measure becomes uniform,
and $\beta_z\to1$. Thus
$\mathbb E\binom\Delta2\le1/2$, which together with
$\mathbb E\Delta=1$ gives $\mathbb E\Delta^2\le2$.
\end{proof}

Write $A\preceq B$ for coefficientwise comparison of power
series with real coefficients. For a binary TC network, let
$\mathcal A(B)$ contain the subsets of tree--tree arcs whose
contraction puts no two parents of a reticulation in one block.
These are exactly the admissible contraction sets. Every block
$D$ of $u_D$ tree vertices has $d_D=u_D+1$ children; write
$b_D$ for its nonreticulation children and $\ell_D$ for its leaf
children. Set
\[
 P_{d,b,\ell}(z)=\sum_{j=0}^{\ell}\binom\ell j w_{d,b-j}z^j.
\]
Refining every normalized skeleton vertex according
to~\eqref{rvm:binaryresolutions}, and dividing by the exact fibre
multiplicity, gives
\begin{equation}
 M_nZ_n(z)=\sum_{k=0}^m z^{m-k}
   \sum_{B\in\mathcal T_{n,k}}\sum_{E\in\mathcal A(B)}
       \prod_D\frac{z^{d_D-2}P_{d_D,b_D,\ell_D}(z)}{B_{d_D,b_D}}.
 \label{rvm:binaryidentity}
\end{equation}
The product is over contraction blocks, including singleton
vertices. A block
has a unique incoming arc, so its contraction preserves
acyclicity and the tree-child condition. The only possible
parallel arcs are those to a reticulation whose parents become
identified, precisely the exclusion defining $\mathcal A(B)$.
The local resolution counts multiply by
Lemma~\ref{rvc:rigidity}. These observations justify the exact
fibre count and identity~\eqref{rvm:binaryidentity}.

Let $Z_n^{(\delta)}$ be the contribution from skeleton defect
$\delta$, and for an integer $L$ put
\[
 F_L(z)=(1+3z)^{2m-L}(1+z)^L.
\]

\begin{lemma}\label{rvm:sectorbound}
For every $0\le\delta\le m$,
\begin{equation}
 Z_n^{(\delta)}(z)\preceq(3z/2)^\delta\,\mathbb E F_L(z).
 \label{rvm:sectorboundeq}
\end{equation}
In particular,
\begin{equation}
 Z_n(z)\preceq
       \left(\sum_{\delta=0}^m(3z/2)^\delta\right)
                 \mathbb E F_L(z).
 \label{rvm:allsectorbound}
\end{equation}
\end{lemma}

\begin{proof}
A binary TC network with $n$ leaves has $2m-L$ tree--tree
arcs, irrespective of its reticulation number. For its contracted
blocks, $\sum_D(b_D-1)=\delta$ and
$\sum_D(d_D-2)=|E|$. Hence
\[
 \prod_D(b_D-1)!\le\delta!.
\]
Lemma~\ref{rvm:threebound}, followed by summation over all
contraction sets and optional marks, bounds this network's
contribution to~\eqref{rvm:binaryidentity} by
$3^\delta\delta!z^\delta F_{L(B)}(z)$.

Iterate the insertion in Lemma~\ref{rvm:freeedges} exactly
$\delta$ times. Every source network has
$2^\delta\delta!$ distinct maximal outputs, and the output
recovers the complete insertion history by successive deletion
of its top vertex. Outputs from different sources or histories
therefore cannot coincide. An insertion preserves $L$ or
reduces it by one, the latter occurring when the chosen free
arc ends at a leaf. Moreover $F_L$ is coefficientwise
nonincreasing in $L$, since
\[
 F_{L-1}-F_L
 =2z(1+3z)^{2m-L}(1+z)^{L-1}\succeq0.
\]
Consequently
\[
 \sum_{B\in\mathcal T_{n,m-\delta}}F_{L(B)}
 \preceq\frac1{2^\delta\delta!}
          \sum_{B\in\mathcal M_n}F_{L(B)}.
\]
Combining the two bounds proves~\eqref{rvm:sectorboundeq}.
\end{proof}

\subsection{A locally uniform complex limit}

\begin{theorem}\label{rvm:microlimit}
For every $R>0$, uniformly for complex $t$ with $|t|\le R$,
\begin{equation}
 Z_n(t/n)=e^{19t/4}+O_R(n^{-1/3}).
 \label{rvm:microlimiteq}
\end{equation}
More explicitly, for $n\ge2$ and $n>3R$,
\begin{equation}
 \left|Z_n(t/n)-e^{19t/4}\right|
 \le2Re^{6R}\mathbb E\left|\frac Ln-\frac58\right|
       +\frac{(30R^2+12R)e^{6R}}n.
 \label{rvm:microerror}
\end{equation}
\end{theorem}

\begin{proof}
For $B\in\mathcal M_n$ the number of individually admissible
tree-arc contractions is
$A=2m-L-\Delta$. Consider
\[
 Q_B(z)=(1+3z)^A(1+z)^L.
\]
In the exact maximal-sector path model, a block of $a$ vertices
has nonfinal weight $a a_a z^{a-1}$; a final block has this
weight or $a_a z^a$ according to its leaf mark. Since
\[
 \frac{(a+1)a_{a+1}}{a a_a}=2+1/a\le3,
\]
every configuration is bounded by its weight in $Q_B$.
Thus $G_0(B;z)\preceq Q_B(z)$.

If selected contraction arcs are pairwise nonadjacent and no
selected terminal arc has its leaf mark omitted, the weights
coincide exactly: an isolated contraction contributes
$w_{3,1}/2=3$. Every discrepancy contains either two adjacent
selected arcs or a selected terminal arc and its omitted leaf
mark. There are at most $2n$ pairs of the first kind and $n$
of the second. With $r=|z|$, their product weights are $9r^2$
and $3r^2$, respectively. A positive coefficientwise union
bound therefore gives
\begin{equation}
 |Z_n^{(0)}(z)-\mathbb E Q_B(z)|
       \le21nr^2e^{6nr}.
 \label{rvm:collisionerror}
\end{equation}
Here every remaining independent-site product is bounded by
$e^{(3A+L)r}\le e^{6nr}$. Configurations which contract both
parents of a nonadjacent pair also contain adjacent selected
arcs, so they are included in this bound.

By Lemma~\ref{rvm:sectorbound}, if $3r/2<1$,
\[
 |Z_n(z)-Z_n^{(0)}(z)|
 \le\frac{3r/2}{1-3r/2}\,e^{6nr}.
\]
Replace $A$ in $Q_B$ by $2n-L$. The difference of exponents
is $2+\Delta$, so telescoping the product and using
$\mathbb E\Delta=1$ costs at most $9re^{6nr}$.
Finally replacing each factor $1+u$ by $e^u$ costs at most
$9nr^2e^{6nr}$, by
$|e^u-1-u|\le|u|^2e^{|u|}/2$. We conclude that
\begin{equation}
 \left|Z_n(z)-\mathbb E e^{(6n-2L)z}\right|
 \le\left(30nr^2+9r+\frac{3r/2}{1-3r/2}\right)e^{6nr}.
 \label{rvm:finitecomplexbound}
\end{equation}
For $z=t/n$ and $n>3R$, the last fraction is at most $3R/n$.
Moreover,
\[
 |e^{t(6-2L/n)}-e^{19t/4}|
       \le2Re^{6R}|L/n-5/8|.
\]
This proves~\eqref{rvm:microerror}, and
\eqref{rvm:leafinput} proves~\eqref{rvm:microlimiteq}.
\end{proof}

\begin{corollary}\label{rvm:defecttwo}
For every fixed integer $q\ge0$,
\begin{equation}
 \frac{V_{n,q}}{M_n n^q}
       =\frac{(19/4)^q}{q!}+O_q(n^{-1/3}).
 \label{rvm:fixedcoefficients}
\end{equation}
In particular the leading coefficient of the second defect
layer is
\begin{equation}
 \boxed{\quad V_{n,2}
       =\frac{361}{32}\,n^2M_n\bigl(1+O(n^{-1/3})\bigr).
       \quad}
 \label{rvm:q2}
\end{equation}
\end{corollary}

\begin{proof}
Apply Cauchy's coefficient formula to the locally uniform
complex estimate~\eqref{rvm:microlimiteq}.
\end{proof}

\begin{corollary}\label{rvm:poisson}
Fix $t\ge0$ and give defect $q$ probability
\begin{equation}
 \mathbb P_{n,t}\{Q=q\}
       =\frac{V_{n,q}(t/n)^q}{M_nZ_n(t/n)}.
 \label{rvm:microscopiclaw}
\end{equation}
Its total variation distance from
$\operatorname{Poisson}(19t/4)$ is $O_t(n^{-1/3})$.
\end{corollary}

\begin{proof}
For $t>0$ the probability generating function is
$Z_n(tu/n)/Z_n(t/n)$. Theorem~\ref{rvm:microlimit} approximates
it by $e^{19t(u-1)/4}$ with error $O_t(n^{-1/3})$ uniformly
on $|u|\le2$. Cauchy's formula bounds the difference of the
degree-$q$ coefficients by $O_t(2^{-q}n^{-1/3})$; sum over $q$.
For $t=0$ both distributions are concentrated at zero.
\end{proof}

The activity in this corollary is $t/n$. In particular, it is
not a Poisson approximation for a uniformly chosen RV network,
which corresponds to activity one.

\subsection{Uniform coefficients for growing defects}

\begin{theorem}\label{rvm:growingcoefficients}
If $Q_n=o(n^{1/3})$, then uniformly for integers
$0\le q\le Q_n$,
\begin{equation}
 V_{n,q}=M_n\frac{((19/4)n)^q}{q!}\,(1+o(1)).
 \label{rvm:growingcoeffeq}
\end{equation}
The layers are strictly increasing throughout these ranges
for all sufficiently large $n$.
\end{theorem}

\begin{proof}
Put
\[
 W=6m-2L,\qquad W'=W-3\Delta=L+3A.
\]
The last component of a maximal word is empty, so $1\le L\le m$.
Also $\Delta\le m$. Hence $W\ge4m$ and $W'\ge m$.

Sample $q$ sites independently with replacement, giving each
admissible contraction arc weight three and each leaf-mark
omission weight one, for total weight $W'$. The probability
that two chosen sites coincide is at most $3/W'$. The total
weighted count of incompatible unordered site pairs is at
most $21n$, by the same adjacent-arc and terminal-mark count
used in~\eqref{rvm:collisionerror}. Thus the probability of
some repeated or incompatible pair of choices is at most
\begin{equation}
 \binom q2\left(\frac3{W'}+\frac{42n}{(W')^2}\right)
       =O(q^2/n),
 \label{rvm:birthdaybound}
\end{equation}
uniformly in $B$. All remaining choices have exactly their
true maximal-sector weights. Conversely the exact model is
coefficientwise bounded by $Q_B$. Therefore, for $q=o(\sqrt n)$,
\begin{equation}
 [z^q]Z_n^{(0)}(z)
       =\frac{\mathbb E(W')^q}{q!}\bigl(1+O(q^2/n)\bigr).
 \label{rvm:sectorcoefficient}
\end{equation}

The coefficient contribution from $\delta\ge1$ is bounded by
Lemma~\ref{rvm:sectorbound} and $F_L(z)\preceq e^{Wz}$:
\begin{equation}
 0\le[z^q](Z_n-Z_n^{(0)})
 \le\sum_{\delta=1}^q(3/2)^\delta
             \frac{\mathbb E W^{q-\delta}}{(q-\delta)!}
 =\frac{\mathbb E W^q}{q!}\,O(q/n)
 \label{rvm:residualcoefficients}
\end{equation}
for $q=o(n)$. The last bound follows pointwise from $W\ge4m$,
by summing a geometric series with ratio at most $3q/(8m)$.

It remains to control the powers of $W$. The tail estimate
in~\eqref{rvm:leafinput} implies that $|L-5n/8|$ is
stochastically bounded by $b_n+Y_n$, where $b_n=O(n^{2/3})$
and $\mathbb P(Y_n>x)\le e^{-cx^2/n}$. Indeed, choose
$b_n$ large enough to absorb $Cn^{1/3}$ in the exponent.
For $a\ge0$, integration of this tail yields
\[
 \mathbb E e^{aY_n}
 \le1+O\bigl(a\sqrt n\,e^{O(a^2n)}\bigr).
\]
Consequently, uniformly for $q\le Q_n=o(n^{1/3})$ and any
fixed $D>0$,
\begin{equation}
 \mathbb E\exp\!\left\{
       \frac{Dq(|L-5n/8|+1)}n\right\}=1+o(1).
 \label{rvm:exponentialleafcontrol}
\end{equation}
The deterministic factor is $e^{O(qn^{-1/3})}$, while the
remaining term is $1+O((q/\sqrt n)e^{O(q^2/n)})$.
Since $W=6n-6-2L$ stays between fixed positive multiples of
$n$,~\eqref{rvm:exponentialleafcontrol} implies
\begin{equation}
 \mathbb E W^q=((19/4)n)^q(1+o(1)),
 \label{rvm:Wmoments}
\end{equation}
also with $q$ replaced by $2q$.

For $q\ge1$, Cauchy--Schwarz and
Lemma~\ref{rvm:freeedges} give
\[
 0\le\mathbb E W^q-\mathbb E(W')^q
 \le3q\mathbb E[\Delta W^{q-1}]
 \le3q\sqrt{2\mathbb E W^{2q-2}}
       =\mathbb E W^q\,O(q/n)
\]
in the stated range. Combine this with
\eqref{rvm:sectorcoefficient},~\eqref{rvm:residualcoefficients},
and~\eqref{rvm:Wmoments}. The case $q=0$ is exact.
Finally the ratio of successive layers is asymptotic,
uniformly within the range, to $19n/[4(q+1)]$, which tends
to infinity. This proves the last assertion.
\end{proof}

The microscopic regime does not extend to $q$ proportional
to $n$. Its Poisson law and its initial increase of coefficients
therefore do not assert global unimodality or identify a mode
at fixed positive activity.

% END SECTION SOURCE: rv_microscopic.tex

% BEGIN SECTION SOURCE: rv_analysis.tex
% Analytic section of Counting tree-child and reticulation visible networks.
% No preamble. Environments theorem, lemma, proposition, corollary and proof
% are supplied by the main document. All labels in this file start with rva:.

\section{A quantitative Laplace principle for weighted words}
\label{rva:laplace}

The word representation of maximal tree-child networks
\cite{FYZ21}, in the form established in
\eqref{rvc:wordaverage}, reduces the analytic problem to a family of
positive kernels. We first treat that family independently of the
network interpretation.

Let
\begin{equation}
 C(x)=\sum_{j\ge0}c_jx^j,\qquad c_0=1,\quad c_j>0,
 \qquad
 c_j=A\rho^{-j}+O(\sigma^{-j}),\quad A>0,\quad 0<\rho<\sigma.
 \label{rva:kernelassumption}
\end{equation}
For uniform statements, the kernels range over a compact parameter
set $\mathcal K$, their coefficients depend continuously on the
parameter, the error in~\eqref{rva:kernelassumption} is uniform, and
$\rho$, $A$, and $\sigma/\rho-1$ have positive uniform lower bounds;
$\rho$ and $A$ also have finite uniform upper bounds. Consequently,
there are common constants $d,D>0$ such that
\begin{equation}
 d\rho^{-j}\le c_j\le D\rho^{-j}\qquad(j\ge0).
 \label{rva:geometriccoefficients}
\end{equation}

A type word has $m$ letters of each of the types $F,S,T$ and satisfies
$i\ge j\ge k$ at every prefix, where $i,j,k$ count the three types.
Each letter $S$ has weight $i-j$ immediately before it is inserted.
Each block of $\ell$ letters of types $F,S$ immediately preceding a
$T$ has weight $c_\ell$. The last letter is necessarily $T$, so there
are exactly $m$ such blocks, including empty ones. Let $A_m(C)$ be
the sum of the resulting weights over these type words.

The function
\[
 \nu_C(x)=\frac{xC'(x)}{C(x)},\qquad 0<x<\rho,
\]
is strictly increasing from zero to infinity. Indeed,
$x\nu_C'(x)$ is the variance of the probability distribution
$c_jx^j/C(x)$. Let $\xi_C$ be its unique point of mean one, and define
\begin{equation}
 \xi_CC'(\xi_C)=C(\xi_C),\qquad
 \Lambda_C^{-1}
 =\frac12\left\{\int_{\xi_C}^{\rho}\frac{dx}{C'(x)}
                   +\frac{\xi_C}{C'(\xi_C)}\right\}.
 \label{rva:Lambda}
\end{equation}

\begin{theorem}[Quantitative word asymptotics]
\label{rva:quantlaplace}
Under~\eqref{rva:kernelassumption},
\begin{equation}
 \log A_m(C)=\log(m!)+m\log\Lambda_C
       +O\!\left(m^{2/3}\{\log(m+2)\}^{2/3}\right).
 \label{rva:wordasymptotic}
\end{equation}
The error is uniform for the compact families described above.
In particular, its implied constant and the threshold in $m$ can
be chosen uniformly on $\mathcal K$.
\end{theorem}

We prove the theorem by solving a continuous variational problem
and then comparing its action quantitatively with the discrete sum.
The estimates give a sublinear logarithmic error; they do not give
a multiplicative asymptotic prefactor.

\subsection{The variational problem and its explicit maximizer}

For $v\ge0$ set
\begin{equation}
 H_C(v)=\inf_{0<x<\rho}\{\log C(x)-v\log x\},\qquad
 F_C(r,h)=hH_C((r+1)/h)+(r+1)\log(r+1)-r\log r.
 \label{rva:entropyfunctions}
\end{equation}
We use $0\log0=0$ and the continuous extension
$hH_C((r+1)/h)=-(r+1)\log\rho$ at $h=0$. This extension follows
from $H_C(v)=-v\log\rho+O(1+\log(1+v))$ as $v\to\infty$.

Time $u$ records the proportion of letters $S$ already inserted.
The two scaled heights are $a=(F-S)/m$ and $b=(S-T)/m$.
Let $\mathcal D$ consist of absolutely continuous pairs on $[0,1]$
satisfying
\begin{equation}
 a(0)=a(1)=b(0)=b(1)=0,\quad a,b\ge0,\quad
 r=1+a'\ge0,\quad h=1-b'\ge0.
 \label{rva:admissiblepaths}
\end{equation}
Then $a(u)\le1-u$, $b(u)\le u$, and $\int_0^1r=\int_0^1h=1$.
The action is
\begin{equation}
 \mathcal J_C(a,b)=\int_0^1
       \{F_C(1+a',1-b')+\log a\}\,du.
 \label{rva:action}
\end{equation}
The integral of $F_C$ is finite: the coefficient comparison bounds
its absolute value by a constant times $1+r+h$.
The action is interpreted as $-\infty$ if $\int\log a=-\infty$.

An equivalent expression for the constant in~\eqref{rva:Lambda} is
\begin{equation}
 \Lambda_C^{-1}
 =\int_{\xi_C}^{\rho}\frac{C(x)C''(x)}{C'(x)^3}\,dx.
 \label{rva:Lambdaalternative}
\end{equation}
To verify it, differentiate $C/C'^2$, use
$C(\xi_C)=\xi_CC'(\xi_C)$, and note that $C/C'^2\to0$ at the pole.
Define a decreasing parametrization $x=x(u)$ and a path by
\begin{equation}
 \begin{split}
 u&=\Lambda_C\int_x^\rho\frac{CC''}{C'^3}\,dt,
 \qquad X(x)=x-\frac{C(x)}{C'(x)},\\
 a_*(u)&=\Lambda_C\frac{X(x(u))}{C'(x(u))},
 \qquad b_*(u)=0.
 \end{split}
 \label{rva:optimalpath}
\end{equation}
Here $X(\xi_C)=0$ and $X'=CC''/C'^2>0$. Direct differentiation gives
\begin{equation}
 a_*'=\nu_C(x)-2,\qquad r_*:=1+a_*'=\nu_C(x)-1>0,
 \qquad h_*=1.
 \label{rva:optimalrates}
\end{equation}
Both endpoints of $a_*$ are zero, so this path belongs to $\mathcal D$.
The pole expansion and the nondegeneracy of the critical point give
\begin{equation}
 \begin{array}{ll}
 a_*(u)\asymp u^{2/3},\quad r_*(u)=O(u^{-1/3}),&u\downarrow0,\\
 a_*(u)=1-u+O((1-u)^2),\quad r_*(u)=O(1-u),&u\uparrow1.
 \end{array}
 \label{rva:profileends}
\end{equation}
All these estimates are uniform on $\mathcal K$. For example,
$C(x)\sim A\rho/(\rho-x)$ implies
$u\asymp(\rho-x)^3$ and $a_*\asymp(\rho-x)^2$ at the first
endpoint. At the other endpoint, $\nu_C'(\xi_C)>0$ and
the integrand in~\eqref{rva:optimalpath} is positive and finite.
The coefficient comparison keeps $\xi_C/\rho$ uniformly away
from both zero and one, which gives the stated uniformity.

\begin{proposition}[Variational maximum]
\label{rva:variationalmaximum}
The path~\eqref{rva:optimalpath} is the unique maximizer of
$\mathcal J_C$ on $\mathcal D$, and
\begin{equation}
 \max_{(a,b)\in\mathcal D}\mathcal J_C(a,b)=\log\Lambda_C-1.
 \label{rva:actionvalue}
\end{equation}
\end{proposition}

\begin{proof}
The function $H_C$ is strictly concave on $(0,\infty)$, with
$H_C'(v)=-\log x(v)$, where $\nu_C(x(v))=v$.
Its perspective in~\eqref{rva:entropyfunctions} is jointly concave.
The remaining term $(r+1)\log(r+1)-r\log r$ is strictly concave
for $r>0$. On the candidate path the two derivatives are
\begin{equation}
 \begin{aligned}
 P(u)&:=F_{C,r}(r_*,1)=-\log X(x(u)),
     & P'(u)&=a_*(u)^{-1},\\
 E(u)&:=F_{C,h}(r_*,1)=\log C(x(u)),
     & E'(u)&<0.
 \end{aligned}
 \label{rva:calibration}
\end{equation}
For any admissible path, tangent inequalities therefore give
\[
 F_C(r,h)-F_C(r_*,1)\le P(r-r_*)+E(h-1),\qquad
 \log a-\log a_*\le(a-a_*)/a_*.
\]
Integrating by parts cancels the terms involving $a-a_*$ because
$r-r_*=(a-a_*)'$ and $P'=1/a_*$. The remaining term is
\begin{equation}
 \int_0^1E(h-1)\,du=\int_0^1E'(u)b(u)\,du\le0.
 \label{rva:barrierterm}
\end{equation}
These integrations may first be performed on
$[\epsilon,1-\epsilon]$. To justify the limit for an arbitrary
admissible pair, the constraints~\eqref{rva:admissiblepaths} give
\[
 \begin{aligned}
 \int_u^1r(t)\,dt&=1-u-a(u)\le1-u,\\
 \int_0^u h(t)\,dt&=u-b(u)\le u.
 \end{aligned}
\]
Tonelli's theorem therefore yields
\[
 \int_0^1r(u)\log\frac1{1-u}\,du\le1,\qquad
 \int_0^1h(u)\log\frac1u\,du\le1.
\]
The coefficient comparison and~\eqref{rva:profileends} also give
$a_*(u)\asymp u^{2/3}(1-u)$ uniformly, so that
$a(u)/a_*(u)=O(u^{-2/3})$. Differentiating $E$ shows that
$E'=O(u^{-1})$ near zero and is bounded near one.
These bounds ensure the absolute integrability of $Pr$, $Eh$,
$P'(a-a_*)$, and $E'b$, using $b\le u$ at zero.
Their boundary terms vanish:
$P$ is bounded at zero and is $O(|\log(1-u)|)$ at one, where
$a,a_*\le1-u$; $E=O(|\log u|)$ at zero, where $b\le u$, and
$E$ is bounded at one. Paths of action $-\infty$ cause no
exception. Thus the candidate is a global maximizer. Equality
requires $b=0$ in~\eqref{rva:barrierterm} and $a=a_*$ in the
strict tangent inequality for the logarithm, proving uniqueness.

For its value, $r_*+1=xC'/C$ yields
\[
 F_C(r_*,1)=\log C'(x)-r_*\log X(x),\qquad
 F_C(r_*,1)+\log a_*=
       \log\Lambda_C-a_*'\log X(x).
\]
Since $(\log X(x(u)))'=-1/a_*$, integration by parts gives
$\int_0^1a_*'\log X=1$, with zero boundary terms by
\eqref{rva:profileends}. This proves~\eqref{rva:actionvalue}.
\end{proof}

The same calculation gives a uniform estimate for regularized paths.
Write $\mathcal J_{C,\eta}$ for the action with $\log(a+\eta)$
in place of $\log a$, where $0<\eta\le1/2$. Tangency at
$a_*+\eta$ leaves, after integration by parts, the residual
\[
 -\frac{\eta(a-a_*)}{a_*(a_*+\eta)}\le\frac\eta{a_*+\eta}.
\]
Consequently, for every path in $\mathcal D$,
\begin{equation}
 \mathcal J_{C,\eta}(a,b)
 \le\log\Lambda_C-1+
 \int_0^1\left\{\log(1+\eta/a_*)+
                     \frac\eta{a_*+\eta}\right\}du
 \le\log\Lambda_C-1+O_{\mathcal K}(\eta(1+|\log\eta|)).
 \label{rva:regularization}
\end{equation}
The last estimate follows from~\eqref{rva:profileends}; its
logarithmic factor comes from the endpoint $u=1$.

\subsection{A uniform local coefficient estimate}

\begin{lemma}
\label{rva:coefficientlemma}
Fix $v_0>0$. There are constants $c_0,c_1>0$ and $N_0$, uniform
on $\mathcal K$, such that, if $N\ge N_0$, $L$ is an integer,
$v=L/N\ge v_0$, and $\log v\le c_0N$, then
\begin{equation}
 \frac{c_1}{v\sqrt N}e^{NH_C(v)}
 \le[x^L]C(x)^N\le e^{NH_C(v)}.
 \label{rva:coefficientbound}
\end{equation}
\end{lemma}

\begin{proof}
The upper bound is Cauchy's positive-coefficient bound. For the
lower bound choose $t\in(0,\rho)$ with $\nu_C(t)=v$ and let
$\pi_j=c_jt^j/C(t)$. With $\tau=t/\rho$,
\begin{equation}
 \frac dD(1-\tau)\tau^j\le\pi_j
      \le\frac Dd(1-\tau)\tau^j.
 \label{rva:geometriclaw}
\end{equation}
The mean condition gives $\tau\ge\tau_0>0$ and
$v\asymp\tau/(1-\tau)$. For $X$ with this law,
\[
 \operatorname{Var}X\asymp v^2,
 \qquad \mathbb E|X-v|^3=O(v^3),
\]
uniformly. For the lower variance estimate, apply the lower bound
in~\eqref{rva:geometriclaw} to two independent copies in
$2\operatorname{Var}X=\mathbb E(X-X')^2$.

Let $\phi(\theta)=\mathbb E e^{i\theta X}$, and let $G$ be a
geometric variable with parameter $\tau$. Since $\pi$ dominates
$(d/D)$ times its law and $1-\cos(\theta(X-X'))\ge0$,
\[
 1-|\phi(\theta)|^2\ge(d/D)^2
       (1-|\phi_G(\theta)|^2)
 \ge c\min(v^2\theta^2,1),\qquad |\theta|\le\pi.
\]
Here the last inequality follows directly from
\[
 1-|\phi_G(\theta)|^2
 =\frac{2\tau(1-\cos\theta)}
        {(1-\tau)^2+2\tau(1-\cos\theta)}.
\]
Thus $|\phi(\theta)|\le e^{-c\min(v^2\theta^2,1)}$.
On $|\theta|\le A/(v\sqrt N)$, the third-moment bound gives
\[
 (e^{-iv\theta}\phi(\theta))^N
 =e^{-N\operatorname{Var}(X)\theta^2/2}+O_A(N^{-1/2}).
\]
The central Fourier integral is therefore bounded below by a
positive constant times $1/(v\sqrt N)$. Choose $A$ large enough
that the Gaussian bound on
$A/(v\sqrt N)<|\theta|\le1/v$ makes the absolute value of this
part less than one quarter of the central lower bound.
The remaining part is $O(e^{-cN})$, also negligible compared with
$1/(v\sqrt N)$ when $\log v\le c_0N$ and $c_0$ is sufficiently
small. Fourier inversion at the integer $L=Nv$ now gives
$\mathbb P(X_1+\cdots+X_N=L)\ge c_1/(v\sqrt N)$.
Finally,
\[
 [x^L]C(x)^N=C(t)^Nt^{-L}
       \mathbb P(X_1+\cdots+X_N=L),
\]
which proves the result.
\end{proof}

Two elementary estimates will handle partial blocks and integer
counts. Set $G(L,T)=TH_C(L/T)$, with $G(L,0)=-L\log\rho$.
Then
\begin{equation}
 0\le G(L,T+1)-G(L,T)\le K+\log(L+1).
 \label{rva:extrafactor}
\end{equation}
Indeed, $\partial_TG=\log C(t(L/T))\le K+\log(1+L/T)$ by
\eqref{rva:geometriclaw}, and this bound can be integrated on
$[T,T+1]$, including $T=0$. Also, if
$B(R,s)=(R+s)\log(R+s)-R\log R-s\log s$, then
\begin{equation}
 B(R,s)-\log(R+s+1)\le\log\binom{R+s}{s}\le B(R,s).
 \label{rva:binomialentropy}
\end{equation}
When $s,T\asymp m/N$, changing $R,s,T$ by bounded amounts changes
$G(R+s,T)+B(R,s)$ by $O(\log(m+1))$. This remains true at $R=0$:
the logarithmic singularity of $\partial_RB=\log(1+s/R)$ is
integrable on an interval of bounded length.

\subsection{Quantitative passage from words to paths}

For a word interval containing $R$ letters $F$, $s$ letters $S$,
and $T$ letters $T$, assign weight $c_j$ to each block of length
$j$, including an auxiliary weight for its terminal partial
block. The matching factors associated with the $S$ steps are
not included here and will be bounded separately. The unrestricted
sum of these block weights is exactly
\begin{equation}
 \binom{R+s}{s}[x^{R+s}]C(x)^{T+1}.
 \label{rva:slabidentity}
\end{equation}
A weak composition of $R+s$ into $T+1$ parts determines the block
lengths, and the binomial factor assigns the types $F,S$ to the
ordered ports. Moreover, \eqref{rva:geometriccoefficients} gives
\[
 \frac d{D^2}\le\frac{c_{a+b}}{c_ac_b}\le\frac D{d^2}
 \qquad(a,b\ge0).
\]
Each cut splits at most one block. This comparison accounts for
the auxiliary partial-block weights when the intervals are
concatenated, even if several cuts split the same original block.
Cutting a word into $N$ intervals therefore changes its block
weight by at most factors $e^{O(N)}$.

\paragraph{Upper bound.}
Cut after each group of approximately $m/N$ letters $S$, including
the final tail of $T$'s in the last interval. The possible vectors
of counts $(R_i,s_i,T_i)$ number at most $(m+1)^{2N}$.
Let $a_N,b_N$ be the linear interpolations of the scaled heights
at these cuts, with time equal to the cumulative number of $S$'s
divided by $m$. They belong to $\mathcal D$. Within interval $i$,
each $S$ has weight at most $m\alpha_i$, where
$\alpha_i=a_{i-1}+R_i/m$. Furthermore,
\[
 0\le\alpha_i-a_N(u)\le(R_i+s_i)/m.
\]
Since each time interval has length at most $1/N+1/m$ and
$\sum R_i=\sum s_i=m$,
\[
 \sum_i\frac{s_i}{m}\log\alpha_i
 \le\int_0^1\log(a_N+\eta)\,du
       +O\!\left(\frac1{N\eta}+\frac1{m\eta}\right).
\]
This estimate does not require a bound on the number of $F$'s
in an individual interval. Use~\eqref{rva:slabidentity}, the
upper coefficient bound, \eqref{rva:extrafactor}, and
\eqref{rva:binomialentropy}, and then sum over count vectors.
Writing $J_*:=\log\Lambda_C-1$, the regularized variational bound
\eqref{rva:regularization} gives
\begin{equation}
 \log A_m(C)-m\log m-mJ_*
 \le K\left\{m\eta(1+|\log\eta|)
        +\frac m{N\eta}+N\log(m+1)+\eta^{-1}\right\}.
 \label{rva:discreteupper}
\end{equation}

\paragraph{An interior lower bound.}
Keep $a_*$ and introduce a positive second margin
$b_\epsilon(u)=\epsilon u(1-u)$, so that
$h_\epsilon=1-\epsilon(1-2u)$. For $0<\epsilon\le1/2$,
\begin{equation}
 \left|\int_0^1F_C(r_*,h_\epsilon)\,du
          -\int_0^1F_C(r_*,1)\,du\right|=O(\epsilon).
 \label{rva:marginloss}
\end{equation}
This follows from
$F_{C,h}=\log C(t((r+1)/h))=O(1+\log(1+r))$ and the
integrability of $\log(1+r_*)$ at zero.

Remove $[0,\delta]$ and $[1-\delta,1]$, and divide the remaining
interval into $N$ intervals of approximately equal $S$-time.
Round the cumulative counts of $F,S,T$, preserving their
monotonicity. Local counts are then nonnegative and differ from
the continuous ones by bounded amounts. In this interior region,
$a_*\ge c\delta$ and $b_\epsilon\ge c\epsilon\delta$.
Provided
\begin{equation}
 N\epsilon\delta\ge D_0,
 \qquad N\delta\longrightarrow\infty,
 \qquad m\epsilon\delta\longrightarrow\infty,
 \label{rva:marginconditions}
\end{equation}
with $D_0$ sufficiently large, every local ordering of the selected
counts is admissible. A decrease in $F-S$ is at most $s_i$, and
a decrease in $S-T$ is at most $T_i$; both are smaller than the
respective initial margins. The rounded errors do not affect this
conclusion.

Require each interior interval to end in $S$, and let $S_i$ be
the prescribed cumulative number of $S$'s at its right endpoint.
The cut is recovered from the concatenated word as the position
immediately after its $S_i$th $S$. The initial endpoint word below
also ends in $S$, so every cut is recovered in this way and
concatenation is injective.

Deleting the final $S$ is a bijection from the local words ending
in $S$ to the unrestricted local words with counts
$(R_i,s_i-1,T_i)$. If the latter word has terminal partial-block
length $j$, appending the deleted $S$ multiplies its block weight
by $c_{j+1}/c_j$. By~\eqref{rva:geometriccoefficients} and the
uniform bounds on $\rho$, this ratio is bounded above and below
by positive constants. The sum of block weights for the words
ending in $S$ is therefore bounded below, up to a constant factor, by
\[
 \binom{R_i+s_i-1}{s_i-1}
       [x^{R_i+s_i-1}]C(x)^{T_i+1}.
\]
Here $T_i+1\asymp s_i\asymp m/N$ and
$v_i=(R_i+s_i-1)/(T_i+1)\ge v_0>0$; also $v_i\le Km$.
If $m/N\gg\log m$, Lemma~\ref{rva:coefficientlemma} applies
uniformly. Its logarithmic factors, the binomial estimates,
the additional factor of $C$, and the rounded counts together
cost $O(N\log(m+1))$.

Concavity of $F_C$ shows that replacing each pair of rates by its
interval average does not decrease its integrated contribution.
Every matching choice has weight at least
$m(a_{i-1}-s_i/m+O(1/m))$. Relative to $\int\log a_*$, the
integrated loss is
$O(1/(N\delta)+1/(m\delta))$: use the lower bound $c\delta$
for the heights, bound their change by
$(R_i+s_i)/m+O(1/m)$, and sum the counts. Multiplication by $m$
gives $O(m/(N\delta)+\delta^{-1})$. Replacing unrounded interval
lengths by the exact integer counts $s_i$ adds
$O(1+|\log\delta|)$ per interval, already included in
$O(N\log(m+1))$.

\paragraph{The two endpoints.}
At the initial endpoint take the single word
$F^RS^{s-1}T^TS$ with the rounded counts at $u=\delta$.
It is admissible because
$R\asymp m\delta^{2/3}\gg s\asymp m\delta$ and
$s-T\asymp m\epsilon\delta\to\infty$.
It ends in $S$, as required for canonical cuts.
For all sufficiently large $m$ we have $T\ge1$, and its block
weight, including the terminal partial block consisting of the
last $S$, is $c_{R+s-1}c_1$; the intervening empty blocks have
weight $c_0=1$. The logarithm of its block weight is therefore
$-(R+s)\log\rho+O(1)$; its matching product is $(R)_s$.
After removing $m^s$, the latter costs
$O(m\delta(1+|\log\delta|))$ in the logarithm. The action on
this interval differs from $-(R+s)\log\rho$ by the same order,
since, uniformly for $h\in[1/2,2]$,
\[
 F_C(r,h)=-(r+1)\log\rho+O(1+\log(1+r)).
\]

At the final endpoint put all remaining $F$'s first, then all
remaining $S$'s, then all remaining $T$'s. This is admissible.
After the $F$'s there are exactly $s$ unmatched first occurrences,
so the matching product is exactly $s!$; the final $T$'s exhaust
the second margin. Stirling's formula and
$a_*(u)=1-u+O((1-u)^2)$ give a loss
$O(m\delta(1+|\log\delta|)+\log m)$.
Block cuts contribute only $e^{O(N)}$ as before.

Combining the interior and endpoint bounds yields
\begin{equation}
 \begin{split}
 \log A_m(C)-m\log m-mJ_*
 \ge-K\biggl\{&m\epsilon+m\delta(1+|\log\delta|)
     +\frac m{N\delta}+N\log(m+1)\\
     &+\delta^{-1}+\log m\biggr\}.
 \end{split}
 \label{rva:discretelower}
\end{equation}

\paragraph{Choice of parameters.}
Let $L_m=\log(m+2)$ and set
\begin{equation}
 \eta=\delta=m^{-1/3}L_m^{-1/3},\qquad
 \epsilon=m^{-1/3}L_m^{2/3},\qquad
 N=\left\lceil D_0m^{2/3}L_m^{-1/3}\right\rceil.
 \label{rva:meshparameters}
\end{equation}
Then $N\epsilon\delta\asymp D_0$, while $N\delta$ and
$m\epsilon\delta$ tend to infinity. Also
$m/N\asymp m^{1/3}L_m^{1/3}\gg\log m$, as required by the
coefficient lemma. Small or zero values of $R_i$ near the final
endpoint are covered by the rounded entropy estimates.
Every term in~\eqref{rva:discreteupper} and
\eqref{rva:discretelower} is
$O(m^{2/3}L_m^{2/3})$.
Finally, $J_* =\log\Lambda_C-1$ and
$\log(m!)=m\log m-m+O(\log m)$ prove
Theorem~\ref{rva:quantlaplace}.

\section{The pressure of reticulation-visible networks}
\label{rva:pressuresection}

Recall the maximal count $M_n=T_{n,n-1}$ and the defect coefficients
\[
 V_{n,q}=\operatorname{RV}_{n,3n-3-q},\qquad
 Z_n(z)=\frac1{M_n}\sum_{q=0}^{3n-3}V_{n,q}z^q.
\]
The upper bound $3n-3$ and the equality $V_{n,0}=M_n$ are the
extremal results of Chang and Fuchs~\cite{CF24}. In the preceding
combinatorial section, the relaxed maximal-skeleton model was
expressed by the kernel, for which we write $C(x;z)=\mathcal C(x;z)$
in this analytic part:
\begin{equation}
 C(x;z)=\frac{1-2zx}{(1-2zx)^{3/2}-x},\qquad z>0.
 \label{rva:RVkernel}
\end{equation}
Write $C^{[z]}(x)=C(x;z)$ and $C_0(x)=(1-x)^{-1}$.
By~\eqref{rvc:wordaverage},
\begin{equation}
 \widetilde Z_n^{(0)}(z)
       =\frac{A_{n-1}(C^{[z]})}{A_{n-1}(C_0)}.
 \label{rva:wordratio}
\end{equation}
The positive-pole description in Lemma~\ref{rvc:stieltjes}
implies~\eqref{rva:kernelassumption}, uniformly when $z$ ranges
over a compact subset of $(0,\infty)$. In particular the pole
$\rho(z)$ is simple and lies strictly before the branch point
$1/(2z)$; the separation is uniform on each such compact set.

Let $\xi(z)\in(0,\rho(z))$ solve
$\xi C_x(\xi;z)=C(\xi;z)$, and set
\begin{equation}
 \begin{split}
 I(z)&=\frac12\left\{\int_{\xi(z)}^{\rho(z)}
                     \frac{dx}{C_x(x;z)}
                     +\frac{\xi(z)}{C_x(\xi(z);z)}\right\},\\
 \Lambda(z)&=I(z)^{-1},\qquad
 \psi(z)=\log\frac{\Lambda(z)}{12}=-\log(12I(z)).
 \end{split}
 \label{rva:RVpressureformula}
\end{equation}
Here $I$ denotes the pressure integral; the large-deviation rate
function introduced below is denoted by $\mathcal I$.

\begin{theorem}[Quantitative RV pressure]
\label{rva:RVpressure}
For every compact $\mathcal K\subset(0,\infty)$,
\begin{equation}
 \log Z_n(z)=n\psi(z)
       +O_{\mathcal K}\!\left(n^{2/3}\{\log(n+2)\}^{2/3}\right),
 \qquad z\in\mathcal K.
 \label{rva:RVpressureerror}
\end{equation}
In particular,
\begin{equation}
 \lim_{n\to\infty}\frac1n\log Z_n(z)=\psi(z),
 \qquad
 \operatorname{RV}_n
 =M_n\exp\!\left\{n\psi(1)+
            O\!\left(n^{2/3}\{\log(n+2)\}^{2/3}\right)\right\}.
 \label{rva:RVtotalrelative}
\end{equation}
Equivalently, for the total number of networks,
\begin{equation}
 \log\operatorname{RV}_n
 =2n\log n+n(\log\Lambda(1)-2)
       +O\!\left(n^{2/3}\{\log(n+2)\}^{2/3}\right).
 \label{rva:RVtotalabsolute}
\end{equation}
\end{theorem}

\begin{proof}
For $C_0=(1-x)^{-1}$, the critical point is $1/2$ and
\eqref{rva:Lambda} gives $\Lambda_{C_0}=12$.
Apply Theorem~\ref{rva:quantlaplace} to numerator and denominator
of~\eqref{rva:wordratio}. This gives
$\log\widetilde Z_n^{(0)}(z)=(n-1)\psi(z)+O_{\mathcal K}(n^{2/3}
\log^{2/3}(n+2))$. The bounded ratios of
Theorem~\ref{rvc:boundedratios} change this logarithm by
$O_{\mathcal K}(1)$. Finally, $\psi$ is bounded on compact sets,
so replacing $n-1$ by $n$ is absorbed by the error.
For~\eqref{rva:RVtotalabsolute}, the exact normalization
$M_n=n!A_{n-1}(C_0)$ in~\eqref{rvc:wordaverage} and
Theorem~\ref{rva:quantlaplace} give
$\log M_n=2n\log n+n(\log12-2)+
O(n^{2/3}\log^{2/3}(n+2))$ by Stirling's formula.
\end{proof}

The bounded comparison of sectors and the quantitative word
estimate play distinct roles here: the former costs only a bounded
logarithmic term, whereas the latter supplies the displayed
sublinear error. Neither statement proves convergence of the
ratio of the full model to the maximal-skeleton model.

\subsection{Analyticity and differentiation of the pressure}

Put $u=\sqrt{1-2zx}$. Direct differentiation yields
\begin{equation}
 \frac1{C_x(x;z)}=\frac{(u^3-x)^2}{1+zu^3},
 \qquad
 C_z(x;z)=\frac{x(u^3+2x)}{(u^3-x)^2}.
 \label{rva:regularintegrand}
\end{equation}
The first expression extends analytically across the simple pole
and vanishes there. The equation for $\rho$ is
$(1-2z\rho)^{3/2}=\rho$ and has a nonzero derivative with
respect to $\rho$. The derivative with respect to $\xi$ of
$\xi C_x-C$ is $\xi C_{xx}>0$.
The implicit function theorem therefore gives analytic dependence
of both endpoints on $z>0$. It follows from
\eqref{rva:RVpressureformula} that $I$ and $\psi$ are real
analytic on $(0,\infty)$.

The following variation identity will also be used for a second
marking variable.

\begin{lemma}[Variation of the pressure integral]
\label{rva:variationlemma}
Suppose that a positive kernel $C(x;\theta)$ depends analytically
on a parameter, has a simple positive pole $\rho(\theta)$, and
has critical point $\xi C_x=C$. With
\[
 I(\theta)=\frac12\left\{
    \frac{\xi}{C_x(\xi;\theta)}+
         \int_\xi^\rho\frac{dx}{C_x(x;\theta)}\right\},
\]
one has
\begin{equation}
 \partial_\theta I
   =-\int_\xi^\rho\frac{C_\theta C_{xx}}{C_x^3}\,dx.
 \label{rva:variation}
\end{equation}
\end{lemma}

\begin{proof}
Differentiating the critical equation gives
$\xi(C_{x\theta}+C_{xx}\xi_\theta)=C_\theta$.
In the derivative of $I$, the terms $\xi_\theta/C_x(\xi)$
cancel, and the term at $\rho$ is zero. Hence
\[
 2I_\theta=-\frac{C_\theta(\xi)}{C_x(\xi)^2}
           -\int_\xi^\rho\frac{C_{x\theta}}{C_x^2}\,dx.
\]
Integrating the last term by parts proves~\eqref{rva:variation}.
The boundary term at $\xi$ cancels the displayed first term;
at the pole, $C_\theta/C_x^2\to0$, also when the pole moves
with the parameter.
\end{proof}

For $z>0$ define
\begin{equation}
 \mu(z)=z\psi'(z)
 =\Lambda(z)z\int_{\xi(z)}^{\rho(z)}
                 \frac{C_z(x;z)C_{xx}(x;z)}{C_x(x;z)^3}\,dx.
 \label{rva:muformula}
\end{equation}
The integrand is positive, so $\mu(z)>0$.
For reproducible numerical evaluation, the endpoints can also be
specified by two algebraic roots. Let $u_0,u_1\in(0,1)$ satisfy
\begin{equation}
 u_0^2+2zu_0^3=1,
 \qquad 3zu_1^5+u_1^4-zu_1^3-1=0,
 \qquad u_0<u_1.
 \label{rva:algebraicroots}
\end{equation}
Then $\rho=u_0^3$, $\xi=(1-u_1^2)/(2z)$, and, writing
$D(u)=2zu^3+u^2-1$, the same integral is
\begin{equation}
 I(z)=\frac{(1-u_1^2)^2D(u_1)}{16z^3u_1^2}
       +\frac1{8z^3}\int_{u_0}^{u_1}
                        \frac{uD(u)^2}{1+zu^3}\,du.
 \label{rva:uquadrature}
\end{equation}
This follows by $x=(1-u^2)/(2z)$ in
\eqref{rva:RVpressureformula}; the critical equation simplifies
the boundary term. The roots in~\eqref{rva:algebraicroots} are
unique in the indicated ranges, as they correspond respectively
to the unique positive pole and the unique point of mean one.

\subsection{Small and large activities}

\begin{proposition}
\label{rva:activitylimits}
The pressure extends analytically to a neighborhood of $z=0$, with
\begin{equation}
 \psi(z)=\frac{19}{4}z-\frac{265}{32}z^2+O(z^3),\qquad
 \mu(z)=\frac{19}{4}z-\frac{265}{16}z^2+O(z^3).
 \label{rva:smallactivity}
\end{equation}
At the other end,
\begin{equation}
 \psi(z)=2\log z+\log\frac{5\sqrt3}{4}+o(1),
 \qquad \mu(z)\longrightarrow2
 \quad(z\to\infty).
 \label{rva:largeactivity}
\end{equation}
Moreover, $\mu$ is strictly increasing and maps $(0,\infty)$ onto
$(0,2)$.
\end{proposition}

\begin{proof}
At $z=0$, the pole and critical-point equations are nondegenerate,
and the regular integrand in~\eqref{rva:regularintegrand} is
analytic on the integration interval. The implicit function theorem
and differentiation under the integral therefore apply in a
neighborhood of zero. Expansion gives
\begin{align*}
 \rho(z)&=1-3z+\frac{21}{2}z^2+O(z^3),&
 \xi(z)&=\frac12-\frac54z+\frac{61}{16}z^2+O(z^3),\\
 I(z)&=\frac1{12}-\frac{19}{48}z+\frac{313}{192}z^2+O(z^3),&
 \Lambda(z)&=12+57z+36z^2+O(z^3).
\end{align*}
For example, the integrand needed to verify these coefficients is
\[
 C_x^{-1}=(1-x)^2+z(-1-4x+5x^2)
                  +z^2(1+7x+x^2)+O(z^3).
\]
Taking the logarithm and differentiating yields
\eqref{rva:smallactivity}.

For large $z$, put $v=zx$. Then
\[
 C(v/z;z)\longrightarrow(1-2v)^{-1/2},\qquad
 z\xi(z)\longrightarrow\frac13,
 \qquad z\rho(z)\longrightarrow\frac12.
\]
The critical-point limit follows either from its equation or
from the mean-one equation for the limiting kernel. On the
integration interval,
\[
 0\le\frac1{\partial_v C(v/z;z)}
 =\frac{((1-2v)^{3/2}-v/z)^2}{(1-2v)^{3/2}+1/z}
 \le(1-2v)^{3/2}.
\]
Dominated convergence, also at the moving endpoints, gives
\[
 z^2I(z)\longrightarrow
 \frac12\left\{\frac13(1-2/3)^{3/2}
                    +\int_{1/3}^{1/2}(1-2v)^{3/2}\,dv\right\}
 =\frac1{15\sqrt3}.
\]
This proves the pressure expansion in~\eqref{rva:largeactivity}.

Finally, $g(t)=\psi(e^t)$ is convex, being the pointwise limit
of the convex functions $n^{-1}\log Z_n(e^t)$.
Its derivative is $g'(t)=\mu(e^t)$. The small-activity expansion
gives $g'(t)\to0$ as $t\to-\infty$, and convexity together with
$g(t)=2t+\log(5\sqrt3/4)+o(1)$ gives $g'(t)\to2$ as
$t\to\infty$. This derivative is strictly increasing: if it
took the same value at two distinct points, convexity would make
$g$ affine between them; analyticity would then make it affine
everywhere, contradicting the two limits. This proves the range
and strict monotonicity.
\end{proof}

The expansion at zero is an expansion of the limiting pressure.
It does not by itself justify an interchange between $n\to\infty$
and coefficient extraction at a defect growing with $n$.

\subsection{Typical defect and location of every mode}

Under activity $z>0$, give an RV network probability proportional
to $z^q$, and denote its defect by $Q_n$. Its distribution is
\begin{equation}
 \mathbb P_{n,z}(Q_n=q)=\frac{V_{n,q}z^q}{M_nZ_n(z)}.
 \label{rva:activitylaw}
\end{equation}

\begin{theorem}
\label{rva:concentrationmodes}
For each $z>0$ and $\epsilon>0$ there is $c_{z,\epsilon}>0$ such
that, for all sufficiently large $n$,
\begin{equation}
 \mathbb P_{n,z}(|Q_n/n-\mu(z)|>\epsilon)
       \le2e^{-c_{z,\epsilon}n},\qquad
 \frac{\mathbb E_{n,z}Q_n}{n}\longrightarrow\mu(z).
 \label{rva:defectconcentration}
\end{equation}
If $q_n^{\rm mode}(z)$ is any mode of this distribution, then
\begin{equation}
 q_n^{\rm mode}(z)/n\longrightarrow\mu(z).
 \label{rva:allmodes}
\end{equation}
\end{theorem}

\begin{proof}
The exact moment identity and the pressure theorem give
\[
 \frac1n\log\mathbb E_{n,z}e^{hQ_n}
 =\frac{\log Z_n(ze^h)-\log Z_n(z)}n
 \longrightarrow\psi(ze^h)-\psi(z).
\]
The limit has derivative $\mu(z)$ at $h=0$. Markov's inequality
with a sufficiently small fixed positive $h$ bounds the upper
tail by an exponential with a strictly negative rate; negative
$h$ gives the lower tail. Since $0\le Q_n/n\le3$, the
convergence of expectations follows. A mode has probability at
least $1/(3n-2)$, whereas all probabilities outside the indicated
linear band are exponentially small. This proves the assertion
for every possible choice of mode.
\end{proof}

Evaluation of~\eqref{rva:uquadrature} and
\eqref{rva:muformula} gives
\begin{equation}
 \psi(1)=2.151698155854262\ldots,
 \qquad \mu(1)=1.196952904766066\ldots.
 \label{rva:pressurenumerics}
\end{equation}
These are numerical evaluations of the analytic formulas, rather
than fits to finite counts. Thus a uniformly chosen RV network
has defect density $\mu(1)$ and reticulation density
$3-\mu(1)=1.803047095233934\ldots$, in the sense of
\eqref{rva:defectconcentration}. The theorem locates all modes
asymptotically but does not assert that the finite distributions
are unimodal.

\begin{figure}[tbp]
 \centering
 \includegraphics[width=\textwidth]{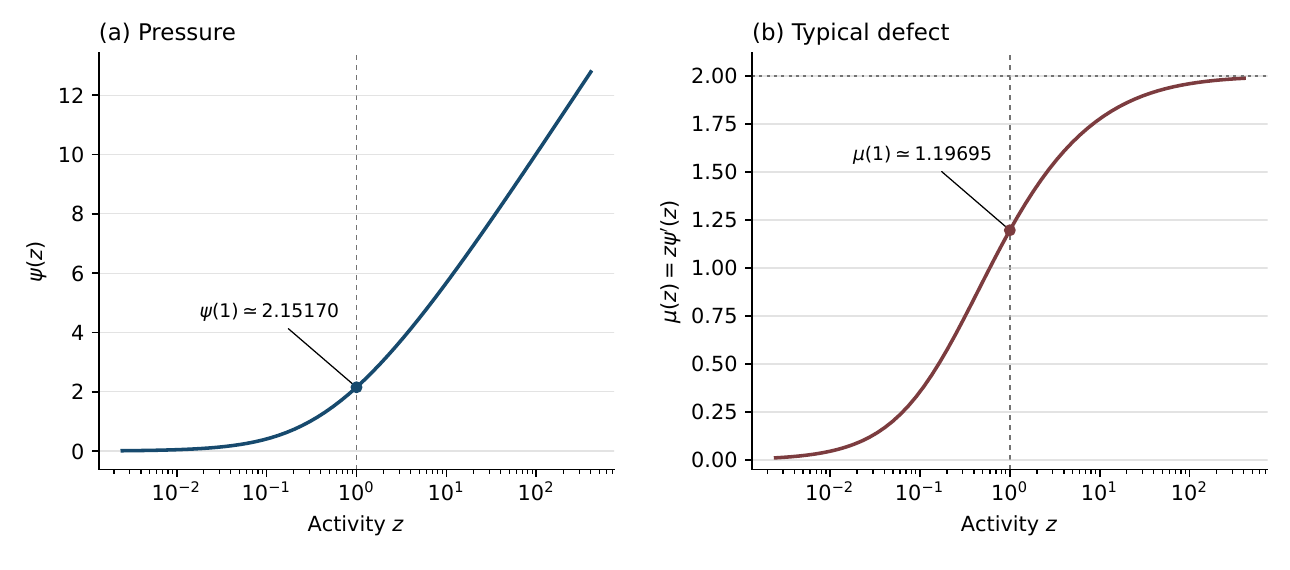}
 \caption{Pressure and typical defect density as functions of the
 activity. The marked points correspond to the uniform RV model
 $z=1$. The curves evaluate the analytic quadrature in
 \eqref{rva:uquadrature}; the horizontal dashed line is the limiting
 density $2$ as $z\to\infty$.}
 \label{fig:pressure}
\end{figure}

\section{Large deviations and counting entropy}
\label{rva:largedeviations}

In this section the RV network is chosen uniformly, so that
$\mathbb P(Q_n=q)=V_{n,q}/(M_nZ_n(1))$.
The normalized logarithmic moment-generating functions satisfy
\begin{equation}
 \ell_n(t):=\frac1n\log\mathbb E e^{tQ_n}
 =\frac{\log Z_n(e^t)-\log Z_n(1)}n
 \longrightarrow\ell(t):=\psi(e^t)-\psi(1)
 \qquad(t\in\mathbb R).
 \label{rva:limitingmgf}
\end{equation}
The derivative $\ell'(t)=\mu(e^t)$ increases strictly from zero
to two.

\begin{theorem}[Full large-deviation principle]
\label{rva:LDP}
The variables $X_n=Q_n/n$ satisfy a full large-deviation principle
with speed $n$ and good rate function
\begin{equation}
 \mathcal I(\alpha)=\sup_{t\in\mathbb R}
          \{\alpha t-\psi(e^t)+\psi(1)\}.
 \label{rva:rateLegendre}
\end{equation}
More explicitly, for $0<\alpha<2$ let $z_\alpha>0$ be the unique
solution of $\mu(z_\alpha)=\alpha$. Then
\begin{equation}
 \mathcal I(\alpha)=\alpha\log z_\alpha
                   -\psi(z_\alpha)+\psi(1),
 \label{rva:rateparametric}
\end{equation}
while
\begin{equation}
 \mathcal I(0)=\psi(1),\qquad
 \mathcal I(2)=\psi(1)-\log\frac{5\sqrt3}{4},\qquad
 \mathcal I(\alpha)=+\infty\quad(\alpha\notin[0,2]).
 \label{rva:rateendpoints}
\end{equation}
The rate is continuous on $[0,2]$, strictly convex and continuously
differentiable on $(0,2)$, with
$\mathcal I'(\alpha)=\log z_\alpha$. Its unique zero is
$\alpha_*:=\mu(1)$. For every Borel set $B$,
\begin{equation}
 -\inf_{B^\circ}\mathcal I
 \le\liminf_{n\to\infty}\frac1n\log\mathbb P(X_n\in B)
 \le\limsup_{n\to\infty}\frac1n\log\mathbb P(X_n\in B)
 \le-\inf_{\overline B}\mathcal I.
 \label{rva:LDPbounds}
\end{equation}
\end{theorem}

\begin{proof}
For $0<\alpha<2$, strict monotonicity of $\ell'$ gives the unique
maximizer $t=\log z_\alpha$ in~\eqref{rva:rateLegendre}.
At $\alpha=0$, the supremum is the limit as $t\to-\infty$;
at $\alpha=2$, it is the limit as $t\to\infty$.
Proposition~\ref{rva:activitylimits} gives exactly
\eqref{rva:rateendpoints}, and makes the supremum infinite outside
$[0,2]$. The rate is convex and lower semicontinuous, as a
supremum of continuous affine functions. Its finite endpoint
values and convexity give continuity there. Uniqueness and
continuity of the maximizing $t$ give the asserted derivative
on the interior and strict convexity. Its level sets are closed
subsets of the compact interval $[0,2]$, so it is good.

Here is a direct proof of the probability bounds. For a closed
set $F$ and a real number $A<\inf_F\mathcal I$, choose for each
$x\in F\cap[0,3]$ a finite $t_x$ such that
$xt_x-\ell(t_x)>A$. This strict inequality persists in a small
neighborhood $U_x$. By exponential Markov,
\[
 \mathbb P(X_n\in U_x)
 \le\exp\!\left\{n\ell_n(t_x)-
                  n\inf_{y\in U_x}yt_x\right\}.
\]
A finite subcover of $F\cap[0,3]$ and
\eqref{rva:limitingmgf} give a limiting upper exponent at most
$-A$. Let $A$ increase to $\inf_F\mathcal I$; if this infimum
is infinite, $A$ is arbitrary. There is no mass outside $[0,3]$.

For the lower bound, fix an interior density $\alpha$ and tilt
to activity $z=z_\alpha$. Theorem~\ref{rva:concentrationmodes}
gives
\[
 \frac1n\log\mathbb P_{n,z}(|X_n-\alpha|<\epsilon)
       \longrightarrow0\qquad(\epsilon>0).
\]
The exact change of measure yields
\begin{align*}
 \mathbb P(|X_n-\alpha|<\epsilon)
 &=\frac{Z_n(z)}{Z_n(1)}
       \mathbb E_{n,z}\!\left[
          z^{-Q_n}\mathbf1_{\{|X_n-\alpha|<\epsilon\}}\right]\\
 &\ge\frac{Z_n(z)}{Z_n(1)}
      e^{-n\alpha\log z-n\epsilon|\log z|}
      \mathbb P_{n,z}(|X_n-\alpha|<\epsilon).
\end{align*}
Thus the limiting lower exponent is at least
$-\mathcal I(\alpha)-\epsilon|\log z_\alpha|$.
For an open set containing $\alpha$, let $\epsilon\downarrow0$
through intervals contained in the open set. If the open set
contains zero or two, it also contains interior points approaching
that endpoint; continuity of $\mathcal I$ gives the same bound.
This proves~\eqref{rva:LDPbounds} without an additional
large-deviation theorem.
\end{proof}

For $0<\alpha<\alpha_*$ the left tail has exponent
$-\mathcal I(\alpha)$, and for $\alpha_*<\alpha<2$ the right
tail has that exponent. For every fixed $\epsilon>0$,
\begin{equation}
 \limsup_{n\to\infty}\frac1n
       \log\mathbb P(Q_n\ge(2+\epsilon)n)=-\infty.
 \label{rva:superexponentialtail}
\end{equation}
Indeed, Markov with $t>0$ gives the upper rate
$\ell(t)-(2+\epsilon)t$, which tends to $-\infty$ as
$t\to\infty$. This is superexponential decay at speed $n$;
it does not identify a finer scale.

At the left endpoint the equality $V_{n,0}=M_n$ also gives the
pointwise statement
\begin{equation}
 \mathbb P(Q_n=0)=Z_n(1)^{-1}
              =\exp\{-n\psi(1)+o(n)\}.
 \label{rva:zerodefectmass}
\end{equation}
At the right endpoint, however, $\mathcal I(2)$ describes
neighborhoods containing densities approaching two from below.
The large-deviation principle alone does not give an exact rate
for $\mathbb P(Q_n\ge2n)$ or $V_{n,2n}$, nor does it imply
\eqref{rva:superexponentialtail} with threshold $2n+1$.
Numerically,
\begin{equation}
 \mathcal I(0)=2.151698155854262\ldots,
 \qquad
 \mathcal I(2)=1.379248460205997\ldots.
 \label{rva:ratenumerics}
\end{equation}

\subsection{Entropy in linear windows}

Define the normalized counting entropy by
\begin{equation}
 S(\alpha)=\psi(1)-\mathcal I(\alpha)
          =\inf_{z>0}\{\psi(z)-\alpha\log z\},
 \qquad 0\le\alpha\le2.
 \label{rva:countingentropy}
\end{equation}
It is continuous on $[0,2]$, strictly concave in the interior,
and has the explicit parametrization
\begin{equation}
 \alpha=\mu(z),\qquad S(\alpha)=\psi(z)-\mu(z)\log z,
 \qquad z>0.
 \label{rva:entropycurve}
\end{equation}
Its maximum is $S(\alpha_*)=\psi(1)$, and its endpoint values are
\begin{equation}
 S(0)=0,
 \qquad S(2)=\log\frac{5\sqrt3}{4}
                  =0.772449695648265\ldots.
 \label{rva:entropyends}
\end{equation}

\begin{corollary}[Counts in linear windows]
\label{rva:windowentropy}
For every fixed $0\le a<b\le2$,
\begin{equation}
 \lim_{n\to\infty}\frac1n
    \log\left(\frac1{M_n}
             \sum_{q:\,a\le q/n\le b}V_{n,q}\right)
       =\max_{\alpha\in[a,b]}S(\alpha).
 \label{rva:windowlimit}
\end{equation}
Either endpoint of the window may be excluded without changing
the limit. For $\alpha\in[0,2]$,
\begin{equation}
 \lim_{\epsilon\downarrow0}\lim_{n\to\infty}\frac1n
    \log\left(\frac1{M_n}
             \sum_{q:\,|q/n-\alpha|<\epsilon}V_{n,q}\right)
       =S(\alpha).
 \label{rva:shrinkingwindow}
\end{equation}
\end{corollary}

\begin{proof}
The normalized sum in~\eqref{rva:windowlimit} equals
$Z_n(1)\mathbb P(X_n\in[a,b])$. Continuity of the rate implies
that its infima over $(a,b)$ and $[a,b]$ agree, including when
$a=0$ or $b=2$. Thus the upper and lower bounds in
Theorem~\ref{rva:LDP} coincide. Add the pressure limit
$n^{-1}\log Z_n(1)\to\psi(1)$. The local statement follows
by continuity as the window contracts.
\end{proof}

The maximizer in~\eqref{rva:windowlimit} is the projection of
$\alpha_*$ onto $[a,b]$. Using the exponential order of $M_n$
already obtained above gives the unnormalized form
\begin{equation}
 \log\sum_{q:\,a\le q/n\le b}V_{n,q}
 =2n\log n+
    n\left(\log12-2+\max_{\alpha\in[a,b]}S(\alpha)\right)+o(n).
 \label{rva:absolutewindow}
\end{equation}
The order of limits in~\eqref{rva:shrinkingwindow} is essential:
it does not prescribe an $n$-dependent sublinear window width.

The small-activity expansion also yields a small-density expansion:
\begin{equation}
 S(\alpha)=\alpha\log\frac{19e}{4\alpha}
                 -\frac{265}{722}\alpha^2+O(\alpha^3)
 \qquad(\alpha\downarrow0).
 \label{rva:smallentropy}
\end{equation}
To see this, write $\psi(z)=az+bz^2+O(z^3)$, where
$a=19/4$ and $b=-265/32$. Then
$z_\alpha=\alpha/a-2b\alpha^2/a^3+O(\alpha^3)$, and substitution
in~\eqref{rva:countingentropy} gives
$S(\alpha)=\alpha\log(ae/\alpha)+(b/a^2)\alpha^2+O(\alpha^3)$.
This connects the entropy to the same constant $19/4$ appearing
at small defect, while retaining the order of limits stated above.

For any integer sequence $q_n/n\to\alpha$, positivity gives
$V_{n,q_n}z^{q_n}/M_n\le Z_n(z)$ and hence
\begin{equation}
 \limsup_{n\to\infty}\frac1n\log\frac{V_{n,q_n}}{M_n}
       \le S(\alpha),\qquad 0\le\alpha\le2.
 \label{rva:individualupper}
\end{equation}
Outside $[0,2]$ this limiting upper exponent is $-\infty$.
Since a window contains at most $3n-2$ coefficients, its largest
coefficient has the same exponential rate as its sum. A diagonal
choice of windows in~\eqref{rva:shrinkingwindow} therefore gives,
for each $\alpha\in[0,2]$, some sequence $q_n/n\to\alpha$
attaining the rate $S(\alpha)$. This does not prove attainment
for a prescribed sequence such as $q_n=\lfloor\alpha n\rfloor$;
a local comparison of coefficients would be needed.

Finally, the number of reticulations is $R_n=3n-3-Q_n$.
The deterministic shift $3/n$ does not change fixed-window
large-deviation rates. Thus $R_n/n$ has good rate
$\mathcal I(3-\beta)$ on $1\le\beta\le3$, and infinite rate
outside that interval.

\begin{figure}[tbp]
 \centering
 \includegraphics[width=\textwidth]{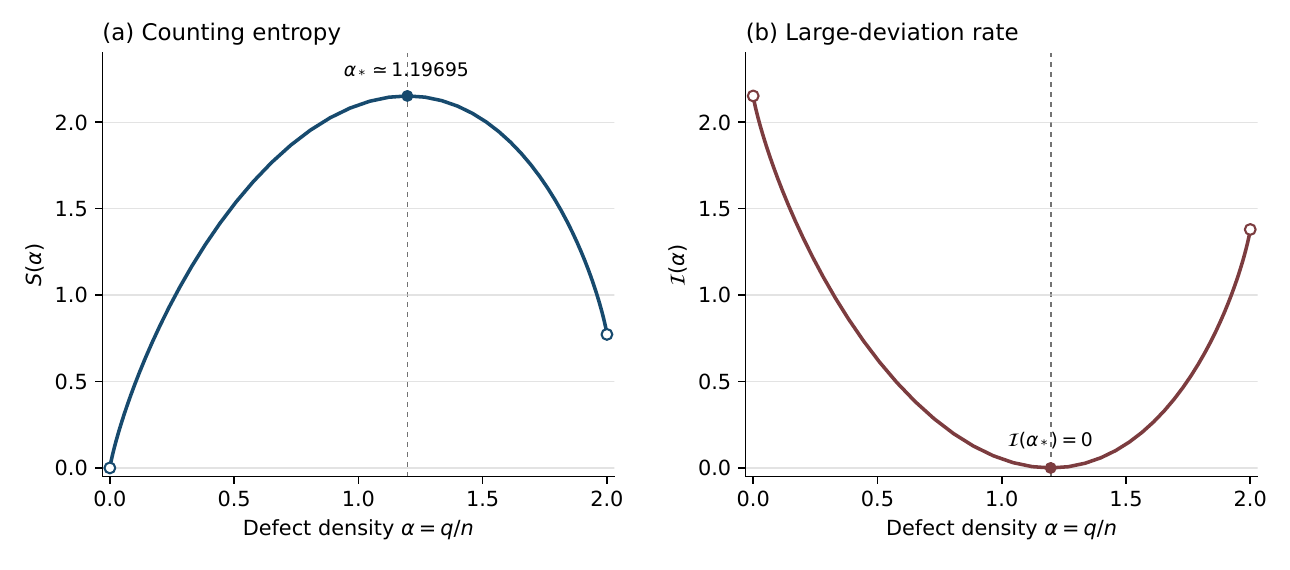}
 \caption{Counting entropy and the large-deviation rate at activity
 one. They satisfy $S(\alpha)+\mathcal I(\alpha)=\psi(1)$.
 Open circles mark the exact limiting endpoint values. The vertical
 line is the common location of the entropy maximum and rate minimum,
 $\alpha_*=\mu(1)$.}
 \label{fig:entropy}
\end{figure}

\section{A structural law for leaves and terminal reticulations}
\label{rva:structurallaw}

Let $S_n$ count the leaves of the original RV network whose parent
is a tree vertex. In the component representation, a leaf with a
tree parent in the normalized skeleton is an optional terminal
port. Marking it reconstructs a leaf with a reticulation parent;
omitting its mark reconstructs an ordinary leaf. A leaf already
below a skeleton reticulation remains below that reticulation.
The artificial edges introduced by normalization are not optional
ports. These cases exhaust the leaves, so
\begin{equation}
 S_n=s,
 \label{rva:leafinterpretation}
\end{equation}
where $s$ is the omitted-mark summand in $q=\delta+e+s$.
The equality concerns the original network. Inverting a maximal
refinement may suppress a reticulation with both parents in one
block; if its child was a leaf, this creates exactly one of the
additional optional ordinary ports, with the same interpretation.

\subsection{The bivariate kernel and its pressure}

Introduce
\begin{equation}
 Z_n(z,y)=\frac1{M_n}\sum_{N\in\mathcal{RV}_n}
                        z^{q(N)}y^{S_n(N)},\qquad z,y>0.
 \label{rva:bivariateZ}
\end{equation}
Only terminal block weights change. With
$v_t=(2t-1)!!/t!$, they become
\[
 \omega_t=t v_tz^{t-1},\qquad
 \tau_t(z,y)=(t+zy)v_tz^{t-1}.
\]
Consequently the chain kernel is
\begin{equation}
 C(x;z,y)=
 \frac{1+y\bigl((1-2zx)^{-1/2}-1\bigr)}
      {1-x(1-2zx)^{-3/2}}.
 \label{rva:bivariatekernel}
\end{equation}
Indeed, the nonterminal block series is
$W=x(1-2zx)^{-3/2}$, and the terminal block series is
$W+y((1-2zx)^{-1/2}-1)$. A nonempty chain consists of any
number of nonterminal blocks followed by a terminal block;
adding the empty chain gives~\eqref{rva:bivariatekernel}.
At $y=1$ this is~\eqref{rva:RVkernel}.

The transfer to the full model is supplied by
Proposition~\ref{rvc:bivariatecomparison}, which proves, uniformly
on positive compact sets,
\begin{equation}
 \log Z_n(z,y)-\log\widetilde Z_n^{(0)}(z,y)=O(1).
 \label{rva:bivariatecomparison}
\end{equation}
This comparison is needed in addition to the relaxed-kernel
calculation.

The bivariate coefficients are positive. Its simple positive
pole is still $\rho(z)$, independently of $y>0$: the numerator
in~\eqref{rva:bivariatekernel} is positive at that pole.
The analytic separation already used for~\eqref{rva:RVkernel}
therefore gives the uniform kernel hypotheses on compact subsets
of the positive quadrant. Let $\xi(z,y)$ solve
$\xi C_x(\xi;z,y)=C(\xi;z,y)$, and define
\begin{equation}
 I(z,y)=\frac12\left\{\frac{\xi}{C_x(\xi;z,y)}
               +\int_\xi^{\rho(z)}\frac{dx}{C_x(x;z,y)}\right\},
 \qquad \psi(z,y)=-\log(12I(z,y)).
 \label{rva:bivariatepressure}
\end{equation}
The critical point is unique, since $xC_x/C$ is a strictly
increasing mean. It is nondegenerate because $\xi C_{xx}>0$.

To check analytic regularity of the integrand explicitly, put
$u=\sqrt{1-2zx}$ and
\[
 D=u^3-x,\qquad N=u^2[y+(1-y)u],\qquad
 B=y(1+zu^3)+(1-y)u(1+zx).
\]
Then $C=N/D$ and $C_x=B/D^2$. The positive coefficient
expansion gives $B>0$ on the real integration interval;
at the pole, $B=N(1+3zu)>0$. Thus $C_x^{-1}=D^2/B$ extends
analytically and vanishes there. The implicit function theorem
and the regular integral prove that $\psi(z,y)$ is real
analytic on $(0,\infty)^2$.

The word representation, Theorem~\ref{rva:quantlaplace}, and
Proposition~\ref{rvc:bivariatecomparison} now give
\begin{equation}
 \log Z_n(z,y)=n\psi(z,y)
       +O_{\mathcal K}\!\left(n^{2/3}\{\log(n+2)\}^{2/3}\right)
 \label{rva:bivariateerror}
\end{equation}
on each positive compact set $\mathcal K$.

\subsection{Leaf densities and other structural consequences}

At $y=1$ one has
$C_y(x;z,1)=(1-\sqrt{1-2zx})C(x;z)$.
Applying Lemma~\ref{rva:variationlemma} therefore gives
\begin{equation}
 \begin{split}
 \sigma(z)&:=\left.\partial_{\log y}\psi(z,y)\right|_{y=1}\\
 &=\Lambda(z)\int_{\xi(z)}^{\rho(z)}
       (1-\sqrt{1-2zx})
              \frac{C(x;z)C_{xx}(x;z)}{C_x(x;z)^3}\,dx.
 \end{split}
 \label{rva:leafdensity}
\end{equation}
The measure with density
$\Lambda(z)CC_{xx}/C_x^3$ on this interval has total mass one
by~\eqref{rva:Lambdaalternative}. Its multiplier
$1-\sqrt{1-2zx}$ is strictly between zero and one, so
$0<\sigma(z)<1$.

\begin{theorem}[Structural densities]
\label{rva:leaflaw}
Under the RV law with any fixed activity $z>0$,
\begin{equation}
 \frac{S_n}{n}\xrightarrow{\mathbb P}\sigma(z),\qquad
 \frac{\mathbb E_{n,z}S_n}{n}\longrightarrow\sigma(z).
 \label{rva:leafLLN}
\end{equation}
For every fixed $\epsilon>0$, the probability of
$|S_n/n-\sigma(z)|>\epsilon$ decays exponentially in $n$.
Jointly, with the same type of exponential concentration,
\begin{equation}
 \left(\frac{Q_n}{n},\frac e n,\frac s n,\frac\delta n\right)
 \xrightarrow{\mathbb P}
 \bigl(\mu(z),\mu(z)-\sigma(z),\sigma(z),0\bigr).
 \label{rva:jointstructure}
\end{equation}
The corresponding expectations converge as well.
\end{theorem}

\begin{proof}
The exact identity
$\mathbb E_{n,z}e^{hS_n}=Z_n(z,e^h)/Z_n(z,1)$ and
\eqref{rva:bivariateerror} give the limiting logarithmic moment
function $\psi(z,e^h)-\psi(z,1)$, differentiable at zero with
derivative $\sigma(z)$. Positive and negative sufficiently small
$h$ in Markov's inequality give both exponential tails.
Boundedness of $S_n/n$ proves convergence of the mean.
The concentration of $Q_n/n$ was proved in
Theorem~\ref{rva:concentrationmodes}. The uniform exponential
moment of $\delta$ in Theorem~\ref{rvc:tightdelta} gives an
exponential bound for $\delta\ge\epsilon n$ and
$\mathbb E\delta=O_z(1)$. Finally,
$s=S_n$ and $e=Q_n-s-\delta$ give the joint assertion by a
union bound. All four scaled variables are bounded, so their
expectations also converge.
\end{proof}

At activity one, evaluation of the analytic integrals gives
\begin{align}
 \sigma(1)&=0.218225585802869\ldots,
 &1-\sigma(1)&=0.781774414197131\ldots,\nonumber\\
 \mu(1)-\sigma(1)&=0.978727318963198\ldots.
 \label{rva:structuralnumerics}
\end{align}
Thus the limiting fraction of leaves with a tree parent is about
$21.82\%$, and the fraction with a reticulation parent is about
$78.18\%$. These statements refer to the original RV network,
not just to its normalized skeleton.

If $R_n^{\rm int}$ counts reticulations whose child is not a
leaf, then each of the $n-s$ remaining leaves has its own
reticulation parent, because a reticulation has exactly one child.
Consequently,
\begin{equation}
 R_n^{\rm int}=(3n-3-Q_n)-(n-s),\qquad
 \frac{R_n^{\rm int}}n\xrightarrow{\mathbb P}
       2-\mu(z)+\sigma(z).
 \label{rva:internalreticulations}
\end{equation}
In the uniform model this last density is
$1.021272681036802\ldots$.

\paragraph{Scope of the estimates.}
The word and pressure errors are uniform with an explicit order,
but the implied constants have not been optimized or numerically
certified. The displayed decimal constants are evaluations of
the analytic integral formulas; they are not interval enclosures.
The results establish exponential concentration, location of
every mode, a large-deviation principle, and entropy in fixed
linear windows. They do not establish a central limit theorem,
convergence of variances divided by $n$, finite-size unimodality,
or a multiplicative prefactor for the full RV count.

% END SECTION SOURCE: rv_analysis.tex

\section{Finite counts and reproducible calculations}\label{sec:verification}
The proofs above establish the asymptotic assertions. Finite
calculations provide independent checks on the conventions and
normalizations. Table~\ref{tab:finite} compares total RV counts and
the maximal-skeleton contribution at activity one. For $n\le5$,
the component expansion and the maximal-refinement identity give
the same complete defect rows; their two implementations enumerate
different combinatorial objects. The $n=6$ row is obtained from the
maximal-refinement identity and is not claimed as an independent
graph census.

\begin{table}[tbp]
 \centering
 \small
 \begin{tabular}{@{}rrrrr@{}}
 \toprule
 $n$ & $M_n$ & $\RV_n$ & $M_n Z_n^{(0)}(1)$ & $Z_n(1)/Z_n^{(0)}(1)$\\
 \midrule
 2 & 2 & 10 & 4 & 2.500000\\
 3 & 42 & 1\,587 & 744 & 2.133065\\
 4 & 2\,544 & 805\,224 & 391\,980 & 2.054248\\
 5 & 309\,000 & 838\,617\,330 & 414\,789\,600 & 2.021790\\
 6 & 63\,068\,400 & 1\,479\,794\,012\,865 & 738\,404\,277\,210 & 2.004043\\
 \bottomrule
 \end{tabular}
 \caption{Exact finite counts and bounded sector ratios. The last
 column is rounded from exact rational numbers. These values do not
 establish convergence of the ratios.}
 \label{tab:finite}
\end{table}

The computational checks use exact integer recurrences, symbolic
Airy coefficient calculations, numerical quadrature formulas, and
figure-generation code. A rational interval certificate is given for
the bound $30$ in the sector comparison. The pressure, densities and
plotted curves use arbitrary-precision numerical quadrature and are
not presented as certified intervals. These computations evaluate
the displayed integral formulas, not parameters fitted to network counts.

\section{Further questions}
The bounded ratio $Z_n(z)/Z_n^{(0)}(z)$ and the exponential tightness
of the skeleton defect leave a concrete next problem: determine
whether these ratios and the skeleton-defect distributions converge,
and identify their limits. Such a result would require control below
the exponential scale of the weighted word model.

A second problem is a local limit theorem for the total RV defect.
The pressure identifies all limiting mode locations and the entropy
of linear windows, but a coefficient equivalent at
$q=\lfloor\alpha n\rfloor$, $0<\alpha<2$, requires local information not supplied
by the present Laplace principle. The same distinction applies to
a Gaussian fluctuation theorem and to a limiting variance. Finally,
the observed uniqueness of the finite modes suggests a unimodality
question which remains separate from their proven asymptotic
location.


\begin{thebibliography}{99}
\bibitem{CF24}
Y.-S.~Chang and M.~Fuchs,
Counting phylogenetic networks with few reticulation vertices:
galled and reticulation-visible networks,
\emph{Bulletin of Mathematical Biology} \textbf{86} (2024), 76.
\url{https://arxiv.org/abs/2401.08958}.

\bibitem{CFLLWY24}
Y.-S.~Chang, M.~Fuchs, H.~Liu, M.~Wallner and G.-R.~Yu,
Enumerative and distributional results for $d$-combining tree-child
networks, \emph{Advances in Applied Mathematics} \textbf{157} (2024),
102704.
\url{https://doi.org/10.1016/j.aam.2024.102704}.

\bibitem{FYZ21}
M.~Fuchs, G.-R.~Yu and L.~Zhang,
On the asymptotic growth of the number of tree-child networks,
\emph{European Journal of Combinatorics} \textbf{93} (2021), 103278.
\url{https://doi.org/10.1016/j.ejc.2020.103278}.

\bibitem{FY26}
M.~Fuchs and H.~Yu,
A cube-root phase transition in tree-child networks and the
enumeration threshold for galled networks,
preprint (2026), arXiv:2608.20860.
\url{https://arxiv.org/abs/2608.20860}.

\bibitem{Janson07}
S.~Janson,
Brownian excursion area, Wright's constants in graph enumeration,
and other Brownian areas,
\emph{Probability Surveys} \textbf{4} (2007), 80--145.
\url{https://doi.org/10.1214/07-PS104}.

\bibitem{Lin26}
Z.~Lin, F.~Liu, J.~Liu, J.~Liu and G.~Xin,
Proof of a conjecture on Young tableaux with walls,
preprint (2026), arXiv:2601.09551v3.
\url{https://arxiv.org/abs/2601.09551}.

\bibitem{DLMF}
F.~W.~J.~Olver, A.~B.~Olde Daalhuis, D.~W.~Lozier,
B.~I.~Schneider, R.~F.~Boisvert, C.~W.~Clark, B.~R.~Miller,
B.~V.~Saunders, H.~S.~Cohl and M.~A.~McClain, eds.,
\emph{NIST Digital Library of Mathematical Functions},
Chapter~9: Airy and related functions.
\url{https://dlmf.nist.gov/9}.

\bibitem{PB21}
M.~Pons and J.~Batle,
Combinatorial characterization of a certain class of words and a
conjectured connection with general subclasses of phylogenetic
tree-child networks,
\emph{Scientific Reports} \textbf{11} (2021), 21875.
\url{https://doi.org/10.1038/s41598-021-01166-w}.

\bibitem{YZ26}
H.~Yu and L.~Zhang,
Asymptotic counting of binary phylogenetic networks,
preprint (2026), arXiv:2605.23126v2.
\url{https://arxiv.org/abs/2605.23126}.
\end{thebibliography}
\end{document}